\documentclass[12pt,leqno]{article}
\usepackage{amssymb}
\usepackage[mathscr]{eucal}
\usepackage{amsmath,amssymb,latexsym,theorem,bbm}
\usepackage{amsmath,amssymb,latexsym,theorem}
\usepackage{color}
\usepackage{appendix}

\newcommand{\NN}{\mathbb{N}}
\newcommand{\QQ}{\mathbb{Q}}
\newcommand{\RR}{\mathbb{R}}

\newcommand{\ZZ}{\mathbb{Z}}

\newcommand{\bA}{{\boldsymbol{A}}}

\newcommand{\bB}{{\boldsymbol{B}}}

\newcommand{\tb}{\widetilde{b}}

\newcommand{\tbB}{\widetilde{\bB}}

\newcommand{\bc}{{\boldsymbol{c}}}

\newcommand{\bD}{{\boldsymbol{D}}}

\newcommand{\be}{{\boldsymbol{e}}}

\newcommand{\Bf}{{\boldsymbol{f}}}

\newcommand{\tcF}{\widetilde{\cF}}
\newcommand{\ttcF}{\widetilde{\tcF}}
\newcommand{\bF}{{\boldsymbol{F}}}

\newcommand{\tcG}{\widetilde{\cG}}
\newcommand{\ttcG}{\widetilde{\tcG}}

\newcommand{\bI}{{\boldsymbol{I}}}

\newcommand{\oM}{\overline{M}}
\newcommand{\oN}{\overline{N}}
\newcommand{\toN}{\widetilde{\oN}}
\newcommand{\hN}{\widehat{N}}
\newcommand{\tcN}{\widetilde{\cN}}
\newcommand{\ttcN}{\widetilde{\tcN}}

\newcommand{\tp}{\widetilde{p}}

\newcommand{\tPP}{\widetilde{\PP}}
\newcommand{\ttPP}{\widetilde{\tPP}}

\newcommand{\br}{{\boldsymbol{r}}}

\newcommand{\bu}{{\boldsymbol{u}}}

\newcommand{\bv}{{\boldsymbol{v}}}

\newcommand{\bx}{{\boldsymbol{x}}}
\newcommand{\bX}{{\boldsymbol{X}}}
\newcommand{\by}{{\boldsymbol{y}}}
\newcommand{\bY}{{\boldsymbol{Y}}}

\newcommand{\tbX}{\widetilde{\bX}}

\newcommand{\tbW}{\widetilde{\bW}}
\newcommand{\bz}{{\boldsymbol{z}}}

\newcommand{\bw}{{\boldsymbol{w}}}
\newcommand{\bW}{{\boldsymbol{W}}}
\newcommand{\Bbeta}{{\boldsymbol{\beta}}}
\newcommand{\tbeta}{\widetilde{\beta}}
\newcommand{\tBbeta}{\widetilde{\Bbeta}}

\newcommand{\blambda}{{\boldsymbol{\lambda}}}

\newcommand{\bpi}{{\boldsymbol{\pi}}}
\newcommand{\bxi}{{\boldsymbol{\xi}}}
\newcommand{\bmu}{{\boldsymbol{\mu}}}

\newcommand{\bbeta}{{\boldsymbol{\beta}}}

\newcommand{\tomega}{\widetilde{\omega}}
\newcommand{\ttomega}{\widetilde{\tomega}}
\newcommand{\tOmega}{\widetilde{\Omega}}
\newcommand{\ttOmega}{\widetilde{\tOmega}}

\newcommand{\bzero}{{\boldsymbol{0}}}

\newcommand{\cA}{{\mathcal A}}
\newcommand{\cB}{{\mathcal B}}

\newcommand{\cD}{{\mathcal D}}

\newcommand{\tE}{\widetilde{E}}
\newcommand{\ttE}{\widetilde{\tE}}

\newcommand{\cG}{{\mathcal G}}
\newcommand{\cF}{{\mathcal F}}
\newcommand{\cH}{{\mathcal H}}

\newcommand{\cM}{{\mathcal M}}

\newcommand{\cN}{{\mathcal N}}

\newcommand{\cR}{{\mathcal R}}

\newcommand{\tX}{\widetilde{X}}

\newcommand{\tbw}{\widetilde{\bw}}
\newcommand{\tw}{\widetilde{w}}

\newcommand{\cc}{\mathrm{c}}
\newcommand{\dd}{\mathrm{d}}
\newcommand{\ee}{\mathrm{e}}
\newcommand{\ii}{\mathrm{i}}

\newcommand{\EE}{\operatorname{\mathbb{E}}}
\newcommand{\tEE}{\operatorname{\mathbb{\tE}}}
\newcommand{\ttEE}{\operatorname{\mathbb{\ttE}}}
\newcommand{\PP}{\operatorname{\mathbb{P}}}

\newcommand{\tN}{\widetilde{N}}

\newcommand{\vare}{\varepsilon}

\renewcommand{\mid}{\,|\,}

\renewcommand{\leq}{\leqslant}
\renewcommand{\geq}{\geqslant}

\newcommand{\as}{\stackrel{{\mathrm{a.s.}}}{\longrightarrow}}

\newcommand{\bbone}{\mathbbm{1}}

\newcommand{\proofend}{\hfill\mbox{$\Box$}}

\numberwithin{equation}{section}

\theoremstyle{change} \theorembodyfont{\em}
\newtheorem{Lem}{Lemma.}[section]
\newtheorem{Thm}[Lem]{Theorem.}

\newtheorem{Def}[Lem]{Definition.}

\theorembodyfont{\rm}
\newtheorem{Rem}[Lem]{Remark.}

\begin{document}

\begin{center}
 {\bfseries\Large Stochastic differential equation with jumps for multi-type
                   continuous state and continuous time \\[2mm]
                   branching processes with immigration}
 \\[6mm]

 {\sc\large
  M\'aty\'as $\text{Barczy}^{*,\diamond}$,
  \ Zenghu $\text{Li}^{**}$,
  \ Gyula $\text{Pap}^{***}$}

\end{center}

\vskip0.2cm

\noindent
 * Faculty of Informatics, University of Debrecen,
   Pf.~12, H--4010 Debrecen, Hungary.

\noindent
 ** School of Mathematical Sciences, Beijing Normal University,
     Beijing 100875,  People's Republic of China.

\noindent
 *** Bolyai Institute, University of Szeged,
     Aradi v\'ertan\'uk tere 1, H--6720 Szeged, Hungary.

\noindent e--mails: barczy.matyas@inf.unideb.hu (M. Barczy),
                    lizh@bnu.edu.cn (Z. Li),
                    papgy@math.u-szeged.hu (G. Pap).

\noindent $\diamond$ Corresponding author.

%\vskip0.2cm

%\centerline{\sl February 17, 2014.}

\renewcommand{\thefootnote}{}
\vspace*{0.2cm}
\footnote{The research of M. Barczy and G. Pap was realized in the frames of
 T\'AMOP 4.2.4.\ A/2-11-1-2012-0001 ,,National Excellence Program --
 Elaborating and operating an inland student and researcher personal support
 system''.
The project was subsidized by the European Union and co-financed by the
 European Social Fund.
Z. Li has been partially supported by NSFC under Grant No.\ 11131003
 and 973 Program under Grant No.\ 2011CB808001.}

\vspace*{-10mm}

\begin{abstract}
A multi-type continuous state and continuous time branching process with
 immigration satisfying some moment conditions is identified as a pathwise
 unique strong solution of certain stochastic differential equation with
 jumps.
\end{abstract}

\section{Introduction}
\label{section_intro}

Continuous state and continuous time branching processes with immigration (CBI
 processes) arise as high density limits of Galton--Watson branching processes
 with immigration, see, e.g., Li \cite[Theorem 3.43]{Li} without immigration
 and Li \cite{Li2} with immigration.
A single-type continuous state and continuous time branching process (CB
 process) is a non-negative Markov process with a branching property.
This class of processes has been first introduced by Ji\v{r}ina \cite{Jir}
 both in discrete and continuous times.
As a generalization of CB processes, Kawazu and Watanabe \cite{KawWat}
 introduced the more general class of CBI processes, where immigrants may come
 from outer sources.
They defined a single-type CBI process as an \ $[0,\infty]$-valued Markov
 process with \ $\infty$ \ as a trap in terms of Laplace transforms, see
 \cite[Definition 1.1]{KawWat}.
An analytic characterization of CBI processes was also presented by giving the
 explicit form of the corresponding non-negative strongly continuous
 contraction semigroup, see \cite[Theorem 1.1']{KawWat}.
Further, limit theorems for Galton-Watson branching processes with immigration
 towards CBI processes were also investigated, see \cite[Section 2]{KawWat}.
Dawson and Li \cite[Theorems 5.1 and 5.2]{DawLi1} proved that a general
 single-type CBI process is the pathwise unique strong solution of a
 stochastic differential equation (SDE) with jumps driven by Wiener processes
 and Poisson random measures.
Watanabe \cite[Definition 1.1]{Wat} introduced two-type CB processes as
 \ $[0,\infty)^2$-valued Markov processes satisfying a branching property.
He characterized them in an analytic way by giving the explicit form of the
 infinitesimal generator of the corresponding non-negative strongly continuous
 contraction semigroup, see Watanabe \cite[Theorem 1]{Wat}.
Fittipaldi and Fontbona \cite[Theorem 2.1]{FitFon} represented a (sub)-critical
 continuous time and continuous state branching process conditioned to never be extinct
 as a pathwise unique strong solution of an appropriate SDE with jumps.
It was also shown that a two-type {\sl diffusion} CB process can be obtained
 as a pathwise unique strong solution of an SDE (without jumps), see Watanabe
 \cite[Theorem 3]{Wat}.
Recently, for a {\sl special} two-type (not necessarily diffusion) CBI process
 (with a {\sl special} immigration mechanism), an SDE with jumps (a special
 case of the SDE \eqref{SDE_X_YW_spec} given later on) has already been
 presented by Ma \cite[Theorem 2.1]{Ma} together with the existence of a
 pathwise unique \ $[0,\infty)^2$-valued strong solution of this SDE.
For a comparison of our results with those of Ma \cite{Ma},
 see Section \ref{section_special_cases}.

The aim of the present paper is to derive and study an SDE with jumps for a
 general multi-type CBI process.
Next, we give an overview of the structure of the paper by explaining some of its
 technical merits and including some sort of preview of the types of results which are proved.

In Section \ref{section_CBI} we recall some facts about CBI processes (e.g.,
 set of admissible parameters, infinitesimal generator) with special emphasis
 on their identification (under some moment conditions) as special immigration
 superprocesses.
This identification turns out to be very important since it is the starting point
 for deriving a formula for the expectation and an SDE with jumps for a general multi-type
 CBI process (see the proofs of Lemma \ref{moment_1} and Theorem \ref{CBI_SDE}).

In Section \ref{section_weak} we formulate an SDE with jumps and, under the same moment
 conditions, we prove that this SDE admits an \ $[0,\infty)^d$-valued weak
 solution which is unique in the sense of probability law among
 \ $[0,\infty)^d$-valued weak solutions.
The idea behind of deriving such an SDE goes back to a result of Li \cite[Theorem 9.18]{Li}
 that an immigration superprocess can be represented as a sum of a continuous local martingale,
 a purely discontinuous local martingale and a drift term.
In our special case, this purely discontinuous local martingale takes the form
 \ $\int_0^t \int_{[0,\infty)^d\setminus\{\bzero\}} \bz \widetilde N_0(\dd s,\dd\bz)$, $t\geq 0$,
 \ with some (not necessarily Poisson) random measure \ $N_0(\dd s,\dd\bz)$ \ on
 \ $(0,\infty)\times ([0,\infty)^d\setminus\{\bzero\})$, \ where \ $\widetilde N_0(\dd s,\dd\bz)$ \
 denotes the compensation of \ $N_0(\dd s,\dd\bz)$.
\ The next key step is that the integral  \ $\int_0^t \int_{[0,\infty)^d\setminus\{\bzero\}} \bz \widetilde N_0(\dd s,\dd\bz)$
 \ can be rewritten as an appropriate sum of integrals with respect to a Possion and compensated Poisson
 random measures, and some additional drift term, due to a representation theorem of right
 continuous martingales, see, e.g., Ikeda and Watanabe \cite[Chapter II, Definition 1.3 and Lemma 1.2]{IkeWat}.
We also prove that any \ $[0,\infty)^d$-valued weak solution of this
 SDE is a CBI process, see Theorem \ref{CBI_SDE}.
For the proof of Theorem \ref{CBI_SDE}, we need a formula for the first moment
 of a CBI process, see Lemma \ref{moment_1}.
The proof of Lemma \ref{moment_1} is based on a formula for expectation of immigration
 superprocesses, see Li \cite[Proposition 9.11]{Li}.

In Section \ref{section_strong} we prove that, under the same moment
 conditions, there is a pathwise unique \ $[0,\infty)^d$-valued strong
 solution to the SDE \eqref{SDE_X_YW_spec} and the solution is a CBI process,
 see Theorem \ref{strong_solution_old}.
For the proof, we need a comparison theorem for the SDE \eqref{SDE_X_YW_spec}
 (see, Lemma \ref{comparison_beta}), which, in particular, yields that
 pathwise uniqueness holds for the SDE \eqref{SDE_X_YW_spec} among
 \ $[0,\infty)^d$-valued weak solutions.
The ideas of the proof of Lemma \ref{comparison_beta} follow those of Theorem 3.1 of
 Ma \cite{Ma}, which are adaptations of those of Theorem 5.5 of Fu and Li \cite{FuLi}.
More precisely, we derive an upper bound for an appropriate deterministic function of the difference of two
 \ $[0,\infty)^d$-valued weak solutions of the SDE \eqref{SDE_X_YW_spec} and then apply Gronwall's inequality.

In Section \ref{section_special_cases} we specialize our SDE
 \eqref{SDE_X_YW_spec} to dimension 1 and 2, respectively, which enables us to
 compare our results with those of Dawson and Li
 \cite[Theorems 5.1 and 5.2]{DawLi1} (single-type) and Ma
\cite[Theorem 2.1]{Ma} (two-type), respectively.
Moreover, we discuss a special case of the SDE \eqref{SDE_X_YW_spec} with
 \ $\nu = 0$, \ $\mu_i = 0$, \ $i \in \{1, \ldots, d\}$, \ i.e., without
 integrals with respect to (compensated) Poisson random measures
 (corresponding to the so-called multi-factor Cox-Ingersoll-Ross process if
 \ $\bB$ \ is diagonal, see, e.g., Jagannathan et al. \cite{JagKapSun}), and
 another special case with \ $\bc = \bzero$, \ i.e., without integral with
 respect to a Wiener process.

In Appendix \ref{section_ext} we present some facts about extensions of
 probability spaces.

Finally, we mention that our work goes beyond that of Ma \cite{Ma}
 in the sense that we consider general multi-type CBI processes with arbitrary
 branching and immigration mechanisms instead of two-type CBI processes with a special
 immigration mechanism, and we carefully present some missing details
 in the proofs of Ma \cite{Ma} for the general multi-type case such as
 the application of Theorem 9.18 in Li \cite{Li} and of
 Theorem 7.4 in Chapter II in Ikeda and Watanabe \cite{IkeWat}.
Further, in a companion paper we established Yamada-Watanabe type results
 for SDEs with jumps that are needed in the proof of Theorem \ref{strong_solution_old}
 (existence of pathwise unique strong solution of the SDE \eqref{SDE_X_YW_spec}).
We point out that Ma \cite{Ma} implicitly used these results without proving or referring to them.

\section{Multi-type CBI processes}
\label{section_CBI}

Let \ $\ZZ_+$, \ $\NN$, \ $\RR$, \ $\RR_+$  \ and \ $\RR_{++}$ \ denote the set
 of non-negative integers, positive integers, real numbers, non-negative real
 numbers and positive real numbers, respectively.
For \ $x , y \in \RR$, \ we will use the notations
 \ $x \land y := \min \{x, y\} $ \ and \ $x^+:= \max \{0, x\} $.
\ By \ $\|\bx\|$ \ and \ $\|\bA\|$, \ we denote the Euclidean norm of a vector
 \ $\bx \in \RR^d$ \ and the induced matrix norm of a matrix
 \ $\bA \in \RR^{d\times d}$, \ respectively.
The natural basis in \ $\RR^d$ \ and the Borel $\sigma$-algebras on \ $\RR^d$
 \ and on \ $\RR_+^d$ \ will be denoted by \ $\be_1$, \ldots, $\be_d$, \ and
 by \ $\cB(\RR^d)$ \ and \ $\cB(\RR_+^d)$, \ respectively.
The $d$-dimensional unit matrix is denoted by \ $\bI_d$.
\ For \ $\bx = (x_i)_{i\in\{1,\ldots,d\}} \in \RR^d$ \ and
 \ $\by = (y_i)_{i\in\{1,\ldots,d\}} \in \RR^d$, \ we will use the notation
 \ $\bx \leq \by$ \ indicating that \ $x_i \leq y_i$ \ for all
 \ $i \in \{1, \ldots, d\}$.
\ By \ $C^2_\cc(\RR_+^d,\RR)$ \ we denote the set of twice continuously
 differentiable real-valued functions on \ $\RR_+^d$ \ with compact support.
Throughout this paper, we make the conventions \ $\int_a^b := \int_{(a,b]}$
 \ and \ $\int_a^\infty := \int_{(a,\infty)}$ \ for any \ $a, b \in \RR$ \ with
 \ $a < b$.

\begin{Def}\label{Def_essentially_non-negative}
A matrix \ $\bA = (a_{i,j})_{i,j\in\{1,\ldots,d\}} \in \RR^{d\times d}$ \ is called
 essentially non-negative if \ $a_{i,j} \in \RR_+$ \ whenever
 \ $i, j \in \{1,\ldots,d\}$ \ with \ $i \ne j$, \ i.e., if \ $\bA$ \ has
 non-negative off-diagonal entries.
The set of essentially non-negative \ $d \times d$ \ matrices will be denoted
 by \ $\RR^{d\times d}_{(+)}$.
\end{Def}

\begin{Def}\label{Def_admissible}
A tuple \ $(d, \bc, \Bbeta, \bB, \nu, \bmu)$ \ is called a set of admissible
 parameters if
 \renewcommand{\labelenumi}{{\rm(\roman{enumi})}}
 \begin{enumerate}
  \item
   $d \in \NN$,
  \item
   $\bc = (c_i)_{i\in\{1,\ldots,d\}} \in \RR_+^d$,
  \item
   $\Bbeta = (\beta_i)_{i\in\{1,\ldots,d\}} \in \RR_+^d$,
  \item
   $\bB = (b_{i,j})_{i,j\in\{1,\ldots,d\}} \in \RR^{d \times d}_{(+)}$,
  \item
   $\nu$ \ is a Borel measure on \ $U_d := \RR_+^d \setminus \{\bzero\}$
    \ satisfying \ $\int_{U_d} (1\wedge \|\bz\|) \, \nu(\dd \bz) < \infty$,
  \item
   $\bmu = (\mu_1, \ldots, \mu_d)$, \ where, for each
    \ $i \in \{1, \ldots, d\}$, \ $\mu_i$ \ is a Borel measure on \ $U_d$
    \ satisfying
    \begin{align}\label{help_Li_page45}
      \int_{U_d}\left[ \Vert\bz\Vert \wedge \Vert\bz\Vert^2
                     + \sum_{j \in \{1, \ldots, d\} \setminus \{i\}} z_j\;
              \right] \mu_i(\dd \bz)<\infty.
    \end{align}
 \end{enumerate}
\end{Def}

\begin{Rem}
Our Definition \ref{Def_admissible} of the set of admissible parameters is a
 special case of Definition 2.6 in Duffie et al.~\cite{DufFilSch}, which is
 suitable for all affine processes.
Namely, one should take \ $m = d$, \ $n = 0$ \ and zero killing rate in
 Definition 2.6 in Duffie et al.~\cite{DufFilSch} noting also that part (v) of
 our Definition \ref{Def_admissible} is equivalent to the corresponding one
 \ $\int_{U_d}  \sum_{i=1}^d (1 \land z_i) \, \nu(\dd \bz  ) < \infty$
 \ in Definition 2.6 in Duffie et al.~\cite{DufFilSch}.
Indeed,
 \begin{align*}
   1 \land \|\bz\|
   \leq 1 \land \left( \sum_{i=1}^d z_i \right)
   \leq \sum_{i=1}^d (1 \land z_i)
   \leq d (1 \land \|\bz\|)
 \end{align*}
 for all \ $\bz = (z_1, \ldots, z_d) \in \RR_+^d$.
\ Further, for all \ $i\in\{1,\ldots,d\}$, \ the condition \eqref{help_Li_page45} is equivalent to
 \begin{align}\label{help_Li_page45_equiv}
  \int_{U_d}
     \left[ (1 \land z_i)^2
             + \sum_{j \in \{1, \ldots, d\} \setminus \{i\}} (1 \land z_j) \right]
      \mu_i(\dd \bz)
    < \infty
   \quad \text{and} \quad
   \int_{U_d} \Vert\bz\Vert \bbone_{\{\Vert\bz\Vert\geq 1\}}\,\mu_i(\dd \bz)<\infty.
 \end{align}
Indeed, if \eqref{help_Li_page45} holds, then
 \ $\int_{U_d} \Vert\bz\Vert \bbone_{\{\Vert\bz\Vert\geq 1\}}\,\mu_i(\dd \bz)
   = \int_{U_d} (\Vert\bz\Vert \wedge \Vert\bz\Vert^2)\bbone_{\{\Vert\bz\Vert\geq 1\}}\,\mu_i(\dd \bz)<\infty$,
 \ and using that \ $z_i\leq \Vert\bz\Vert$ \ and
 \ $(1\wedge z_i)^2 = (1\wedge z_i)^2\bbone_{\{\Vert \bz\Vert\leq 1\}} + (1\wedge z_i)^2\bbone_{\{\Vert \bz\Vert > 1\}}
 \leq \Vert\bz\Vert^2 \bbone_{\{\Vert \bz\Vert\leq 1\}} + \Vert\bz\Vert\bbone_{\{\Vert \bz\Vert > 1\}}
 = \Vert\bz\Vert \wedge \Vert\bz\Vert^2$, $i\in\{1,\ldots,d\}$, \ we have \eqref{help_Li_page45_equiv}.
If \eqref{help_Li_page45_equiv} holds, then, using again \ $z_j\leq \Vert\bz\Vert$, $j\in\{1,\ldots,d\}$,
 \ we have
 \begin{align*}
  &\int_{U_d}\left[ \Vert\bz\Vert \wedge \Vert\bz\Vert^2
                   + \sum_{j \in \{1, \ldots, d\} \setminus \{i\}} z_j\;
             \right] \mu_i(\dd \bz) \\
  & = \int_{U_d}\left[ \Vert\bz\Vert^2 + \sum_{j \in \{1, \ldots, d\} \setminus \{i\}} z_j\;
             \right] \bbone_{\{ \Vert \bz\Vert< 1 \}} \mu_i(\dd \bz)
       + \int_{U_d}\left[ \Vert\bz\Vert + \sum_{j \in \{1, \ldots, d\} \setminus \{i\}} z_j\;
             \right] \bbone_{\{ \Vert \bz\Vert \geq 1 \}} \mu_i(\dd \bz) \\
  &\leq \int_{U_d}\left[ z_i^2 + 2\sum_{j \in \{1, \ldots, d\} \setminus \{i\}} z_j\;
             \right] \bbone_{\{ \Vert \bz\Vert < 1 \}} \mu_i(\dd \bz)
       + \int_{U_d} \Vert\bz\Vert \bbone_{\{\Vert\bz\Vert\geq 1\}}\,\mu_i(\dd \bz)\\
  &\quad +  \sum_{j \in \{1, \ldots, d\} \setminus \{i\}} \int_{U_d}  \Vert\bz\Vert  \bbone_{\{\Vert\bz\Vert\geq 1\}}
             \mu_i(\dd \bz) < \infty, \qquad i\in\{1,\ldots,d\},
 \end{align*}
 yielding \eqref{help_Li_page45}.
Note that, here the finiteness of the first integral in \eqref{help_Li_page45_equiv} is nothing else but condition
 (2.11) in Definition 2.6 in Duffie et al.~\cite{DufFilSch}, and the finiteness of the second integral
 in \eqref{help_Li_page45_equiv} is an additional condition that we assume compared to Duffie et al.~\cite{DufFilSch},
 its role is explained in Remark \ref{Rem_Duffie_explaining}.
\proofend
\end{Rem}

\begin{Thm}\label{CBI_exists}
Let \ $(d, \bc, \Bbeta, \bB, \nu, \bmu)$ \ be a set of admissible parameters
 in the sense of Definition \ref{Def_admissible}.
Then there exists a unique conservative transition semigroup \ $(P_t)_{t\in\RR_+}$ \ acting on
 the Banach space (endowed with the supremum norm) of real-valued bounded
 Borel-measurable functions on the state space \ $\RR_+^d$ \ such that its
 infinitesimal generator is
 \begin{align}\label{CBI_inf_gen}
  \begin{split}
   (\cA f)(\bx)
   &= \sum_{i=1}^d c_i x_i f_{i,i}''(\bx)
      + \langle \Bbeta + \bB \bx, \Bf'(\bx) \rangle
      + \int_{U_d} \bigl( f(\bx + \bz) - f(\bx) \bigr) \, \nu(\dd \bz) \\
   &\phantom{\quad}
      + \sum_{i=1}^d
         x_i
         \int_{U_d}
          \bigl( f(\bx + \bz) - f(\bx) - f'_i(\bx) (1 \land z_i) \bigr)
          \, \mu_i(\dd \bz)
 \end{split}
 \end{align}
 for \ $f \in C^2_\cc(\RR_+^d,\RR)$ \ and \ $\bx \in \RR_+^d$, \ where \ $f_i'$
\ and \ $f_{i,i}''$, \ $i \in \{1, \ldots, d\}$, \ denote the first and second
 order partial derivatives of \ $f$ \ with respect to its \ $i$-th variable,
 respectively, and \ $\Bf'(\bx) := (f_1'(\bx), \ldots, f_d'(\bx))^\top$.
\ Moreover, the Laplace transform of the transition semigroup
 \ $(P_t)_{t\in\RR_+}$ \ has a representation
 \begin{align*}
  \int_{\RR_+^d} \ee^{- \langle \blambda, \by \rangle} P_t(\bx, \dd \by)
  = \ee^{- \langle \bx, \bv(t, \blambda) \rangle - \int_0^t \psi(\bv(s, \blambda)) \, \dd s} ,
  \qquad \bx \in \RR_+^d, \quad \blambda \in \RR_+^d , \quad t \in \RR_+ ,
 \end{align*}
 where, for any \ $\blambda \in \RR_+^d$, \ the continuously differentiable function
 \ $\RR_+ \ni t \mapsto \bv(t, \blambda)
    = (v_1(t, \blambda), \ldots, v_d(t, \blambda))^\top \in \RR_+^d$
 \ is the unique locally bounded solution to the system of differential
 equations
 \begin{equation}\label{EES}
   \partial_t v_i(t, \blambda) = - \varphi_i(\bv(t, \blambda)) , \qquad
   v_i(0, \blambda) = \lambda_i , \qquad i \in \{1, \ldots, d\} ,
 \end{equation}
 with
 \[
   \varphi_i(\blambda)
   := c_i \lambda_i^2 -  \langle \bB \be_i, \blambda \rangle
      + \int_{U_d}
         \bigl( \ee^{- \langle \blambda, \bz \rangle} - 1
                + \lambda_i (1 \land z_i) \bigr)
         \, \mu_i(\dd \bz)
 \]
 for \ $\blambda \in \RR_+^d$ \ and \ $i \in \{1, \ldots, d\}$, \ and
 \[
   \psi(\blambda)
   := \langle \bbeta, \blambda \rangle
      - \int_{U_d}
         \bigl( \ee^{- \langle \blambda, \bz \rangle} - 1 \bigr)
         \, \nu(\dd \bz) , \qquad
   \blambda \in \RR_+^d .
 \]
Further, the function \ $\RR_+\times\RR_+^d\ni(t,\blambda) \mapsto \bv(t, \blambda)$
 \ is continuous.
\end{Thm}

\begin{Rem}\label{Rem_Duffie_explaining}
This theorem is a special case of Theorem 2.7 of Duffie et
 al.~\cite{DufFilSch} with \ $m = d$, \ $n = 0$ \ and zero killing rate.
The unique existence of a locally bounded solution to the system of differential equations \eqref{EES}
 is proved by Li \cite[page 45]{Li}.
Here, we point out that the moment condition given in part (vi) in Definition \ref{Def_admissible}
 (which is stronger than the one (2.11) in Definition 2.6 in Duffie et al.~\cite{DufFilSch})
 ensures that the semigroup \ $(P_t)_{t\in\RR_+}$ \ is conservative
 (we do not need the one-point compactification of \ $\RR_+^d$),
 \ see Duffie et al.~\cite[Lemma 9.2]{DufFilSch} and Li \cite[page 45]{Li}.
For the continuity of the function \ $\RR_+\times\RR_+^d\ni(t,\blambda) \mapsto \bv(t, \blambda)$,
 \ see Duffie et al.~\cite[Proposition 6.4]{DufFilSch}.
Finally, we note that the infinitesimal generator \eqref{CBI_inf_gen} can be rewritten
 in another equivalent form, see formula \eqref{CBI_inf_gen_1} in Lemma \ref{superprocess}.
\proofend
\end{Rem}

\begin{Def}\label{Def_CBI}
A conservative Markov process with state space \ $\RR_+^d$ \ and with transition semigroup
 \ $(P_t)_{t\in\RR_+}$ \ given in Theorem \ref{CBI_exists} is called a multi-type
 CBI process with parameters \ $(d, \bc, \Bbeta, \bB, \nu, \bmu)$.
\end{Def}

In what follows, we will identify a multi-type CBI process
 \ $(\bX_t)_{t\in\RR_+}$ \ with parameters \ $(d, \bc, \Bbeta, \bB, \nu, \bmu)$
 \ under a moment condition as a special immigration superprocess.
First we parametrize the family of immigration superprocesses for which
 Theorem 9.18 in Li \cite{Li} is valid.
We will use some notations of the book of Li \cite{Li}.
For a locally compact separable metric space \ $E$, \ let us introduce the
 following function spaces:
 \begin{itemize}
  \item
   $B(E)$ \ is the space of bounded real-valued Borel functions on \ $E$,
  \item
   $B(E)^+$ \ is the space of bounded non-negative real-valued Borel functions
    on \ $E$,
  \item
   $C(E)$ \ is the space of bounded continuous real-valued functions on \ $E$,
  \item
   $C(E)^+$ \ is the space of bounded continuous non-negative real-valued
    functions on \ $E$,
  \item
   $C_0(E)$ \ is the space of continuous real-valued functions on \ $E$ \ vanishing
    at infinity.
 \end{itemize}
Let \ $M(E)$ \ denote the space of finite Borel measures on \ $E$.
\ We write \ $\mu(f) := \int_E f(x) \, \mu(\dd x)$ \ for the integral of a
 function \ $f : E \to \RR$ \ with respect to a measure \ $\mu \in M(E)$ \ if
 the integral exists.

\renewcommand{\labelenumi}{{\rm(\roman{enumi})}}
\begin{Def}\label{Def_admissible_IS}
A tuple \ $\bigl( E, (R_t)_{t\in\RR_+}, c, \beta, b, B, H_1, H_2 \bigr)$ \ is
 called a set of admissible parameters if
 \begin{enumerate}
  \item
   $E$ \ is a locally compact separable metric space,
  \item
   $(R_t)_{t\in\RR_+}$ \ is the transition semigroup of a Hunt process
    \[
      \xi = \bigl( \Omega, \cG, (\cG_t)_{t\in\RR_+}, (\xi_t)_{t\in\RR_+},
                    (\theta_t)_{t\in\RR_+}, (\PP_x)_{x\in E} \bigr)
    \]
    with values in \ $E$ \ (see, e.g., Li \cite[page 314]{Li}) such that
    \ $(R_t)_{t\in\RR_+}$ \ preserves \ $C_0(E)$, \ and
    \ $\RR_+ \ni t \mapsto R_tf \in C_0(E)$ \ is continuous in the supremum
    norm for every \ $f \in C_0(E)$,
  \item
   $c \in C(E)^+$,
  \item
   $\beta \in M(E)$,
  \item
   $b \in C(E)$,
  \item
   $H_1$ \ is a finite measure on \ $M(E)^\circ := M(E) \setminus \{0\}$
    \ (where \ $0$ \ denotes the null measure) satisfying
    \ $\int_{M(E)^\circ} \kappa(1) \, H_1(\dd\kappa) < \infty$,
  \item
   $B(x, \dd y)$ \ is a bounded kernel on \ $E$ \ (i.e., from \ $E$ \ to
    \ $E$) \ and \ $H_2(x, \dd\kappa)$ \ is a \ $\sigma$-finite kernel from
    \ $E$ \ to \ $M(E)^\circ$ \ such that
    \ $E \ni x \mapsto
       \big( \kappa(1) \land \kappa(1)^2 \big) H_2(x, \dd\kappa)$
    \ is continuous with respect to the topology of weak convergence in
    \ $M(E)^\circ$, \ and the operators
    \[
      f \mapsto
      \int_{M(E)^\circ}
       \big( \kappa(f) \land \kappa(f)^2 \big) H_2(\cdot, \dd\kappa)
      \qquad \text{and} \qquad
      f \mapsto \gamma(\cdot, f)
    \]
    preserve \ $C_0(E)^+$, \ where the kernel \ $\gamma(x, \dd y)$ \ on \ $E$
    \ is defined by
    \[
      \gamma(x, \dd y)
      := B(x, \dd y) + \int_{M(E)^\circ} \kappa_x(\dd y) \, H_2(x, \dd\kappa) ,
    \]
    where \ $\kappa_x(\dd y)$ \ denotes the restriction of \ $\kappa(\dd y)$
    \ to \ $E \setminus \{x\}$, \ and by \ $\gamma(\cdot, f)$ \ we mean the
    function
    \ $E \ni x \mapsto \gamma(x, f) := \int_E f(y) \, \gamma(x, \dd y)$.
  \end{enumerate}
\end{Def}

\begin{Rem}
Note that Condition (2.25) in Li \cite{Li} readily follows from (vii) of
 Definition \ref{Def_admissible_IS}, since a function in \ $C_0(E)$ \ is
 bounded, hence
 \[
   \sup_{x\in E}
    \int_{M(E)^\circ}
     \big[ \kappa(1) \land \kappa(1)^2 \big] H_2(x, \dd\kappa)
    < \infty , \qquad
   \sup_{x\in E} \int_{M(E)^\circ} \kappa_x(1) \, H_2(x, \dd\kappa)
   \leq \sup_{x\in E} \gamma(x, 1) < \infty ,
 \]
 where we used that \ $B(x, 1) \in \RR_+$ \ for all \ $x \in E$.
\proofend
\end{Rem}

\begin{Thm}\label{IS_exists}
Let \ $\bigl( E, (R_t)_{t\in\RR_+}, c, \beta, b, B, H_1, H_2 \bigr)$ \ be a set
 of admissible parameters in the sense of Definition \ref{Def_admissible_IS}.
Then there exists a unique transition semigroup \ $(Q_t)_{t\in\RR_+}$ \ acting
 on the Banach space (endowed with the supremum norm) of real-valued bounded
 Borel-measurable functions on the state space \ $M(E)$ \ such that its
 infinitesimal generator is
 \begin{align}\label{IS_inf_gen}
  \begin{split}
   (\cA F)(\mu)
   &= \int_E c(x) F''(\mu; x) \, \mu(\dd x)
      + \int_E
         \bigl( AF'(\mu; x) + \gamma(x, F'(\mu; \cdot))
                - b(x) F'(\mu; x) \bigr)
         \mu(\dd x) \\
   &\quad
      + \int_E F'(\mu; x) \, \beta(\dd x)
      + \int_{M(E)^\circ}
         \bigl( F(\mu + \kappa) - F(\mu) \bigr) H_1(\dd\kappa) \\
   &\quad
      + \int_E
         \left( \int_{M(E)^\circ}
                 \bigl( F(\mu + \kappa) - F(\mu)
                        - \kappa(F'(\mu; \cdot))\bigr)
                 \, H_2(x, \dd\kappa) \right)
         \mu(\dd x) ,
 \end{split}
 \end{align}
 for \ $\mu \in M(E)$ \ and functions \ $F : M(E) \to \RR$ \ of the form
 \ $F(\mu) = G(\mu(f_1), \ldots, \mu(f_n))$, \ where \ $n \in \NN$,
 \ $G \in C^2(\RR^n, \RR)$, \ and \ $f_1, \ldots, f_n \in D_0(A)$, \ where
 \ $A$ \ denotes the strong generator of \ $(R_t)_{t\in\RR_+}$ \ defined by
 \[
   Af(x) := \lim_{t\downarrow0} \frac{R_tf(x) - f(x)}{t} , \qquad x \in E ,
 \]
 where the limit is taken in the supremum norm, and the domain \ $D_0(A)$ \ of
 \ $A$ \ is the totality of functions \ $f \in C_0(E)$ \ for which the above
 limit exists,
 \[
   F'(\mu; x)
   := \lim_{\vare\downarrow0}
       \frac{F(\mu + \vare \delta_x) - F(\mu)}{\vare} , \qquad
   \mu \in M(E) , \quad x \in E ,
 \]
 and \ $F''(\mu; x)$ \ is defined by the limit with \ $F(\cdot)$ \ replaced
 by \ $F'(\cdot; x)$.

Moreover, the Laplace transform of the transition semigroup
 \ $(Q_t)_{t\in\RR_+}$ \ has a representation
 \begin{align}\label{IS_SG}
  \int_{M(E)} \ee^{- \kappa(f)} Q_t(\mu, \dd\kappa)
  = \ee^{- \mu(V_tf) - \int_0^t I(V_sf) \, \dd s} ,
  \qquad \mu \in M(E), \quad f \in B(E)^+ , \quad t \in \RR_+ ,
 \end{align}
 where, for any \ $x \in E$ \ and \ $f \in B(E)^+$, \ the continuously
 differentiable function \ $\RR_+ \ni t \mapsto V_tf(x) \in \RR_+$ \ is the
 unique locally bounded solution to the integral evolution equation
 \begin{equation*}
  V_tf(x) = R_tf(x)
           - \int_0^t
              \left( \int_E \phi(y, V_sf) \, R_{t-s}(x, \dd y) \right) \dd s ,
  \qquad t \in \RR_+ ,
 \end{equation*}
 with
 \[
   \phi(x, f)
   := c(x) f(x)^2 + b(x) f(x) - \int_E f(y) \, B(x, \dd y)
      + \int_{M(E)^\circ}
         \bigl( \ee^{- \kappa(f)} - 1 + \kappa(\{x\}) f(x) \bigr)
         \, H_2(x, \dd\kappa)
 \]
 for \ $x \in E$ \ and \ $f \in B(E)^+$, \ and
 \[
   I(f)
   := \beta(f)
      + \int_{M(E)^\circ}
         \bigl( 1 - \ee^{- \kappa(f)} \bigr) \, H_1(\dd\kappa) , \qquad
   f \in B(E)^+ .
 \]
\end{Thm}

\noindent
\textbf{Proof.} \
Formula \eqref{IS_SG}, which is, in fact, formula (9.18) in Li \cite{Li},
 defines a transition semigroup of an immigration superprocess corresponding
 to the skew convolution semigroup given by (9.7) in Li \cite{Li}.
Theorem 9.18 in Li \cite{Li} yields that the infinitesimal generator of the
 immigration superprocess in question has the form given in
 \eqref{IS_inf_gen}, and the unicity of the transition semigroup.
\proofend

\begin{Def}\label{Def_IS}
A Markov process with state space \ $M(E)$ \ and with transition semigroup
 \ $(Q_t)_{t\in\RR_+}$ \ given in Theorem \ref{IS_exists} is called an
 immigration superprocess with state space \ $M(E)$ \ with parameters
 \ $\bigl( E, (R_t)_{t\in\RR_+}, c, \beta, b, B, H_1, H_2 \bigr)$.
\end{Def}

In what follows, we identify a multi-type CBI process \ $(\bX_t)_{t\in\RR_+}$
 \ with parameters \ $(d, \bc, \Bbeta, \bB, \nu, \bmu)$ \ under the moment
 condition
 \begin{equation}\label{moment_condition_1}
 \int_{U_d} \|\bz\| \bbone_{\{\|\bz\|\geq1\}} \, \nu(\dd \bz) < \infty,
 \end{equation}
 as a special immigration superprocess.

First we introduce the modified parameters
 \ $\tBbeta := (\tbeta_i)_{i\in\{1,\ldots,d\}}$,
 \ $\tbB := (\tb_{i,j})_{i,j\in\{1,\ldots,d\}}$ \ and
 \ $\bD := (d_{i,j})_{i,j\in\{1,\ldots,d\}}$ \ given by
 \begin{gather}\label{tbeta_tB}
  \tBbeta := \Bbeta + \int_{U_d} \bz \, \nu(\dd\bz) , \qquad
  \tb_{i,j} := b_{i,j} + \int_{U_d} (z_i - \delta_{i,j})^+ \, \mu_j(\dd \bz), \\
  d_{i,j} := \tb_{i,j} - \int_{U_d} z_i \bbone_{\{\|\bz\|\geq1\}} \, \mu_j(\dd\bz) ,
  \label{dij}
 \end{gather}
 with \ $\delta_{i,j} := 1$ \ if \ $i = j$, \ and \ $\delta_{i,j} := 0$ \ if
 \ $i \ne j$.
\ The moment condition \eqref{moment_condition_1} together with the fact that
 \ $\nu$ \ and \ $\bmu$ \ satisfy Definition \ref{Def_admissible} imply
 \ $\tBbeta \in \RR_+^d$, \ $\tbB \in \RR^{d \times d}_{(+)}$ \ and
 \ $\bD \in \RR^{d \times d}_{(+)}$.
\ Indeed,
 \begin{align}\label{help4}
  \int_{U_d} \|\bz\| \, \nu(\dd\bz)
  = \int_{U_d} (1 \land \|\bz\|) \bbone_{\{\|\bz\|<1\}} \, \nu(\dd\bz)
    + \int_{U_d} \|\bz\| \bbone_{\{\|\bz\|\geq1\}} \, \nu(\dd\bz)
  < \infty
 \end{align}
 by part (v) of Definition \ref{Def_admissible} and
 \eqref{moment_condition_1}.
Moreover, for all \ $i \in \{1, \ldots, d\}$,
 \begin{align}\label{help2}
  \int_{U_d} (z_i - 1)^+ \, \mu_i(\dd\bz)
  \leq \int_{U_d} z_i \bbone_{\{z_i\geq1\}} \, \mu_i(\dd\bz)
  \leq \int_{U_d} \|\bz\| \bbone_{\{\|\bz\|\geq1\}} \, \mu_i(\dd\bz)
  < \infty
 \end{align}
 by \ $z_i \leq \|\bz\|$, \ $\bz \in \RR_+^d$, \ and
 \eqref{help_Li_page45_equiv}.
Further, for all \ $i, j \in \{1, \ldots, d\}$, \ $i \ne j$,
 \begin{align}\label{help3}
  \begin{split}
   \int_{U_d} z_i \, \mu_j(\dd\bz)
   &= \int_{U_d} z_i \bbone_{\{z_i<1\}} \, \mu_j(\dd\bz)
      + \int_{U_d} z_i \bbone_{\{z_i\geq1\}} \, \mu_j(\dd\bz) \\
   &\leq \int_{U_d} (1 \land z_i) \, \mu_j(\dd\bz)
         + \int_{U_d} \|\bz\| \bbone_{\{\|\bz\|\geq1\}} \, \mu_j(\dd\bz)
    < \infty
  \end{split}
 \end{align}
 by \ $z_i \leq \|\bz\|$, \ $\bz \in \RR_+^d$, \ part (vi) of Definition
 \ref{Def_admissible} and \eqref{help_Li_page45_equiv}.
Finally, \ $d_{i,j}$ \ is well-defined for all \ $i, j \in \{1, \ldots, d\}$
 \ because of  \eqref{help_Li_page45_equiv}, and, for all
 \ $i, j \in \{1, \ldots, d\}$, \ $i \ne j$,
 \[
   d_{i,j} = b_{i,j} + \int_{U_d} z_i \, \mu_j(\dd z)
            - \int_{U_d} z_i \bbone_{\{\|\bz\|\geq1\}} \, \mu_j(\dd \bz)
          = b_{i,j} + \int_{U_d} z_i \bbone_{\{\|\bz\|<1\}} \, \mu_j(\dd \bz)
          \in \RR_+ .
 \]
Note also that for all \ $j \in \{1, \ldots, d\}$,
 \begin{equation}\label{z2}
  \int_{U_d} \|\bz\|^2 \bbone_{\{\|\bz\|<1\}} \, \mu_j(\dd \bz)
  \leq \int_{U_d}
        \biggl( z_j^2 + \sum_{k\in\{1,\ldots,d\}\setminus\{j\}} z_k \biggr)
        \bbone_{\{\|\bz\|<1\}} \, \mu_j(\dd \bz)
  < \infty
 \end{equation}
 by \ $z_i \leq \|\bz\|$, \ $\bz \in \RR_+^d$, \ part (vi) of
 Definition \ref{Def_admissible} and \eqref{help_Li_page45_equiv}.

For the discrete metric space \ $E := \{1, \ldots, d\}$, \ we have the
 following identifications:
 \begin{itemize}
  \item
   $B(E)$, \ $C(E)$ \ and \ $C_0(E)$ \ can be identified with \ $\RR^d$,
    \ since a function \ $f : E \to \RR$ \ can be identified with the vector
    \ $(f(1), \ldots, f(d))^\top \in \RR^d$,
  \item
   $B(E)^+$ \ and \ $C(E)^+$ \ can be identified with \ $\RR_+^d$,
  \item
   $M(E)$ \ can be identified with \ $\RR_+^d$, \ since a finite Borel
    measure \ $\mu$ \ on \ $E$ \ can be identified with the vector
    \ $(\mu(\{1\}), \ldots, \mu(\{d\}))^\top \in \RR_+^d$,
  \item
   for \ $\mu \in M(E)$ \ and \ $f \in B(E)$, \ the integral
    \ $\mu(f) = \int_E f(x) \, \mu(\dd x) = \sum_{i=1}^d f(i) \mu(\{i\})$ \ can be
    identified with the usual Euclidean inner product
    \ $\langle \mu, f \rangle$ \ in \ $\RR^d$,
  \item
   $M(E)^\circ$ \ can be identified with \ $U_d$.
 \end{itemize}
If \ $(\Omega, \cF, \PP)$ \ is a probability space, then, by \ $\PP$-null sets
 from a sub $\sigma$-algebra \ $\cH\subset \cF$, \ we mean the elements of the
 set
 \[
   \{ A \subset \Omega \,
      : \, \exists B \in \cH \;
        \text{\ such that \ $A \subset B$ \ and \ $\PP(B) = 0$} \} .
 \]
A filtered probability space \ $(\Omega, \cF, (\cF_t)_{t\in\RR_+}, \PP)$ \ is
 said to satisfy the usual hypotheses if \ $(\cF_t)_{t\in\RR_+}$ \ is right
 continuous and \ $\cF_0$ \ contains all the \ $\PP$-null sets in \ $\cF$.

\begin{Lem}\label{superprocess}
Let \ $(d, \bc, \Bbeta, \bB, \nu, \bmu)$ \ be a set of admissible parameters
 in the sense of Definition \ref{Def_admissible} satisfying the moment condition
 \eqref{moment_condition_1}.
Then \ $\bigl( E, (R_t)_{t\in\RR_+}, c, \beta, b, B, H_1, H_2 \bigr)$ \ is a set
 of admissible parameters in the sense of Definition \ref{Def_admissible_IS},
 where
 \begin{enumerate}
  \item
   $E := \{1, \ldots, d\}$ \ with the discrete metric,
  \item
   $(R_t)_{t\in\RR_+}$ \ is the transition semigroup given by \ $R_tf := f$,
    \ $f \in B(E)$, \ $t \in \RR_+$,
  \item
   $c \in B(E)^+$ \ is given by \ $c(i) := c_i$, \ $i \in E$,
  \item
   $\beta \in M(E)$ \ is given by \ $\beta(\{i\}) := \beta_i$, \ $i \in E$,
  \item
   $b \in B(E)$, \ is given by \ $b(i) := - \tb_{i,i}$, \ $i \in E$,
  \item
   $B(x, \dd y)$ \ is the kernel on \ $E$ \ given by \ $B(i, \{i\}) := 0$
    \ for \ $i \in \{1, \ldots, d\}$ \ and \ $B(i, \{j\}) := b_{j,i}$ \ for
    \ $i, j \in \{1, \ldots, d\}$ \ with \ $i \ne j$,
  \item
   $H_1$ \ is the measure on \ $M(E)^\circ$ \ identified with the measure
    \ $\nu$ \ on \ $U_d$,
  \item
   $H_2(x, \dd\kappa)$ \ is the kernel from \ $E$ \ to \ $M(E)^\circ$ \ such
    that the measure \ $H_2(i, \cdot)$ \ on \ $M(E)^\circ$ \ is identified with
    the measure \ $\mu_i$ \ on \ $U_d$ \ for each \ $i \in \{1, \ldots, d\}$.
  \end{enumerate}
If \ $(\Omega, \cF, (\cF_t)_{t\in\RR_+}, \PP)$ \ is a filtered probability space
 satisfying the usual hypotheses and \ $(Y_t)_{t\in\RR_+}$ \ is a c\`adl\`ag
 immigration superprocess with parameters
 \ $\bigl( E, (R_t)_{t\in\RR_+}, c, \beta, b, B, H_1, H_2 \bigr)$ \ satisfying
 \ $\EE(Y_0(1)) < \infty$ \ and adapted to \ $(\cF_t)_{t\in\RR_+}$, \ then
 \ $\bX_t := (Y_t(\{1\}), \ldots, Y_t(\{d\}))^\top$, \ $t \in \RR_+$, \ is a
 multi-type CBI process with parameters \ $(d, \bc, \Bbeta, \bB, \nu, \bmu)$
 \ satisfying \ $\EE(\|\bX_0\|) < \infty$.
\ The infinitesimal generator \eqref{CBI_inf_gen} of \ $(\bX_t)_{t\in\RR_+}$
 \ can also be written in the form
 \begin{align}\label{CBI_inf_gen_1}
  \begin{split}
   (\cA_{\bX} f)(\bx)
   &= \sum_{i=1}^d c_i x_i f_{i,i}''(x)
      + \sum_{i=1}^d
         x_i
         \int_{U_d}
          \bigl(f(\bx + \bz) - f(\bx) - \langle \bz, \Bf'(\bx) \rangle\bigr)
          \, \mu_i(\dd \bz) \\
   &\quad
      + \langle \Bbeta + \tbB \bx, \Bf'(\bx) \rangle
      + \int_{U_d} \bigl(f(\bx + \bz) - f(\bx)\bigr) \, \nu(\dd \bz)
  \end{split}
 \end{align}
 for \ $f \in C^2_\cc(\RR_+^d,\RR)$ \ and \ $\bx \in \RR_+^d$.
\end{Lem}

\noindent
\textbf{Proof.} \
The discrete metric space \ $\{1, \ldots, d\}$ \ is trivially a locally
 compact separable metric space.
Clearly, \ $R_tf := f$, \ $f \in B(E)$, \ $t \in \RR_+$, \ is the transition
 semigroup of the Hunt process
 \[
   \xi = \bigl( \Omega, \cG, (\cG_t)_{t\in\RR_+}, (\xi_t)_{t\in\RR_+},
                (\theta_t)_{t\in\RR_+}, (\PP_x)_{x\in E} \bigr)
 \]
 with \ $\Omega = \{1, \ldots, d\}$, \ $\cG = \cG_t = 2^\Omega$,
 \ $t \in \RR_+$, \ $\xi_t(\omega) = \theta_t(\omega) = \omega$,
 \ $\omega \in \Omega$, \ $t \in \RR_+$, \ $\PP_x = \delta_x$,
 \ $x \in \{1, \ldots, d\}$.
\ Moreover, \ $(R_t)_{t\in\RR_+}$ \ trivially satisfies (ii) of Definition
 \ref{Def_admissible_IS}, and (iii), (iv) and (v) of Definition
 \ref{Def_admissible_IS} trivially hold.
Further (vi) of Definition \ref{Def_admissible_IS} also holds, since
 \ $\int_{U_d} \left(\sum_{i=1}^d z_i\right) \nu(\dd\bz) < \infty$ \ follows from
 \eqref{help4} by \ $z_i \leq \|\bz\|$, \ $\bz \in \RR_+^d$,
 \ $i \in \{1, \ldots, d\}$.
\ The kernel \ $B(x, \dd y)$ \ on \ $E$ \ is bounded, since
 \ $\sup_{x\in E} B(x, E)
    = \max_{i\in\{1,\ldots,d\}} \sum_{j\in\{1,\ldots,d\}\setminus\{i\}} b_{j,i}
    < \infty$.
\ On the dicrete metric space \ $\{1, \ldots, d\}$ \ every function is
 continuous, hence
 \ $E \ni x \mapsto \big( \kappa(1) \land \kappa(1)^2 \big) H_2(x, \dd\kappa)$
 \ is continuous with respect to the topology of weak convergence in
 \ $M(E)^\circ$.
\ In order to show that the operator
 \[
   f \mapsto
   \int_{M(E)^\circ}
    \big( \kappa(f) \land \kappa(f)^2 \big) H_2(\cdot, \dd\kappa)
 \]
 preserve \ $C_0(E)^+$, \ it suffices to observe that for each
 \ $\blambda \in \RR_+^d$ \ and \ $i \in \{1, \ldots, d\}$, \ we have
 \[
   \int_{U_d}
    \left( \langle \blambda, \bz \rangle \land
           \langle \blambda, \bz \rangle^2 \right)
    \mu_i(\dd\bz) \in \RR_+ ,
 \]
 which follows from the estimate
 \begin{align*}
  \int_{U_d}
   &\left( \langle \blambda, \bz \rangle \land
           \langle \blambda, \bz \rangle^2 \right)
    \mu_i(\dd\bz)
  \leq \int_{U_d}
        \bigl[ (\|\blambda\| \|\bz\|) \land
               (\|\blambda\| \|\bz\|)^2 \bigr]
        \, \mu_i(\dd\bz)
  \leq c_\blambda
       \int_{U_d} (\|\bz\| \land \|\bz\|^2) \, \mu_i(\dd\bz) \\
   &= c_\blambda
      \int_{U_d}
       (\|\bz\| \land \|\bz\|^2) \bbone_{\{\|\bz\|\leq1\}}
       \, \mu_i(\dd\bz)
      + c_\blambda
        \int_{U_d}
         (\|\bz\| \land \|\bz\|^2) \bbone_{\{\|\bz\|>1\}}
         \, \mu_i(\dd\bz) \\
   &= c_\blambda \int_{U_d} \|\bz\|^2 \bbone_{\{\|\bz\|\leq1\}} \, \mu_i(\dd\bz)
      + c_\blambda \int_{U_d} \|\bz\| \bbone_{\{\|\bz\|>1\}} \, \mu_i(\dd\bz)
    < \infty
  \end{align*}
 with \ $c_\blambda := \max\{\|\blambda\|, \|\blambda\|^2\}$ \ by \eqref{z2}
 and \eqref{help_Li_page45_equiv}.
In order to show that the operator \ $f \mapsto \gamma(\cdot, f)$ \ preserves
 \ $C_0(E)^+$, \ it suffices to observe that for each
 \ $\blambda = (\lambda_1, \ldots, \lambda_d)^\top \in \RR_+^d$ \ and
 \ $i \in \{1, \ldots, d\}$, \ we have
 \begin{align*}
  \sum_{j=1}^d \lambda_j B(i, \{j\})
  + \sum_{j\in\{1,\ldots,d\}\setminus\{i\}} \lambda_j \int_{U_d} z_j \, \mu_j(\dd\bz)
  \in \RR_+ ,
 \end{align*}
 which follows from \eqref{help3}.
Consequently, \ $\bigl( E, (R_t)_{t\in\RR_+}, c, \beta, b, B, H_1, H_2 \bigr)$
 \ is a set of admissible parameters in the sense of Definition
 \ref{Def_admissible_IS}.

By Theorem \ref{IS_exists}, we have
 \[
   \EE(\ee^{-Y_t(f)} \mid Y_0 = \mu)
   = \int_{M(E)^\circ} \ee^{-\kappa(f)} \, Q_t(\mu, \dd\kappa)
   = \ee^{- \mu(V_tf) - \int_0^t I(V_sf) \, \dd s}
 \]
 for \ $\mu \in M(E)$, \ $f \in B(E)^+$ \ and \ $t \in \RR_+$, \ hence we
 obtain
 \[
   \EE(\ee^{-\langle\blambda,\bX_t\rangle} \mid \bX_0 = \bx)
   = \ee^{-\langle\bx,\bv(t,\blambda)\rangle-\int_0^t\psi(\bv(s,\blambda))\,\dd s} ,
   \qquad \bx, \blambda \in \RR_+^d, \quad t \in \RR_+ ,
 \]
 where, for any \ $i \in \{1, \ldots, d\}$ \ and \ $\blambda \in \RR_+^d$,
 \ the function
 \ $\RR_+ \ni t \mapsto
    \bv(t, \blambda) = (v_1(t, \blambda), \ldots, v_d(t, \blambda))$
 \ is the unique locally bounded solution to the integral evolution equation
 \[
   v_i(t, \blambda)
   = \lambda_i - \int_0^t \varphi_i(\bv(s, \blambda)) \, \dd s , \qquad
   t \in \RR_+ , \quad i \in \{1, \ldots, d\} ,
 \]
 with
 \[
   \varphi_i(\blambda)
   := c_i \lambda_i^2 - \tb_{i,i} \lambda_i
      - \sum_{j\in\{1,\ldots,d\}\setminus\{i\}} \lambda_j b_{j,i}
      + \int_{U_d}
         (\ee^{-\langle\blambda,\bz\rangle} - 1 + \lambda_i z_i) \, \mu_i(\dd\bz)
 \]
 for \ $\blambda \in \RR_+^d$ \ and \ $i \in \{1, \ldots, d\}$, \ and
 \[
   \psi(\blambda)
   := \langle \Bbeta, \blambda \rangle
      + \int_{U_d} (1 - \ee^{-\langle\blambda,\bz\rangle}) \, \nu(\dd\bz) , \qquad
   \blambda \in \RR_+^d .
 \]
We have
 \begin{equation}\label{phi_i}
  \varphi_i(\blambda)
  = c_i \lambda_i^2 - \langle \tbB \be_i, \blambda \rangle
    + \int_{U_d}
       (\ee^{-\langle\blambda,\bz\rangle} - 1 + \langle \blambda, \bz \rangle)
       \, \mu_i(\dd\bz) ,
 \end{equation}
 since, by \eqref{help3},
 \begin{align*}
  &\varphi_i(\blambda) - c_i \lambda_i^2 + \langle \tbB \be_i, \blambda \rangle
   - \int_{U_d}
      (\ee^{-\langle\blambda,\bz\rangle} - 1 + \langle \blambda, \bz \rangle)
      \, \mu_i(\dd\bz) \\
  &= - \tb_{i,i} \lambda_i - \sum_{j\in\{1,\ldots,d\}\setminus\{i\}} \lambda_j b_{j,i}
     + \sum_{j=1}^d \lambda_j \tb_{j,i}
     - \sum_{j\in\{1,\ldots,d\}\setminus\{i\}} \lambda_j \int_{U_d} z_j \, \mu_i(\dd\bz)
     = 0 .
 \end{align*}
Moreover, we can write the functions \ $\varphi_i$, \ $i \in \{1, \ldots, d\}$,
 \ in the form
 \[
   \varphi_i(\blambda)
   = c_i \lambda_i^2 -  \langle \bB \be_i, \blambda \rangle
     + \int_{U_d}
        \bigl( \ee^{- \langle \blambda, \bz \rangle} - 1
               + \lambda_i (1 \land z_i) \bigr)
        \, \mu_i(\dd \bz)
 \]
 for \ $\blambda = (\lambda_1, \ldots, \lambda_d)^\top \in \RR_+^d$ \ and
 \ $i \in \{1, \ldots, d\}$.
\ Indeed, by \eqref{help2} and \eqref{help3},
 \begin{align*}
  \varphi_i(\blambda) - c_i \lambda_i^2 &+ \langle \bB \be_i, \blambda \rangle
  - \int_{U_d}
     \bigl( \ee^{- \langle \blambda, \bz \rangle} - 1
            + \lambda_i (1 \land z_i) \bigr)
     \, \mu_i(\dd \bz) \\
  &= \langle (\bB - \tbB) \be_i, \blambda \rangle
     - \int_{U_d}
        \bigl( \lambda_i(1 \land z_i) - \langle \blambda, \bz \rangle \bigr)
        \, \mu_i(\dd \bz) \\
  &= - \lambda_i \int_{U_d} (z_i-1)^+ \, \mu_i(\dd \bz)
     - \int_{U_d}
        \bigl( \lambda_i(1 \land z_i) - \lambda_i z_i \bigr)
        \, \mu_i(\dd \bz)
   = 0 .
 \end{align*}
By Theorem \ref{CBI_exists}, \ $(\bX_t)_{t\in\RR_+}$ \ is a multi-type CBI with
 parameters \ $(d, \bc, \Bbeta, \bB, \nu, \bmu)$
 \ satisfying \ $\EE(\|\bX_0\|) < \infty$.

Finally, \eqref{CBI_inf_gen_1} follows from
 \begin{align*}
  &(\cA_{\bX} f)(\bx) - \sum_{i=1}^d c_i x_i f_{i,i}''(\bx)
    - \sum_{i=1}^d
       x_i
       \int_{U_d}
        \bigl(f(\bx + \bz) - f(\bx) - \langle \bz, \Bf'(\bx) \rangle\bigr)
        \, \mu_i(\dd \bz) \\
  &- \langle \Bbeta + \tbB \bx, \Bf'(\bx) \rangle
   - \int_{U_d}
      \bigl(f(\bx + \bz) - f(\bx)\bigr) \, \nu(\dd \bz) \\
 &= \sum_{i=1}^d
      x_i
      \int_{U_d}
       \bigl( \langle \bz, \Bf'(\bx) \rangle - f_i'(\bx) (1 \land z_i) \bigr)
       \, \mu_i(\dd \bz)
     - \langle (\tbB - \bB) \bx, \Bf'(\bx) \rangle  \\
  &= \sum_{i=1}^d
      x_i \int_{U_d}
           \biggl( f_i'(\bx) ( z_i - (1 \land z_i) )
                   + \sum_{j\in\{1,\ldots,d\}\setminus\{i\}} z_j f_j'(\bx) \biggr)
           \mu_i(\dd \bz) \\
  &\phantom{=\;}
     - \sum_{i=1}^d
        \sum_{j=1}^d
         x_j f_i'(\bx) \int_{U_d} (z_i-\delta_{i,j})^+ \, \mu_j(\dd \bz)
   = 0 .
 \end{align*}
 using \eqref{help4}, \eqref{help2} and \eqref{help3}.
\proofend

\section{Multi-type CBI process as a weak solution of an SDE}
\label{section_weak}

Let \ $\cR := \bigcup_{j=0}^d \cR_j$, \ where
 \ $\cR_j$, \ $j \in \{0, 1, \ldots, d\}$, \ are disjoint sets given by
 \[
   \cR_0 := U_d \times \{ (\bzero, 0) \}^d
         \subset \RR_+^d \times (\RR_+^d \times \RR_+)^d ,
 \]
 and
 \[
   \cR_j := \{\bzero\} \times H_{j,1} \times \cdots \times H_{j,d}
         \subset \RR_+^d \times (\RR_+^d \times \RR_+)^d , \qquad
   j \in \{1, \ldots, d\} ,
 \]
 where
 \[
   H_{j,i} := \begin{cases}
              U_d \times U_1 & \text{if \ $i = j$,} \\
              \{ (\bzero, 0) \} & \text{if \ $i \ne j$.}
             \end{cases}
 \]
(Recall that \ $U_1 = \RR_{++}$.)
Let \ $m$ \ be the uniquely defined measure on
 \ $V := \RR_+^d \times (\RR_+^d \times \RR_+)^d$ \ such that
 \ $m(V \setminus \cR) = 0$ \ and its restrictions on \ $\cR_j$,
 \ $j \in \{0, 1, \ldots, d\}$, \ are
 \begin{equation}\label{m}
  m|_{\cR_0}(\dd\br) = \nu(\dd\br) , \qquad
  m|_{\cR_j}(\dd\bz, \dd u) = \mu_j(\dd\bz) \, \dd u ,
  \quad j \in \{1, \ldots, d\} ,
 \end{equation}
 where we identify \ $\cR_0$ \ with \ $U_d$ \ and \ $\cR_1$, \ldots, $\cR_d$
 \ with \ $U_d \times U_1$ \ in a natural way.
Using again this identification, let \ $h : \RR^d \times V \to \RR_+^d$ \ be
 defined by
 \[
   h(\bx, \br)
   := \begin{cases}
       \br, & \text{if \ $\bx \in \RR_+^d$, \ $\br \in \cR_0$,} \\
       \bz \bbone_{\{u \leq x_j\}} ,
        & \text{if \ $\bx = (x_1, \ldots, x_d)^\top \in \RR_+^d$,
                \ $\br = (\bz, u) \in \cR_j$, \ $j \in \{1, \ldots, d\}$,} \\
       \bzero , & \text{otherwise.}
      \end{cases}
 \]
Consider the decomposition \ $\cR = V_0 \cup V_1$, \ where
 \ $V_0 := \bigcup_{j=1}^d \cR_{j,0}$ \ and
 \ $V_1 := \cR_0 \cup \bigl( \bigcup_{j=1}^d \cR_{j,1} \bigr)$ \ with
 \ $\cR_{j,k} := \{\bzero\} \times H_{j,1,k} \times \cdots \times H_{j,d,k}$,
 \ $j \in \{1, \ldots, d\}$, \ $k \in \{0, 1\}$, \ and
 \[
   H_{j,i,k} := \begin{cases}
                U_{d,k} \times U_1 & \text{if \ $i = j$,} \\
                \{ (\bzero, 0) \} & \text{if \ $i \ne j$,}
               \end{cases} \qquad
   U_{d,k} := \begin{cases}
              \{ \bz \in U_d : \|\bz\| < 1 \} & \text{if \ $k = 0$,} \\
              \{ \bz \in U_d : \|\bz\| \geq 1 \} & \text{if \ $k = 1$.}
             \end{cases}
 \]
Then the sets \ $V_0$ \ and \ $V_1$ \ are disjoint, and the function \ $h$
 \ can be decomposed in the form \ $h = f + g$ \ with
 \[
   f(\bx, \br) := h(\bx, \br) \bbone_{V_0}(\br) , \qquad
   g(\bx, \br) := h(\bx, \br) \bbone_{V_1}(\br) , \qquad
   (\bx, \br) \in \RR^d \times V .
 \]
Let \ $(d, \bc, \Bbeta, \bB, \nu, \bmu)$ \ be a set of admissible parameters
 in the sense of Definition \ref{Def_admissible} such that the moment
 condition \eqref{moment_condition_1} holds.
Let us consider the \ $d$-dimensional SDE
 \begin{equation}\label{SDE_X_YW_spec}
  \begin{aligned}
   \bX_t
   &= \bX_0 + \int_0^t b(\bX_s) \, \dd s
      + \int_0^t \sigma(\bX_s) \, \dd \bW_s \\
   &\quad
      + \int_0^t \int_{V_0} f(\bX_{s-}, \br) \, \tN(\dd s, \dd \br)
      + \int_0^t \int_{V_1} g(\bX_{s-}, \br) \, N(\dd s, \dd \br) , \qquad
   t \in \RR_+ ,
  \end{aligned}
 \end{equation}
 where the functions \ $b : \RR^d \to \RR^d$ \ and
 \ $\sigma : \RR^d \to \RR^{d\times d}$ \ are defined by
 \begin{align*}
  b(\bx) := \Bbeta + \bD \bx  , \qquad
  \sigma(\bx) := \sum_{i=1}^d \sqrt{2 c_i x_i^+} \, \be_i \be_i^\top , \qquad
  \bx \in \RR^d ,
 \end{align*}
 $\bD$ \ is defined in \eqref{dij}, \ $(\bW_t)_{t\in\RR_+}$ \ is a
 \ $d$-dimensional standard Brownian motion, \ $N(\dd s, \dd\br)$ \ is a
 Poisson random measure on \ $\RR_{++} \times V$ \ with intensity measure
 \ $\dd s \, m(\dd\br)$, \ and
 \ $\tN(\dd s, \dd\br) := N(\dd s, \dd\br) - \dd s \, m(\dd\br)$.
\ For a short review on point measures and point processes needed for this
 paper, see, e.g., Barczy et al. \cite[Section 2]{BarLiPap}.

\begin{Def}\label{Def_weak_solution}
Let \ $n$ \ be a probability measure on \ $(\RR_+^d, \cB(\RR_+^d))$.
\ An \ $\RR_+^d$-valued weak solution of the SDE \eqref{SDE_X_YW_spec} with
 initial distribution \ $n$ \ is a tuple
 \ $\bigl( \Omega, \cF, (\cF_t)_{t\in\RR_+}, \PP, \bW, p, \bX \bigr)$, \ where
 \begin{enumerate}
  \item[\textup{(D1)}]
   $(\Omega, \cF, (\cF_t)_{t\in\RR_+}, \PP)$ \ is a filtered probability space
    satisfying the usual hypotheses;
  \item[\textup{(D2)}]
   $(\bW_t)_{t\in\RR_+}$ \ is a \ $d$-dimensional standard
    \ $(\cF_t)_{t\in\RR_+}$-Brownian motion;
  \item[\textup{(D3)}]
   $p$ \ is a stationary \ $(\cF_t)_{t\in\RR_+}$-Poisson point process on \ $V$
    \ with characteristic measure \ $m$ \ given in \eqref{m};
  \item[\textup{(D4)}]
   $(\bX_t)_{t\in\RR_+}$ \ is an \ $\RR_+^d$-valued \ $(\cF_t)_{t\in\RR_+}$-adapted
    c\`adl\`ag process such that
     \begin{itemize}
      \item[\textup{(a)}]
       the distribution of \ $\bX_0$ \ is \ $n$, \
      \item[\textup{(b)}]
       $\PP\bigl(\int_0^t
                  \big(\|b(\bX_s)\| + \|\sigma(\bX_s)\|^2 \big)\, \dd s
                 < \infty\bigr) = 1$
        \ for all \ $t \in \RR_+$,
      \item[\textup{(c)}]
       $\PP\bigl( \int_0^t
                   \int_{V_0} \|f(\bX_s, \br)\|^2 \, \dd s \, m(\dd \br)
                  < \infty \bigr)
        = 1$
        \ for all \ $t \in \RR_+$,
      \item[\textup{(d)}]
       $\PP\bigl( \int_0^t \int_{V_1} \|g(\bX_{s-}, \br)\| \, N(\dd s, \dd \br)
                   < \infty \bigr)
         = 1$
        \ for all \ $t \in \RR_+$, \ where \ $N(\dd s, \dd \br)$ \ is the
        counting measure of \ $p$ \ on \ $\RR_{++} \times V$,
      \item[\textup{(e)}]
       equation \eqref{SDE_X_YW_spec} holds \ $\PP$-a.s., where
        \ $\tN(\dd s, \dd \br) := N(\dd s, \dd \br) - \dd s \, m(\dd \br)$.
    \end{itemize}
 \end{enumerate}
\end{Def}

For the definitions of an $(\cF_t)_{t\in\RR_+}$-Brownian motion and an
 $(\cF_t)_{t\in\RR_+}$-Poisson point process, see, e.g., Ikeda and Watanabe
 \cite[Chapter I, Definition 7.2 and Chapter II, Definition 3.2]{IkeWat}.

\begin{Rem}\label{intexist}
If conditions (D1)--(D3) and (D4)(b)--(d) are satisfied, then the mappings
 \ $\RR_+ \times V_0 \times \Omega \ni (s, \br, \omega)
    \mapsto f(\bX_{s-}(\omega), \br) \in \RR^d$
 \ and
 \ $\RR_+ \times V_1 \times \Omega \ni (s, \br, \omega)
    \mapsto g(\bX_{s-}(\omega), \br) \in \RR^d$
 \ are in the (multidimensional versions of the) classes \ $\bF_p^{2,loc}$ \ and
 \ $\bF_p$, \ respectively, defined in Ikeda and Watanabe
 \cite[pages 61, 62]{IkeWat}, the integrals in \eqref{SDE_X_YW_spec} are
 well-defined and have c\`adl\`ag modifications as functions of \ $t$, \ see,
 e.g., Barczy et al.\ \cite[Remark 3.2]{BarLiPap}.

Moreover, if \ $\EE\bigl(\int_0^t \|\bX_s\| \, \dd s\bigr) < \infty$ \ for all
 \ $t \in \RR_+$, \ and the moment condition \eqref{moment_condition_1} holds,
 then conditions (D4)(b)--(d) are satisfied, and the mappings
 \ $\RR_+ \times V_0 \times \Omega \ni (s, \br, \omega)
    \mapsto f(\bX_{s-}(\omega), \br) \in \RR^d$
 \ and
 \ $\RR_+ \times V_1 \times \Omega \ni (s, \br, \omega)
    \mapsto g(\bX_{s-}(\omega), \br) \in \RR^d$
 \ are in the (multidimensional versions of the) smaller classes \ $\bF_p^2$
 \ and \ $\bF_p^1$, \ respectively, defined in Ikeda and Watanabe
 \cite[page 62]{IkeWat}.
Indeed, with the notation \ $\bX_s = (X_{s,1}, \ldots, X_{s,d})^\top$,
 \ $s \in \RR_+$,
 \begin{align*}
  \EE\left( \int_0^t \int_{V_0}
             \|f(\bX_s, \br)\|^2 \, \dd s \, m(\dd \br) \right)
  &= \sum_{j=1}^d
      \EE\left( \int_0^t \int_{U_d} \int_{U_1}
                 \|\bz\|^2 \bbone_{\{\|\bz\|<1\}} \bbone_{\{u\leq X_{s,j} \}}
                 \, \dd s \, \mu_j(\dd\bz) \, \dd u \right) \\
  &= \sum_{j=1}^d
      \EE\left( \int_0^t X_{s,j} \, \dd s \right)
      \int_{U_d} \|\bz\|^2 \bbone_{\{\|\bz\|<1\}} \, \mu_j(\dd\bz)
   < \infty
 \end{align*}
 by \eqref{z2}, and
  \begin{align*}
   &\EE\left( \int_0^t \int_{V_1}
             \|g(\bX_s, \br)\| \, \dd s \, m(\dd \br) \right) \\
   &= \int_0^t \int_{U_d} \|\br\| \, \dd s \, \nu(\dd\br)
      + \sum_{j=1}^d
         \EE\left( \int_0^t \int_{U_d} \int_{U_1}
                    \|\bz\| \bbone_{\{\|\bz\|\geq1\}}
                    \bbone_{\{u\leq X_{s,j} \}}
                    \, \dd s \, \mu_j(\dd\bz) \, \dd u \right)
  \end{align*}
  \begin{align*}
   = t \int_{U_d} \|\br\| \, \nu(\dd\br)
      + \sum_{j=1}^d \EE\left( \int_0^t X_{s,j} \, \dd s \right)
        \int_{U_d} \|\bz\| \bbone_{\{\|\bz\|\geq1\}} \, \mu_j(\dd\bz)
    < \infty
  \end{align*}
 by \eqref{help4} and \eqref{help_Li_page45_equiv}.
Note that if \ $(\bX_t)_{t\in\RR_+}$ \ is a CBI process with
 \ $\EE(\|\bX_0\|) < \infty$ \ satisfying the moment condition
 \eqref{moment_condition_1}, then
 \ $\EE\bigl(\int_0^t \|\bX_s\| \, \dd s\bigr) < \infty$ \ for all
 \ $t \in \RR_+$, \ see Lemma \ref{moment_1}.
\proofend
\end{Rem}

\begin{Rem}\label{dRM}
Note that if conditions (D1)--(D3) are satisfied, then \ $\bW$ \ and \ $p$
 \ are automatically independent according to Theorem 6.3 in Chapter II of
 Ikeda and Watanabe \cite{IkeWat}, since the intensity measure
 \ $\dd s \, m(\dd \br)$ \ of \ $p$ \ is deterministic.
Moreover, if \ $\bigl( \Omega, \cF, (\cF_t)_{t\in\RR_+}, \PP, \bW, p, \bX \bigr)$
 \ is an \ $\RR_+^d$-valued weak solution of the SDE \eqref{SDE_X_YW_spec},
 then \ $\cF_0$, \ $\bW$ \ and \ $p$ \ are mutually independent, and hence
 \ $\bX_0$, \ $\bW$ \ and \ $p$ \ are mutually independent as well, see, e.g.,
 Barczy et al.\ \cite[Remark 3.4]{BarLiPap}.
\proofend
\end{Rem}

\begin{Lem}\label{moment_1}
Let \ $(\bX_t)_{t\in\RR_+}$ \ be a CBI process with parameters
 \ $(d, \bc, \Bbeta, \bB, \nu, \bmu)$ \ and with initial distribution \ $n$
 \ satisfying \ $\int_{\RR_+^d} \|\bz\| \, n(\dd\bz) < \infty$.
\ Suppose that the moment condition \eqref{moment_condition_1} holds.
Then
 \begin{align*}
  \EE(\bX_t)
  = \ee^{t\tbB} \EE(\bX_0)
    + \left( \int_0^t \ee^{u\tbB} \, \dd u \right) \tBbeta , \qquad
  t \in \RR_+ ,
 \end{align*}
 where \ $\tbB \in \RR^{d \times d}_{(+)}$ \ and \ $\tBbeta \in \RR_+^d$ \ are
 defined in \eqref{tbeta_tB}.
In particular, \ $\int_0^t \EE(\|\bX_s\|) \, \dd s < \infty$ \ for all
 \ $t \in \RR_+$.
\end{Lem}

\noindent
\textbf{Proof.} \
By the tower rule for conditional expectations, it suffices to show
 \begin{equation}\label{MC1}
  \EE(\bX_t \mid \bX_0)
  = \ee^{t\tbB} \bX_0 + \left( \int_0^t \ee^{u\tbB} \, \dd u \right) \tBbeta ,
  \qquad t \in \RR_+ ,
 \end{equation}
 where the conditional expectation
 \ $\EE(\bX_t \mid \bX_0) \in [0, \infty]^d$ \ is meant in the generalized
 sense, see, e.g., Stroock \cite[Theorem 5.1.6]{Str}.
In order to show \eqref{MC1}, it is enough to check that for a CBI process
 \ $(\bX_t)_{t\in\RR_+}$ \ with initial value \ $\bX_0 = \bx \in \RR_+^d$, \ we
 have
 \begin{equation}\label{MC2}
  \EE(\bX_t)
  = \ee^{t\tbB} \bx + \left( \int_0^t \ee^{u\tbB} \, \dd u \right) \tBbeta ,
  \qquad t \in \RR_+ , \quad \bx \in \RR_+^d .
 \end{equation}
Indeed, let \ $\phi_n : \RR_+^d \to \RR_+^d$, \ $n \in \NN$, \ be simple
 functions such that \ $\phi_n(\by) \uparrow \by$ \ as \ $n \to \infty$
 \ for all \ $\by \in \RR_+^d$.
\ Then, by the (multidimensional version of the) monotone convergence theorem
 for (generalized) conditional expectations, see, e.g., Stroock
 \cite[Theorem 5.1.6]{Str}, we obtain
 \ $\EE(\phi_n(\bX_t) \mid \bX_0) \uparrow \EE(\bX_t \mid \bX_0)$
 \ as \ $n \to \infty$ \ $\PP$-almost surely.
For each \ $B \in \cB(\RR^d)$, \ we have
 \[
   \EE(\bbone_B(\bX_t) \mid \bX_0)
   = \PP(\bX_t \in B \mid \bX_0)
   = \int_{\RR_+^d} \bbone_B(\by) \, P_t(\bX_0, \dd \by) ,
 \]
 hence
 \ $\EE(\phi_n(\bX_t) \mid \bX_0)
    = \int_{\RR_+^d} \phi_n(\by) \, P_t(\bX_0, \dd \by)$.
\ By the (multidimensional version of the) monotone convergence theorem,
 \ $\int_{\RR_+^d} \phi_n(\by) \, P_t(\bX_0, \dd \by)
    \uparrow \int_{\RR_+^d} \by \, P_t(\bX_0, \dd \by)$
 \ as \ $n \to \infty$.
\ By \eqref{MC2}, we get
 \[
   \EE(\bX_t \mid \bX_0)
   = \int_{\RR_+^d} \by \, P_t(\bX_0, \dd \by)
    = \ee^{t\tbB} \bX_0 + \left( \int_0^t \ee^{u\tbB} \, \dd u \right) \tBbeta ,
 \]
 hence we conclude \eqref{MC1}.

In order to show \eqref{MC2}, we are going to apply Proposition 9.11 of Li
 \cite{Li} for the immigration superprocess given in Lemma \ref{superprocess}.
For each \ $f \in B(E)$ \ and \ $i \in E$, \ the function
 \ $\RR_+ \ni t \mapsto \pi_tf(i)$ \ is the unique locally bounded solution to
 the linear evolution equation (2.35) in Li \cite{Li} taking the form
 \begin{align*}
  &\pi_tf(i)
   = f(i) + \int_0^t \gamma(i, \pi_sf) \, \dd s
     - \int_0^t b(i) \pi_sf(i) \, \dd s \\
  &= f(i)
     + \int_0^t \left( \sum_{j=1}^d \pi_sf(j) \gamma(i, \{j\}) \right) \dd s
     - \int_0^t b(i) \pi_sf(i) \, \dd s
   = f(i)
     + \int_0^t \left( \sum_{j=1}^d \pi_sf(j) \tb_{j,i} \right) \dd s ,
 \end{align*}
 where we used \ $R_tf = f$ \ for \ $f \in B(E)$ \ and \ $t \in \RR_+$,
 \ $b(i) = - \tb_{i,i}$ \ and \ $\gamma(i, \{i\}) = B(i, \{i\}) = 0$ \ for
 \ $i \in \{1, \ldots, d\}$, \ and
 \begin{align}\label{tgammaij}
   \gamma(i, \{j\})
   = B(i, \{j\}) + \int_{U_d} z_j \, \mu_i(\dd\bz)
   = b_{j,i} + \int_{U_d} z_j \, \mu_i(\dd\bz)
   = \tb_{j,i}
 \end{align}
 for \ $i, j \in \{1, \ldots, d\}$ \ with \ $i \ne j$.
\ The functions \ $\RR_+ \ni t \mapsto \pi_tf(i)$, \ $f \in B(E)$,
 \ $i \in \{1, \ldots, d\}$, \ can be identified with the functions
 \ $\RR_+ \ni t \mapsto \pi_i(t, \blambda)$, \ $\blambda \in \RR^d$,
 \ $i \in \{1, \ldots, d\}$, \ which are the unique locally bounded solution
 to the linear evolution equations
 \[
   \pi_i(t,\blambda)
   = \lambda_i
     + \int_0^t \langle \tbB \be_i, \pi(s, \blambda) \rangle \, \dd s ,
   \qquad t \in \RR_+ , \quad i \in \{1, \ldots, d\} , \quad
   \blambda \in \RR^d .
 \]
Consequently, the functions
 \ $\RR_+ \ni t \mapsto
    \pi(t, \blambda) := (\pi_1(t, \blambda), \ldots, \pi_d(t, \blambda))$,
 \ $\blambda \in \RR^d$, \ satisfies
 \[
   \bpi(t, \blambda)
   = \blambda + \int_0^t \tbB^\top \!\! \bpi(s, \blambda) \, \dd s ,
   \qquad t \in \RR_+ , \quad \blambda \in \RR^d ,
 \]
 and hence
 \begin{equation*}
  \bpi(t, \blambda) = \ee^{t\tbB^\top} \! \blambda ,
  \qquad t \in \RR_+ , \quad \blambda \in \RR^d .
 \end{equation*}
The functional
 \ $B(E) \ni f \mapsto \Gamma(f)
    = \eta(f) + \int_{M(E)^\circ} \kappa(f) \, H_1(\dd\kappa)$
 \ of \cite[formula (9.20)]{Li} can be identified with the functional
 \ $\RR^d \ni \bx \mapsto
    \bx^\top \Bbeta + \int_{U_d} \bx^\top \bz \, \nu(\dd\bz)
    = \bx^\top \tBbeta$.
\ Hence Proposition 9.11 of Li \cite{Li} implies
 \begin{align*}
  \langle \blambda, \EE(\bX_t)\rangle
     = \langle \ee^{t\tbB^\top} \! \blambda , \bx \rangle
      + \left(\int_0^t (\ee^{s\tbB^\top}\blambda)^\top\,\dd s \right) \tBbeta
     = \left\langle \blambda , \ee^{t\tbB} \bx + \left(\int_0^t \ee^{s\tbB}\,\dd s \right) \tBbeta \right\rangle
 \end{align*}
 for \ $t \in \RR_+$ \ and \ $\blambda \in \RR^d$, \ which yields \eqref{MC2}.
\proofend

\begin{Rem}
We call the attention that in the proof of the forthcoming Theorem
 \ref{CBI_SDE}, which states existence of an \ $\RR_+^d$-valued weak solution
 of the SDE \eqref{SDE_X_YW_spec}, we will extensively use that for a CBI
 process \ $(\bX_t)_{t\in\RR_+}$ \ with parameters
 \ $(d, \bc, \Bbeta, \bB, \nu, \bmu)$ \ satisfying \ $\EE(\|\bX_0\|) < \infty$
 \ and the moment condition \eqref{moment_condition_1}, we have
 \ $\int_0^t \EE(\|\bX_s\|) \, \dd s < \infty$, \ $t \in \RR_+$, \ proved in
 Lemma \ref{moment_1}.
We point out that in the proof of Lemma \ref{moment_1} we can not use the SDE
 \eqref{SDE_X_YW_spec}, since at that point it has not yet been proved that a
 CBI process is a solution of this SDE.
This drives us back to Definition \ref{Def_CBI} of CBI processes in the proof
 of Lemma \ref{moment_1}.
Having proved that a CBI process is a solution of the SDE
 \eqref{SDE_X_YW_spec}, one could give another proof of Lemma \ref{moment_1}
 (roughly speaking by taking expectations via localization argument).
\proofend
\end{Rem}

\begin{Def}\label{Def_uniqueness_in_law}
We say that uniqueness in the sense of probability law holds for the SDE
 \eqref{SDE_X_YW_spec} among \ $\RR_+^d$-valued weak solutions if whenever
 \ $\bigl( \Omega, \cF, (\cF_t)_{t\in\RR_+}, \PP, \bW, p, \bX \bigr)$ \ and
 \ $\bigl( \tOmega, \tcF, (\tcF_t)_{t\in\RR_+}, \tPP, \tbW, \tp, \tbX \bigr)$
 \ are \ $\RR_+^d$-valued weak solutions of the SDE \eqref{SDE_X_YW_spec} such
 that \ $\PP(\bX_0 \in B) = \tPP(\tbX_0 \in B)$ \ for all
 \ $B \in \cB(\RR^d)$, \ then \ $\PP(\bX \in C) = \tPP(\tbX \in C)$ \ for all
 \ $C \in \cD(\RR_+, \RR^d)$.
\end{Def}

\begin{Thm}\label{CBI_SDE}
Let \ $(d, \bc, \Bbeta, \bB, \nu, \bmu)$ \ be a set of admissible parameters
 in the sense of Definition \ref{Def_admissible} such that the moment
 condition \eqref{moment_condition_1} holds.
Then for any probability measure \ $n$ \ on \ $(\RR_+^d, \cB(\RR_+^d))$
 \ with \ $\int_{\RR^d_+} \|\bz\| \, n(\dd \bz) < \infty$, \ the SDE
 \eqref{SDE_X_YW_spec} admits an \ $\RR_+^d$-valued weak solution with initial
 distribution \ $n$ \ which is unique in the sense of probability law among
 \ $\RR_+^d$-valued weak solutions.
Moreover, any \ $\RR_+^d$-valued weak solution is a CBI process with
 parameters \ $(d, \bc, \Bbeta, \bB, \nu, \bmu)$.
\end{Thm}

\noindent
\textbf{Proof.} \
Suppose that \ $(\bX_t)_{t\in\RR_+}$ \ is a c\`adl\`ag realization of a CBI
 process with parameters \ $(d, \bc, \Bbeta, \bB, \nu, \bmu)$ \ on a
 probability space \ $(\Omega, \cF, \PP)$ \ having initial distribution \ $n$,
 \ i.e., \ $(\bX_t)_{t\in\RR_+}$ \ is a time homogeneous Markov process having
 c\`adl\`ag trajectories and the same finite dimensional distributions as a
 CBI process with parameters \ $(d, \bc, \Bbeta, \bB, \nu, \bmu)$ \ having
 initial distribution \ $n$ \ (such a realization exists due to Theorem 9.15
 in Li \cite{Li}).
Let
 \[
   \cF_t := \bigcap_{\vare>0} \sigma\left(\cF_{t+\vare}^\bX \cup \cN\right) ,
   \qquad t \in \RR_+ ,
 \]
 where \ $\cN$ \ denotes the collection of null sets under the probability
 measure \ $\PP$, \ and \ $(\cF_t^\bX)_{t\in\RR_+}$ \ stands for the natural
 filtration generated by the process \ $(\bX_t)_{t\in\RR_+}$, \ hence the
 filtered probability space \ $(\Omega, \cF, (\cF_t)_{t\in\RR_+}, \PP)$
 \ satisfies the usual hypotheses.

By the equivalence of parts (3) and (4) of Theorem 9.18 of Li \cite{Li}
 applied to the immigration superprocess given in Lemma \ref{superprocess},
 we conclude that the process \ $(\bX_t)_{t\in\RR_+}$ \ has no negative jumps,
 the (not necessarily Poisson) random measure
 \[
   N_0(\dd s, \dd\bz)
   := \sum_{u\in\RR_{++}}
       \bbone_{\{\bX_u\ne\bX_{u-}\}} \delta_{(u,\bX_u-\bX_{u-})} (\dd s, \dd\bz)
 \]
 on \ $\RR_{++} \times U_d$ \ has predictable compensator
 \[
   \hN_0(\dd s, \dd \bz)
   := \sum_{j=1}^d X_{s-,j} \, \dd s \, \mu_j(\dd\bz) + \dd s \, \nu(\dd\bz) ,
 \]
 and
 \[
   \bX_t - \bX_0
   - \int_0^t \bigl( \tBbeta + \tbB \bX_s \bigr) \, \dd s
   - \int_0^t \int_{U_d} \bz \, \tN_0(\dd s, \dd\bz) , \qquad
   t \in \RR_+ ,
 \]
 is a continuous locally square integrable martingale starting from
 \ $\bzero \in \RR^d$ \ with quadratic variation process
 \[
   \left( 2 \delta_{i,j} c_i
          \int_0^t X_{s,i} \, \dd s \right)_{i,j\in\{1,\ldots,d\}} , \qquad
   t \in \RR_+ ,
 \]
 where \ $\tN_0(\dd s, \dd\bz) := N_0(\dd s, \dd\bz) - \hN_0(\dd s, \dd\bz)$.
\ Indeed, first, note that \ $R_tf=f$, $t\in\RR_+$, $f\in B(E)$, \ yields that
 the strong generator of \ $(R_t)_{t\in\RR_+}$ \ is identically \ $0$, \ i.e.,
 \ $A=0$, \ see Li \cite[(7.1)]{Li}.
Using \ $b(i) = - \tb_{i,i}$ \ and
 \ $\gamma(i, \{i\}) = B(i, \{i\}) = 0$ \ for \ $i \in \{1, \ldots, d\}$ \ and
 \ $\gamma(i, \{j\}) = \tb_{j,i}$ \ for \ $i, j \in \{1, \ldots, d\}$ \ with \ $i \ne j$ \ (see, \eqref{tgammaij}),
 \ the function \ $B(E)\ni f\mapsto Af+\gamma f-bf$ \ of Li \cite[page 218]{Li} can be identified with
 the function
 \begin{align}\label{Agammab}
    E\ni i \mapsto \sum_{j=1}^d f(j)\gamma(i,\{j\}) - b(i)f(i)
                    = \sum_{j=1}^d \tb_{j,i}f(j).
 \end{align}
Recalling that the functional
 \ $B(E) \ni f \mapsto \Gamma(f)
    = \eta(f) + \int_{M(E)^\circ} \kappa(f) \, H(\dd\kappa)$
 \ is identified with the functional \ $\RR^d \ni \bx \mapsto  \bx^\top \tBbeta$ \ (see, the end of the proof
 of Lemma \ref{moment_1}), Theorem 9.18 of Li \cite{Li} yields that for each
 \ $\bw = (w_1, \ldots, w_d)^\top \in \RR^d$, \ the process
 \ $(\bw^\top \bX_t)_{t\in\RR_+}$ \ has no negative jumps, and
 \[
   \bw^\top \bX_t - \bw^\top \bX_0
   - \int_0^t
      \bigl( \bw^\top \tBbeta + \bw^\top \tbB \bX_s \bigr) \, \dd s
   - \int_0^t \int_{U_d} \bw^\top \bz \, \tN_0(\dd s, \dd\bz) ,
   \qquad t \in \RR_+ ,
 \]
 is a continuous locally square integrable martingale strating from
 \ $0 \in \RR$ \ with quadratic variation process
 \[
   \langle \bw^\top \bX \rangle_t
   = 2 \sum_{i=1}^d c_i w_i^2 \int_0^t X_{s,i} \, \dd s , \qquad t \in \RR_+ .
 \]
Further, by polarization identity, for all \ $\bw, \tbw \in \RR^d$, \ the
 cross quadratic variation process of \ $(\bw^\top \bX_t)_{t\in\RR_+}$ \ and
 \ $(\tbw^\top \bX_t)_{t\in\RR_+}$ \ takes the form
 \begin{align*}
  \langle \bw^\top \bX, \tbw^\top \bX \rangle_t
  &= \frac{1}{4}\left( \langle (\bw + \tbw)^\top \bX \rangle_t
                       - \langle (\bw - \tbw)^\top \bX \rangle_t \right) \\
  &= \frac{1}{4}
     \biggl( 2 \sum_{i=1}^d c_i (w_i + \tw_i)^2 \int_0^t X_{s,i} \, \dd s
             - 2 \sum_{i=1}^d c_i (w_i - \tw_i)^2 \int_0^t X_{s,i} \, \dd s
     \biggr)\\
 &= 2 \sum_{i=1}^d c_i w_i \tw_i \int_0^t X_{s,i} \, \dd s, \qquad t \in \RR_+ .
 \end{align*}
We note that the integral \ $\int_0^t \int_{U_d} \bz \, \tN_0(\dd s, \dd\bz)$
 \ is well-defined, since
 \ $\bz = \bz \bbone_{\{\|\bz\|<1\}} + \bz \bbone_{\{\|\bz\|\geq1\}}$,
 \ $\bz \in U_d$, \ and the functions
 \ $\RR_+ \times U_d \times \Omega \ni (s, \bz, \omega)
    \mapsto \bz \bbone_{\{\|\bz\|<1\}}$
 \ and
 \ $\RR_+ \times U_d \times \Omega \ni (s,\bz,\omega)
    \mapsto \bz \bbone_{\{\|\bz\|\geq1\}}$
 \ belong to the classes \ $\bF_{p_0}^2$ \ and \ $\bF_{p_0}^1$, \ respectively,
 where \ $p_0$ \ denotes the point process on \ $U_d$ \ with counting measure
 \ $N_0(\dd s, \dd\bz)$,
 \ i.e., \ $p_0(u) := \bX_u - \bX_{u-}$ \ for \ $u \in D(p_0)$ \ with
 \ $D(p_0) := \{ u \in \RR_{++} : \bX_u \ne \bX_{u-} \}$.
\ Indeed,
 \begin{align*}
  &\EE\left(\int_0^t \int_{U_d}
             \|\bz\|^2 \bbone_{\{\|\bz\|<1\}} \, \hN_0(\dd s, \dd\bz)\right) \\
  &= \int_0^t \int_{U_d}
             \|\bz\|^2 \bbone_{\{\|\bz\|<1\}} \, \dd s \, \nu(\dd\bz)
     + \sum_{j=1}^d
        \int_0^t \int_{U_d}
         \|\bz\|^2 \bbone_{\{\|\bz\|<1\}} \EE(X_{s,j}) \, \dd s \, \mu_j(\dd\bz) \\
  &\leq t \int_{U_d} \|\bz\| \, \nu(\dd\bz)
        + \sum_{j=1}^d
           \int_0^t \EE(X_{s,j}) \, \dd s
           \int_{U_d} \|\bz\|^2 \bbone_{\{\|\bz\|<1\}} \, \mu_j(\dd\bz)
   < \infty
 \end{align*}
 by Lemma \ref{moment_1} and the inequalities \eqref{help4} and \eqref{z2},
 and
 \begin{align*}
  &\EE\left(\int_0^t \int_{U_d}
             \|\bz\| \bbone_{\{\|\bz\|\geq1\}} \, \hN_0(\dd s, \dd\bz)\right) \\
  &= \int_0^t \int_{U_d}
             \|\bz\| \bbone_{\{\|\bz\|\geq1\}} \, \dd s \, \nu(\dd\bz)
     + \sum_{j=1}^d
        \int_0^t \int_{U_d}
         \|\bz\| \bbone_{\{\|\bz\|\geq1\}} \EE(X_{s,j}) \, \dd s \, \mu_j(\dd\bz)
 \end{align*}
 \begin{align*}
  &\leq t \int_{U_d} \|\bz\| \, \nu(\dd\bz)
        + \sum_{j=1}^d
           \int_0^t \EE(X_{s,j}) \, \dd s
           \int_{U_d} \|\bz\| \bbone_{\{\|\bz\|\geq1\}} \, \mu_j(\dd\bz)
   < \infty
 \end{align*}
 by Lemma \ref{moment_1} and the inequalities \eqref{help4} and
 \eqref{help_Li_page45_equiv}.

Using that \ $\PP\big(\int_0^t X_{s,i} \, \dd s < \infty\big) = 1$,
 \ $i \in \{1, \ldots, d\}$ \ (since \ $\bX$ \ has c\`adl\`ag trajectories
 almost surely), by choosing \ $\bw = \be_j$, \ $j \in \{1, \ldots, d\}$,
 \ a representation theorem for continuous locally square integrable
 martingales (see, e.g.,
 Ikeda and Watanabe \cite[Chapter II, Theorem 7.1']{IkeWat}) yields
 \begin{align*}
  \bX_t
  = \bX_0 + \int_0^t \bigl( \tBbeta + \tbB \bX_s \bigr) \, \dd s
    + \sum_{i=1}^d \be_i \int_0^t \sqrt{2c_i X_{s,i}} \, \dd W_{s,i}
    + \int_0^t \int_{U_d} \bz \, \tN_0(\dd s, \dd\bz)
 \end{align*}
 for all \ $t \in \RR_+$, \ $\tPP$-almost surely on an extension
 \ $\bigl(\tOmega, \tcF, (\tcF_t)_{t\in\RR_+}, \tPP\bigr)$ \ of the filtered
 probability space \ $(\Omega, \cF, (\cF_t)_{t\in\RR_+}, \PP)$ \ (see Definition
 \ref{Def_extension}),
 and \ $(W_{t,1}, \ldots, W_{t,d})_{t\in\RR_+}$ \ is a $d$-dimensional
 \ $(\tcF_t)_{t\in\RR_+}$-Brownian motion.
We note that, with a little abuse of notation, the extended random variables
 on the extension \ $\bigl(\tOmega, \tcF, (\tcF_t)_{t\in\RR_+}, \tPP\bigr)$
 \ are denoted in the same way as the original ones.
Let
 \[
   \tcG_t := \bigcap_{\vare>0} \sigma\left(\tcF_{t+\vare} \cup \tcN\right) ,
   \qquad t \in \RR_+ ,
 \]
 where \ $\tcN$ \ denotes the collection of null sets under the probability
 measure \ $\tPP$.
\ Then the filtered probability space
 \ $(\tOmega, \tcF, (\tcG_t)_{t\in\RR_+}, \tPP)$ \ satisfies the usual
 hypotheses, and by Lemma \ref{ext_BM_PRN},
 \ $(W_{t,1}, \ldots, W_{t,d})_{t\in\RR_+}$ \ is a \ $d$-dimensional
 \ $(\tcG_t)_{t\in\RR_+}$-Brownian motion.

The aim of the following discussion is to show, by the representation theorem
 of Ikeda and Watanabe \cite[Chapter II, Theorem 7.4]{IkeWat}, that the SDE
 \eqref{SDE_X_YW_spec} holds on an extension of the original probability space.
The predictable compensator of the random measure \ $N_0(\dd s, \dd\bz)$ \ can
 be written in the form \ $\hN_0(\dd s, \dd\bz) = \dd s \, q(s, \dd\bz)$,
 \ where
 \[
   q(s, \dd\bz) := \sum_{j=1}^d X_{s-,j} \, \mu_j(\dd\bz) + \nu(\dd\bz) .
 \]
Let
 \ $\Theta : \RR_+ \times V \times \tOmega \to U_d \cup \{\bzero\} = \RR_+^d$
 \ be defined by
 \[
   \Theta(s, \br, \tomega) := h(\bX_{s-}(\tomega), \br) , \qquad
   (s, \br, \tomega) \in \RR_+ \times V \times \tOmega .
 \]
(Note, that \ $\Delta = \bzero$ \ in the notation of Ikeda and Watanabe
 \cite[Chapter II, Theorem 7.4]{IkeWat}.)
Then condition (7.26) on page 93 in Ikeda and Watanabe \cite{IkeWat} holds,
 since for all \ $s\in\RR_+$, \ $\tomega\in\tOmega$, \ and \ $B \in \cB(U_d)$,
 \ we have
 \begin{align*}
  &m(\{ \br \in V : \Theta(s, \br, \tomega) \in B \})
   = \sum_{i=0}^d m(\{ \br \in \cR_i : \Theta(s, \br, \tomega) \in B \}) \\
  &= \sum_{i=1}^d
      (\mu_i \times \ell)
      \big( \{ (\bz, u) \in \cR_i
               : \bz \bbone_{\{u\leq X_{s-,i}(\tomega)\}} \in B \} \big)
     + \nu\big( \{ \br \in \cR_0  : \br \in B \} \big)
  \end{align*}
  \begin{align*}
  &= \sum_{i=1}^d X_{s-,i}(\tomega) \, \mu_i(B) + \nu(B)
   = q(s, B)(\tomega) ,
 \end{align*}
 where \ $\ell$ \ denotes the Lebesgue measure on \ $\RR_{++}$, \ and we used
 that \ $\bzero \notin B$.
\ By Theorem II.7.4 in Ikeda and Watanabe \cite{IkeWat}, on an extension
 \ $\bigl(\ttOmega, \ttcF, (\ttcF_t)_{t\in\RR_+}, \ttPP\bigr)$ \ of
 \ $(\tOmega, \tcF, (\tcG_t)_{t\in\RR_+}, \tPP)$, \ there is a stationary
 \ $(\ttcF_t)_{t\in\RR_+}$-Poisson point process \ $p$ \ on \ $V$
 \ with characteristic measure \ $m$ \ such that
 \begin{align*}
  N_0\bigl( (0, t] \times B \bigr)
  &= \int_0^t \int_V \bbone_B(\Theta(s, \br)) \, N(\dd s, \dd\br) \\
  &= \# \{ s \in D(p) : s \in (0, t] , \, \Theta(s, p(s)) \in B \}
  \qquad \text{$\ttPP$-a.s.}
 \end{align*}
 for all \ $B \in \cB(U_d)$, \ where \ $N(\dd s, \dd\br)$ \ denotes the
 counting measure of \ $p$, \ and \ $D(p)$ \ is the domain of \ $p$ \ being a
 countable subset of \ $\RR_{++}$ \ such that
 \ $\{ s \in D(p) : s \in (0, t] , \, p(s) \in B\}$ \ is finite for all
 \ $t \in \RR_+$ \ and compact subsets \ $B \in \cB(U_d)$.
\ Then, by Lemma \ref{ext_BM}, \ $(W_{t,1}, \ldots, W_{t,d})_{t\in\RR_+}$ \ is a
 \ $d$-dimensional
 \ $(\ttcF_t)_{t\in\RR_+}$-Brownian motion.
Let
 \[
   \ttcG_t := \bigcap_{\vare>0} \sigma\left(\ttcF_{t+\vare} \cup \ttcN\right) ,
   \qquad t \in \RR_+ ,
 \]
 where \ $\ttcN$ \ denotes the collection of null sets under the probability
 measure \ $\ttPP$.
\ Then the filtered probability space
 \ $(\ttOmega, \ttcF, (\ttcG_t)_{t\in\RR_+}, \ttPP)$ \ satisfies the usual
 hypotheses.
By Lemma \ref{ext_BM_PRN}, \ $(W_{t,1}, \ldots, W_{t,d})_{t\in\RR_+}$ \ is a
 \ $d$-dimensional \ $(\ttcG_t)_{t\in\RR_+}$-Brownian motion, and \ $p$ \ is a
 stationary \ $(\ttcG_t)_{t\in\RR_+}$-Poisson point process \ on \ $V$ \ with
 characteristic measure \ $m$.
\ Consequently,
 \begin{equation}\label{II.7.4}
  \# \{ s \in D(p_0) : s \in (0, t] , \, p_0(s) \in B \}
  = \# \{ s \in D(p) : s \in (0, t] , \, h(\bX_{s-}, p(s)) \in B \}
 \end{equation}
 for all \ $B \in \cB(U_d)$.
\ Using this representation, we will calculate
 \ $\int_0^t \int_{U_d} \bz \, \tN_0(\dd s, \dd\bz)$, \ $t \in \RR_+$.
\ First observe that
 \[
   \int_0^t \int_{U_d} \bz \, \tN_0(\dd s, \dd\bz)
   = \int_0^t \int_{U_d} \bz \bbone_{\{\|\bz\|\geq1\}} \, \tN_0(\dd s, \dd\bz)
     + \int_0^t \int_{U_d} \bz \bbone_{\{\|\bz\|<1\}} \, \tN_0(\dd s, \dd\bz) .
 \]
Since the function
 \ $\RR_+ \times U_d \times \Omega \ni (s, \bz, \omega)
    \mapsto \bz \bbone_{\{\|\bz\|\geq1\}}$ \ belongs to the class \ $\bF_{p_0}^1$,
 \ by Ikeda and Watanabe \cite[Chapter II, (3.8)]{IkeWat}, we obtain
 \[
   \int_0^t \int_{U_d} \bz \bbone_{\{\|\bz\|\geq1\}} \, \tN_0(\dd s, \dd\bz)
   = \int_0^t \int_{U_d} \bz \bbone_{\{\|\bz\|\geq1\}} \, N_0(\dd s, \dd\bz)
     - \int_0^t \int_{U_d} \bz \bbone_{\{\|\bz\|\geq1\}} \, \hN_0(\dd s, \dd\bz) .
 \]
Applying \eqref{II.7.4}, we obtain
 \begin{align*}
  &\int_0^t \int_{U_d} \bz \bbone_{\{\|\bz\|\geq1\}} \, N_0(\dd s, \dd\bz)
   = \sum_{s\in D(p_0)\cap(0,t]} p_0(s) \bbone_{\{\|p_0(s)\|\geq1\}} \\
  &= \sum_{s\in D(p)\cap(0,t]} h(\bX_{s-}, p(s)) \bbone_{\{\|h(\bX_{s-}, p(s))\|\geq1\}}
   = \int_0^t \int_V
      h(\bX_{s-}, \br) \bbone_{\{\|h(\bX_{s-}, \br)\|\geq1\}}
      \, N(\dd s, \dd\br) \\
  &= \int_0^t \int_{\cR_0} \br \bbone_{\{\|\br\|\geq1\}} \, N(\dd s, \dd\br)
     + \sum_{j=1}^d
        \int_0^t \int_{\cR_j}
         \bz \bbone_{\{\|\bz\|\geq1\}} \bbone_{\{u\leq X_{s-,j}\}}
         \, N(\dd s, \dd\br) \\
  &= \int_0^t \int_{V_1} g(\bX_{s-}, \br) \, N(\dd s, \dd\br)
     - \int_0^t \int_{\cR_0} \br \bbone_{\{\|\br\|<1\}} \, N(\dd s, \dd\br) .
 \end{align*}
Here we used that the function
 \ $\RR_+ \times U_d \times \ttOmega \ni (s, \bz, \ttomega)
    \mapsto \bz \bbone_{\{\|\bz\|\geq1\}}$
 \ belongs to the class \ $\bF_{p_0}^1$, \ hence the function
 \ $\RR_+ \times V \times \ttOmega \ni (s, \br, \ttomega)
    \mapsto h(\bX_{s-}(\ttomega), \br) \bbone_{\{\|h(\bX_{s-}(\ttomega), \br)\|\geq1\}}$
 \ belongs to the class \ $\bF_p^1$, \ and function
 \ $\RR_+ \times V \times \ttOmega \ni (s, \br, \ttomega)
    \mapsto \br \bbone_{\{\|\br\|<1\}} \bbone_{\cR_0}(\br)$
 \ also belongs to the class \ $\bF_p^1$ \ (due to \eqref{help4}), thus the
 function
 \ $\RR_+ \times V \times \ttOmega \ni (s, \br, \ttomega)
     \mapsto g(\bX_{s-}(\ttomega), \br)$
 \ belongs to the class \ $\bF_p^1$ \ as well.
Moreover,
 \begin{align*}
  &\int_0^t \int_{U_d} \bz \bbone_{\{\|\bz\|\geq1\}} \, \hN_0(\dd s, \dd\bz)
   = \int_0^t \int_{U_d} \bz \bbone_{\{\|\bz\|\geq1\}} \, \dd s \, \nu(\dd\bz)
     + \sum_{j=1}^d
        \int_0^t \int_{U_d}
         \bz \bbone_{\{\|\bz\|\geq1\}} X_{s,j} \, \dd s \, \mu_j(\dd\bz) \\
  &= \int_0^t \int_{U_d} \br \bbone_{\{\|\br\|\geq1\}} \, \dd s \, \nu(\dd\br)
     + \sum_{j=1}^d
        \int_0^t X_{s,j}\,\dd s \int_{U_d}
         \bz \bbone_{\{\|\bz\|\geq1\}} \, \mu_j(\dd\bz) .
 \end{align*}
Let \ $\cM_2$ \ denote the complete metric space of square integrable
 right continuous \ $d$-dimensional martingales on
 \ $(\ttOmega, \ttcF, \ttPP)$ \ with respect to \ $(\ttcF_t)_{t\in\RR_+}$
 \ starting from \ $\bzero$, \ see, e.g., Ikeda and Watanabe
 \cite[Chapter II, Definition 1.3 and Lemma 1.2]{IkeWat}.
The function
 \ $\RR_+ \times U_d \times \ttOmega \ni (s, \bz, \ttomega)
    \mapsto \bz \bbone_{\{\|\bz\|<1\}}$
 \ belongs to the class \ $\bF_{p_0}^2$, \ hence, by Ikeda and Watanabe
 \cite[Chapter II, (3.9)]{IkeWat}, the process
 \ $\bigl(\int_0^t \int_{U_d}
           \bz \bbone_{\{\|\bz\|<1\}} \, \tN_0(\dd s, \dd\bz)\bigr)_{t\in\RR_+}$
 \ belongs to the space \ $\cM_2$.
\ Moreover, by Ikeda and Watanabe \cite[page 63]{IkeWat},
 \ $\int_0^t \int_{U_d} \bz \bbone_{\{\|\bz\|<1\}} \, \tN_0(\dd s, \dd\bz)$
 \ is the limit of the sequence
 \ $\int_0^t \int_{U_d}
     \bz \bbone_{\{\frac{1}{n}\leq\|\bz\|<1\}} \, \tN_0(\dd s, \dd\bz)$,
 \ $n \in \NN$, \ in \ $\cM_2$ \ as \ $n \to \infty$.
\ For all \ $n \in \NN$, \ the mapping
 \ $\RR_+ \times U_d \times \ttOmega \ni (s, \bz, \ttomega)
    \mapsto \bz \bbone_{\{\frac{1}{n}\leq\|\bz\|<1\}}$
 \ belongs to the class \ $\bF_{p_0}^1 \cap \bF_{p_0}^2$, \ hence we obtain
 \[
   \int_0^t \int_{U_d}
    \bz \bbone_{\{\frac{1}{n}\leq\|\bz\|<1\}} \, \tN_0(\dd s, \dd\bz)
   = \int_0^t \int_{U_d} \bz
      \bbone_{\{\frac{1}{n}\leq\|\bz\|<1\}} \, N_0(\dd s, \dd\bz)
      - \int_0^t \int_{U_d}
         \bz \bbone_{\{\frac{1}{n}\leq\|\bz\|<1\}} \, \hN_0(\dd s, \dd\bz) .
 \]
Similarly as above, \begin{align*}
  &\int_0^t \int_{U_d}
    \bz \bbone_{\{\frac{1}{n}\leq\|\bz\|<1\}} \, N_0(\dd s, \dd\bz) \\
  &= \int_0^t \int_{\cR_0}
      \br \bbone_{\{\frac{1}{n}\leq\|\br\|<1\}} \, N(\dd s, \dd\br)
     + \sum_{j=1}^d
        \int_0^t \int_{\cR_j}
         \bz \bbone_{\{\frac{1}{n}\leq\|\bz\|<1\}} \bbone_{\{u\leq X_{s-,j}\}}
         \, N(\dd s, \dd\br)
 \end{align*}
 and
 \begin{align*}
  &\int_0^t \int_{U_d}
    \bz \bbone_{\{\frac{1}{n}\leq\|\bz\|<1\}} \, \hN_0(\dd s, \dd\bz) \\
  &= \int_0^t \int_{U_d}
      \br \bbone_{\{\frac{1}{n}\leq\|\br\|<1\}} \, \dd s \, \nu(\dd\br)
     + \sum_{j=1}^d
        \int_0^t \int_{U_d} \int_{U_1}
         \bz \bbone_{\{\frac{1}{n}\leq\|\bz\|<1\}} \bbone_{\{u\leq X_{s,j}\}}
         \, \dd s \, \mu_j(\dd\bz) \, \dd u .
 \end{align*}
Consequently,
 \begin{align*}
  &\int_0^t \int_{U_d}
    \bz \bbone_{\{\frac{1}{n}\leq\|\bz\|<1\}} \, \tN_0(\dd s, \dd\bz) \\
  &= \int_0^t \int_{\cR_0}
      \br \bbone_{\{\frac{1}{n}\leq\|\br\|<1\}} \, \tN(\dd s, \dd\br)
     + \sum_{j=1}^d
        \int_0^t \int_{\cR_j}
         \bz \bbone_{\{\frac{1}{n}\leq\|\bz\|<1\}} \bbone_{\{u\leq X_{s-,j}\}}
         \, \tN(\dd s, \dd\br) .
 \end{align*}
Taking the limit in \ $\cM_2$ \ as \ $n \to \infty$, \ we conclude
 \begin{align*}
  &\int_0^t \int_{U_d} \bz \bbone_{\{\|\bz\|<1\}} \, \tN_0(\dd s, \dd\bz) \\
  &= \int_0^t \int_{\cR_0}
      \br \bbone_{\{\|\br\|<1\}} \, \tN(\dd s, \dd\br)
     + \sum_{j=1}^d
        \int_0^t \int_{\cR_j}
         \bz \bbone_{\{\|\bz\|<1\}} \bbone_{\{u\leq X_{s-,j}\}}
         \, \tN(\dd s, \dd\br) \\
  &= \int_0^t \int_{\cR_0}
      \br \bbone_{\{\|\br\|<1\}} \, \tN(\dd s, \dd\br)
     + \int_0^t \int_{V_0} f(\bX_{s-}, \br) \, \tN(\dd s, \dd\br) .
 \end{align*}
Summarizing, we conclude
 \begin{align*}
  \int_0^t \int_{U_d} \bz \, \tN_0(\dd s, \dd\bz)
  &= \int_0^t \int_{V_1} g(\bX_{s-}, \br) \, N(\dd s, \dd\br)
    - \int_0^t \int_{\cR_0} \br \bbone_{\{\|\br\|<1\}} \, N(\dd s, \dd\br) \\
  &\quad
    - \int_0^t \int_{U_d} \br \bbone_{\{\|\br\|\geq1\}} \, \dd s \, \nu(\dd\br)
    - \sum_{j=1}^d
       \int_0^t X_{s,j} \, \dd s
       \int_{U_d} \bz \bbone_{\{\|\bz\|\geq1\}} \, \mu_j(\dd\bz) \\
  &\quad
    + \int_0^t \int_{\cR_0}
       \br \bbone_{\{\|\br\|<1\}} \, \tN(\dd s, \dd\br)
    + \int_0^t \int_{V_0} f(\bX_{s-}, \br) \, \tN(\dd s, \dd\br) \\
  &= \int_0^t \int_{V_0} f(\bX_{s-}, \br) \, \tN(\dd s, \dd\br)
     + \int_0^t \int_{V_1} g(\bX_{s-}, \br) \, N(\dd s, \dd\br) \\
  &\quad
     - \int_0^t \int_{U_d} \br \, \dd s \, \nu(\dd\br)
     - \sum_{j=1}^d
        \int_0^t X_{s,j} \, \dd s
        \int_{U_d} \bz \bbone_{\{\|\bz\|\geq1\}} \, \mu_j(\dd\bz) .
 \end{align*}
This proves that the SDE \eqref{SDE_X_YW_spec} holds \ $\ttPP$-almost surely,
 since
 \begin{equation*}
  \begin{aligned}
  &\int_0^t \bigl( \tBbeta + \tbB \bX_s \bigr) \dd s
   - \int_0^t \int_{U_d} \br \, \dd s \, \nu(\dd\br)
   - \sum_{j=1}^d
      \int_0^t X_{s,j} \, \dd s
      \int_{U_d} \bz \bbone_{\{\|\bz\|\geq1\}} \, \mu_j(\dd\bz) \\
  &= \tBbeta t + \tbB \int_0^t \bX_s \, \dd s
     - t \int_{U_d} \br \, \nu(\dd\br)
     - \sum_{j=1}^d
        \int_0^t X_{s,j} \, \dd s
        \int_{U_d} \bz \bbone_{\{\|\bz\|\geq1\}} \, \mu_j(\dd\bz) \\
  &= \left( \Bbeta + \int_{U_d} \br \, \nu(\dd\br) \right) t
     + \bD \int_0^t \bX_s \, \dd s
     + \sum_{j=1}^d
        \int_{U_d} \bz \bbone_{\{\|\bz\|\geq1\}} \, \mu_j(\dd\bz)
        \int_0^t X_{s,j} \, \dd s \\
  &\quad
     - t \int_{U_d} \br \, \nu(\dd\br)
     - \sum_{j=1}^d
      \int_0^t X_{s,j} \, \dd s
      \int_{U_d} \bz \bbone_{\{\|\bz\|\geq1\}} \, \mu_j(\dd\bz)
   = \int_0^t \bigl( \Bbeta + \bD \bX_s \bigr) \dd s .
  \end{aligned}
 \end{equation*}

The aim of the following discussion is to show that
 \ $\bigl( \ttOmega, \ttcF, (\ttcG_t)_{t\in\RR_+}, \ttPP, \bW, p, \bX \bigr)$
 \ is an \ $\RR_+^d$-valued weak solution to the SDE \eqref{SDE_X_YW_spec}.
Recall that the filtered probability space
 \ $(\ttOmega, \ttcF, (\ttcG_t)_{t\in\RR_+}, \ttPP)$ \ satisfies the usual
 hypotheses, and by Lemma \ref{ext_BM_PRN},
 \ $(W_{t,1}, \ldots, W_{t,d})_{t\in\RR_+}$ \ is a \ $d$-dimensional
 \ $(\ttcG_t)_{t\in\RR_+}$-Brownian motion, and \ $p$ \ is a stationary
 \ $(\ttcG_t)_{t\in\RR_+}$-Poisson point process on \ $V$ \ with characteristic
 measure \ $m$.
\ Since \ $(\bX_t)_{t\in\RR_+}$ \ is \ $\RR_+^d$-valued and has c\`adl\`ag
 trajectories on the original probability space \ $(\Omega, \cF, \PP)$, \ by
 the definition of an extension of a probability space (see Definition
 \ref{Def_extension}), the extended process (which is denoted by \ $\bX$ \ as
 well) on the extended probability space is \ $\RR_+^d$-valued and admits
 c\`adl\`ag trajectories as well.
By Remark \ref{Rem_extension}, the process \ $(\bX_t)_{t\in\RR_+}$ \ is
 \ $(\ttcG_t)_{t\in\RR_+}$-adapted, and clearly, the distribution of \ $\bX_0$
 \ is \ $n$.
\ Since \ $(\bX_t)_{t\in\RR_+}$ \ has c\`adl\`ag trajectories, (D4)(b) holds.
Since the process
 \ $\bigl(\int_0^t \int_{V_0}
           f(\bX_{s-}, \br) \, \tN(\dd s, \dd\br)\bigr)_{t\in\RR_+}$
 \ belongs to the space \ $\cM_2$, \ we have
\vspace*{-1mm}
 \begin{align*}
    &\ttEE\left( \left\| \int_0^t \int_{V_0}
                          f(\bX_{s-}, \br) \, \tN(\dd s, \dd\br)
                 \right\|^2 \right) \\
    &= \sum_{j=1}^d
        \ttEE\left( \int_0^t \int_{U_d} \int_{U_1}
                     \|\bz\|^2 \bbone_{\{\|\bz\|<1\}} \bbone_{\{u\leq X_{s,j}\}}
                     \, \dd s \, \mu_j(\dd\bz) \, \dd u \right)
     < \infty
\vspace*{-3mm}
 \end{align*}
 by Ikeda and Watanabe \cite[Chapter II, (3.9)]{IkeWat}, which yields (D4)(c).
We have already checked that (D4)(d) and (D4)(e) are satisfied.

Now we turn to prove the uniqueness in the sense of probability law for the
 SDE \eqref{SDE_X_YW_spec} among \ $\RR_+^d$-valued weak solutions.
If \ $\bigl( \Omega, \cF, (\cF_t)_{t\in\RR_+}, \PP, \bW, p, \bX \bigr)$ \ is an
 \ $\RR_+^d$-valued weak solution to the SDE \eqref{SDE_X_YW_spec}, then for
 each \ $G \in C_\cc^2(\RR, \RR)$ \ and
 \ $\bw = (w_1, \ldots, w_d)^\top \in\RR^d$, \ by It\^o's formula for
 \ $F(\bx) := G(\bw^\top \bx)$, \ $\bx = (x_1, \ldots, x_d)^\top \in \RR^d$,
 \ with \ $\partial_{x_k} F(\bx) = G'(\bw^\top \bx) w_k$,
 \ $\partial_{x_k} \partial_{x_\ell} F(\bx) = G''(\bw^\top \bx) w_k w_\ell$,
 \ $k, \ell \in \{1, \ldots, d\}$, \ we have
\vspace*{-1mm}
 \[
   G(\bw^\top \bX_t) = G(\bw^\top \bX_0) + \sum_{\ell=1}^6 I_\ell(t) , \qquad
   t \in \RR_+ ,
\vspace*{-1mm}
 \]
 where
\vspace*{-1mm}
 \begin{align*}
  I_1(t)
  &:= \int_0^t
       G'(\bw^\top \bX_s) \bw^\top (\Bbeta + \bD \bX_s) \, \dd s , \\
  I_2(t)
  &:= \sum_{j=1}^d
       \int_0^t w_j G'(\bw^\top \bX_s) \sqrt{2 c_j X_{s,j}} \, \dd W_{s,j} , \\
  I_3(t)
  &:= \sum_{j=1}^d \int_0^t w_j^2 G''(\bw^\top \bX_s) c_j X_{s,j} \, \dd s , \\
  I_4(t)
  &:= \int_0^t \int_{V_0}
       \bigl[ G(\bw^\top \bX_{s-} + \bw^\top f(\bX_{s-}, \br))
              - G(\bw^\top \bX_{s-}) \bigr]
       \, \tN(\dd s, \dd\br) , \\
  I_5(t)
  &:= \int_0^t \int_{V_0}
       \bigl[ G(\bw^\top \bX_s + \bw^\top f(\bX_{s-}, \br))
              - G(\bw^\top \bX_s) \\
  &\phantom{:= \int_0^t \int_{\cR_j} \bigl(}
              - G'(\bw^\top \bX_s) \bw^\top f(\bX_{s-}, \br) \bigr]
        \, \dd s \, m(\dd\br) ,
 \end{align*}
 \begin{align*}
  I_6(t)
  &:= \int_0^t \int_{V_1}
       \bigl[ G(\bw^\top \bX_{s-} + \bw^\top g(\bX_{s-}, \br))
              - G(\bw^\top \bX_{s-}) \bigr]
       \, N(\dd s, \dd\br) .
 \end{align*}
The last integral can be written as \ $I_6(t) = I_{6,1}(t) + I_{6,2}(t)$,
 \ where
 \begin{align*}
  I_{6,1}(t)
  &:= \int_0^t \int_{V_1}
       \bigl[ G(\bw^\top \bX_{s-} + \bw^\top g(\bX_{s-}, \br))
              - G(\bw^\top \bX_{s-}) \bigr]
       \, \tN(\dd s, \dd\br) , \\
  I_{6,2}(t)
  &:= \int_0^t \int_{V_1}
       \bigl[ G(\bw^\top \bX_s + \bw^\top g(\bX_s, \br))
              - G(\bw^\top \bX_s) \bigr]
       \,\dd s \, m(\dd\br) ,
 \end{align*}
 since
 \begin{align}\label{help_int61}
 \begin{split}
  &\EE\left( \int_0^t \int_{V_1}
              \bigl| G(\bw^\top \bX_{s-} + \bw^\top g(\bX_{s-}, \br))
                     - G(\bw^\top \bX_{s-}) \bigr|
              \, \dd s \, m(\dd\br) \right) \\
  &= \EE\left( \int_0^t \int_{U_d}
              \bigl| G(\bw^\top \bX_s + \bw^\top \br)
                     - G(\bw^\top \bX_s) \bigr|
              \, \dd s \, \nu(\dd\br) \right) \\
  &\quad
     + \sum_{j=1}^d
        \EE\left( \int_0^t \int_{U_d} \int_{U_1}
                   \bigl| G(\bw^\top \bX_s + \bw^\top \bz \bbone_{\{u\leq X_{s,j}\}})
                          - G(\bw^\top \bX_s) \bigr| \bbone_{\{\|\bz\|\geq1\}}
                   \, \dd s \, \mu_j(\dd\bz) \, \dd u \right)\\
  &< \infty ,
  \end{split}
 \end{align}
 i.e., for all \ $\bw\in\RR^d$, \ the function
 \ $\RR_+\times V_1\times\Omega\ni(s,\br,\omega)\mapsto
  G(\bw^\top \bX_{s-}(\omega) + \bw^\top g(\bX_{s-}(\omega), \br)) - G(\bw^\top \bX_{s-}(\omega))$
 \ belongs to the class \ $\bF_p^1$.
\ Indeed, by mean value theorem and \eqref{help4}, there exists some
 \ $\theta_0 = \theta_0(\bw, \bX_s, \br) \in [0, 1]$ \ such that
 \begin{align*}
  &\EE\left( \int_0^t \int_{U_d}
              \bigl| G(\bw^\top \bX_s + \bw^\top \br)
                     - G(\bw^\top \bX_s) \bigr|
              \, \dd s \, \nu(\dd\br) \right) \\
  & = \EE\left( \int_0^t \int_{U_d}
                 \bigl| G'(\bw^\top \bX_s + \theta_0 \bw^\top \br) \bigr|
                 |\bw^\top \br|
                 \, \dd s \, \nu(\dd\br) \right)
   \leq \|\bw\| \sup_{x\in\RR} |G'(x)| \int_{U_d} \|\br\| \, \nu(\dd\br)
    < \infty
 \end{align*}
 due to that \ $G'$ \ is bounded.
In a similar way, there exists some
 \ $\theta = \theta(\bw, \bX_s, \bz) \in [0, 1]$ \ such that for each
 \ $j \in \{1, \ldots, d\}$,
 \begin{align*}
  &\EE\left( \int_0^t \int_{U_d} \int_{U_1}
              \bigl| G(\bw^\top \bX_s + \bw^\top \bz \bbone_{\{u\leq X_{s,j}\}})
                     - G(\bw^\top \bX_s) \bigr|
              \bbone_{\{\|\bz\|\geq1\}}
              \, \dd s \, \mu_j(\dd\bz) \, \dd u \right) \\
  &\qquad
   =\EE\left( \int_0^t \int_{U_d} \int_{U_1}
              \bigl| G(\bw^\top \bX_s + \bw^\top \bz)
                     - G(\bw^\top \bX_s) \bigr|
              \bbone_{\{u\leq X_{s,j}\}} \bbone_{\{\|\bz\|\geq1\}}
              \, \dd s \, \mu_j(\dd\bz) \, \dd u \right) \\
  &\qquad
   = \EE\left( \int_0^t \int_{U_d} \int_{U_1}
                \bigl| G'(\bw^\top \bX_s
                          + \theta \bw^\top \bz ) \bigr|
                |\bw^\top \bz| \bbone_{\{u\leq X_{s,j}\}} \bbone_{\{\|\bz\|\geq1\}}
                 \, \dd s \, \mu_j(\dd\br) \, \dd u \right) \\
  &\qquad
   \leq \|\bw\| \sup_{x\in\RR} |G'(x)| \int_0^t \EE(X_{s,j}) \, \dd s
        \int_{U_d} \|\bz\| \bbone_{\{\|\bz\|\geq1\}} \, \mu_j(\dd\bz)
    < \infty
 \end{align*}
 due to that \ $G'$ \ is bounded, Lemma \ref{moment_1} (which can be applied
 since \ $\int_{\RR_+^d} \|\bz\| \, n(\dd\bz) < \infty$) \ and the moment
 condition \eqref{help_Li_page45_equiv}.

In what follows, we identify some of these integrals with some terms in part
 (5) of Theorem 9.18 of Li \cite{Li}.
We have
 \begin{align*}
  I_1(t) &= \int_0^t G'(\bw^\top \bX_s) \bw^\top \Bbeta \, \dd s
            + \int_0^t G'(\bw^\top \bX_s) \bw^\top \tbB \bX_s \, \dd s \\
         &\quad
            - \sum_{i=1}^d \sum_{j=1}^d
               \int_0^t G'(\bw^\top \bX_s) w_i X_{s,j} \, \dd s
               \int_{U_d} z_i \bbone_{\{ \|\bz\|\geq1 \}} \, \mu_j(\dd\bz) ,
 \end{align*}
 where the first two terms on the right hand side can be identified with
 \ $\int_0^t G'(Y_s(f)) \eta(f) \, \dd s$ \ and
 \ $\int_0^t G'(Y_s(f)) Y_s(Af+\gamma f - bf) \, \dd s$
 \ (see, \eqref{Agammab}).
The sum of the third term on the right hand side and \ $I_{6,2}(t) + I_5(t)$
 \ can be written in the form
 \begin{align*}
  &- \sum_{j=1}^d
      \int_0^t \int_{U_d} \int_{U_1}
       G'(\bw^\top \bX_s) \bw^\top \bz \bbone_{\{\|\bz\|\geq1\}}
       \bbone_{\{u\leq X_{s,j}\}}
       \, \dd s \, \mu_j(\dd\bz) \, \dd u \\
  &+ \int_0^t \int_{U_d}
      \bigl[ G(\bw^\top \bX_s + \bw^\top \br) - G(\bw^\top \bX_s) \bigr]
      \, \dd s \, \nu(\dd\br) \\
  &+ \sum_{j=1}^d
      \int_0^t \int_{U_d} \int_{U_1}
       \bigl[ G(\bw^\top \bX_s + \bw^\top \bz \bbone_{\{u\leq X_{s,j}\}})
              - G(\bw^\top \bX_s) \bigr]
       \bbone_{\{\|\bz\|\geq1\}}
       \, \dd s \, \mu_j(\dd\bz) \, \dd u \\
  &+ \sum_{j=1}^d
      \int_0^t \int_{U_d} \int_{U_1}
       \bigl[ G(\bw^\top \bX_s + \bw^\top \bz \bbone_{\{u\leq X_{s,j}\}})
              - G(\bw^\top \bX_s) \\
  &\phantom{+ \sum_{j=1}^d \int_0^t \int_{U_d} \int_{U_1} \bigl[}
              - G'(\bw^\top \bX_s) \bw^\top \bz \bbone_{\{u\leq X_{s,j}\}} \bigr]
       \bbone_{\{\|\bz\|<1\}}
       \, \dd s \, \mu_j(\dd\bz) \, \dd u \\
  &= \int_0^t \int_{U_d}
      \bigl[ G(\bw^\top \bX_s + \bw^\top \br) - G(\bw^\top \bX_s) \bigr]
      \, \dd s \, \nu(\dd\br) \\
  &\quad
     + \sum_{j=1}^d
        \int_0^t \int_{U_d} \int_{U_1}
         \bigl[ G(\bw^\top \bX_s + \bw^\top \bz \bbone_{\{u\leq X_{s,j}\}})
                - G(\bw^\top \bX_s) \\
  &\phantom{\quad + \sum_{j=1}^d \int_0^t \int_{U_d} \int_{U_1} \bigl[}
                - G'(\bw^\top \bX_s) \bw^\top \bz \bbone_{\{u\leq X_{s,j}\}} \bigr]
       \, \dd s \, \mu_j(\dd\bz) \, \dd u ,
 \end{align*}
 which can be identified with
 \begin{align*}
  &\int_0^t \int_{M(E)^\circ}
       [G(Y_s(f) + \kappa(f)) - G(Y_s(f))] \, H(\dd\kappa) \, \dd s \\
  &+ \int_0^t \int_E Y_s(\dd x) \int_{M(E)^\circ}
      [G(Y_s(f) + \kappa(f)) - G(Y_s(f)) - \kappa(f) G'(Y_s(f))]
      \, H(x, \dd\kappa) \, \dd s .
 \end{align*}
The integral \ $I_3(t)$ \ can be identified with
 \ $\int_0^t G''(Y_s(f)) Y_s(cf^2) \, \dd s$.

Next we show that the process \ $(I_2(t) + I_4(t) + I_{6,1}(t))_{t\in\RR_+}$ \ is
 a continuous local martingale.
Since \ $G'$ \ is bounded and \ $\bX$ \ has c\`adl\`ag trajectories, we have
 \ $\PP(\int_0^t w_j^2 G'(\bw^\top \bX_s)^2 \, 2 c_j X_{s,j} \, \dd s < \infty)
    = 1$ \ for all
 \ $t \in \RR_+$ \ and \ $j \in \{1, \ldots, d\}$, \ hence
 \ $(I_2(t))_{t\in\RR_+}$ \ is a continuous local martingale
 (see, e.g., Karatzas and Shreve \cite[Definition 3.2.23]{KarShr}).
In order to prove that \ $(I_4(t))_{t\in\RR_+}$ \ is a martingale, by page 62 in
 Ikeda and Watanabe \cite{IkeWat}, it is enough to check that
 \[
   \EE\left( \int_0^t \int_{V_0}
              \bigl| G(\bw^\top \bX_s + \bw^\top f(\bX_s, \br))
                     - G(\bw^\top \bX_s) \bigr|^2
              \,\dd s \, m(\dd\br) \right)
   < \infty .
 \]
By mean value theorem, there exists some
 \ $\vartheta_0 = \vartheta_0(\bw, \bX_s, \bz) \in [0, 1]$ \ such that for
 each \ $j \in \{1, \ldots, d\}$,
 \begin{align*}
  &\EE\left( \int_0^t \int_{U_d} \int_{U_1}
              \bigl| G(\bw^\top \bX_s + \bw^\top \bz \bbone_{\{u\leq X_{s,j}\}})
                     - G(\bw^\top \bX_s) \bigr|^2
              \bbone_{\{\|\bz\|<1\}}
              \, \dd s \, \mu_j(\dd\bz) \, \dd u \right) \\
  &\qquad
   = \EE\left( \int_0^t \int_{U_d} \int_{U_1}
                \bigl| G'(\bw^\top \bX_s
                          + \vartheta_0 \bw^\top \bz) \bigr|^2
                (\bw^\top \bz)^2 \bbone_{\{u\leq X_{s,j}\}}
                \bbone_{\{\|\bz\|<1\}}
                \, \dd s \, \mu_j(\dd\bz) \, \dd u \right) \\
  &\qquad
   \leq \|\bw\|^2 \sup_{x\in\RR} |G'(x)|^2 \int_0^t \EE(X_{s,j}) \, \dd s
        \int_{U_d} \|\bz\|^2 \bbone_{\{\|\bz\|<1\}} \, \mu_j(\dd\bz)
   < \infty
 \end{align*}
 due to that \ $G'$ \ is bounded, Lemma \ref{moment_1} and \eqref{z2}.
Hence \ $(I_4(t))_{t\in\RR_+}$ \ is a martingale.
Further, by \eqref{help_int61} and page 62 in Ikeda and Watanabe
 \cite{IkeWat}, we get \ $(I_{6,1}(t))_{t\in\RR_+}$ \ is a martingale.
Consequently, by Theorem 9.18 of Li \cite{Li}, \ $(\bX_t)_{t\in\RR_+}$ \ is a CBI
 process with parameters \ $(d, \bc, \Bbeta, \bB, \nu, \bmu)$.
\ This yields the uniqueness in the sense of probability law for the SDE
 \eqref{SDE_X_YW_spec} among \ $\RR_+^d$-valued weak solutions, and that any
 \ $\RR_+^d$-valued weak solution is a CBI
 process with parameters \ $(d, \bc, \Bbeta, \bB, \nu, \bmu)$ \ as well.
\proofend

\section{Multi-type CBI process as a strong solution of an SDE}
\label{section_strong}

\begin{Def}\label{Def_pathwise_uniqueness}
We say that pathwise uniqueness holds for the SDE \eqref{SDE_X_YW_spec} among
 \ $\RR_+^d$-valued weak solutions if
 whenever \ $\bigl( \Omega, \cF, (\cF_t)_{t\in\RR_+}, \PP, \bW, p, \bX \bigr)$
 \ and \ $\bigl( \Omega, \cF, (\cF_t)_{t\in\RR_+}, \PP, \bW, p, \tbX \bigr)$
 \ are \ $\RR_+^d$-valued weak solutions of the SDE \eqref{SDE_X_YW_spec} such
 that \ $\PP(\bX_0 = \tbX_0) = 1$, \ then
 \ $\PP(\text{$\bX_t = \tbX_t$ \ for all \ $t \in \RR_+$}) = 1$.
\end{Def}

Next we prove a comparison theorem for the SDE \eqref{SDE_X_YW_spec} in
 \ $\Bbeta$.

\begin{Lem}\label{comparison_beta}
Let \ $(d, \bc, \Bbeta, \bB, \nu, \bmu)$ \ be a set of admissible parameters
 in the sense of Definition \ref{Def_admissible} such that the moment
 condition \eqref{moment_condition_1} holds.
Suppose that \ $\Bbeta' \in \RR_+^d$ \ with \ $\Bbeta \leq \Bbeta'$.
\ Let \ $\bigl( \Omega, \cF , (\cF_t)_{t\in\RR_+}, \PP, \bW, p, \bX \bigr)$ \ and
 \ $\bigl( \Omega, \cF , (\cF_t)_{t\in\RR_+}, \PP, \bW, p, \bX' \bigr)$ \ be
 \ $\RR_+^d$-valued weak solutions of the SDE \eqref{SDE_X_YW_spec} with
 \ $\Bbeta$ \ and \ $\Bbeta'$, \ respectively.
Then \ $\PP(\bX_0 \leq \bX_0') = 1$ \ implies
 \ $\PP(\text{$\bX_t \leq \bX_t'$ \ for all \ $t \in \RR_+$}) = 1$.
\ Particularly, pathwise uniqueness holds for the SDE \eqref{SDE_X_YW_spec}
 among \ $\RR_+^d$-valued weak solutions.
\end{Lem}

\noindent
\textbf{Proof.} \
We follow the ideas of the proof of Theorem 3.1 of Ma \cite{Ma}, which is an
 adaptation of that of Theorem 5.5 of Fu and Li \cite{FuLi}.
There is a sequence \ $\phi_k : \RR \to \RR_+$, \ $k \in \NN$, \ of twice
 continuously differentiable functions such that
 \begin{enumerate}
  \item[(i)]
   $\phi_k(z) \uparrow z^+$ \ as \ $k \to \infty$;
  \item[(ii)]
   $\phi_k'(z) \in [0, 1]$ \ for all \ $z \in \RR_+$ \ and \ $k \in \NN$;
  \item[(iii)]
   $\phi_k'(z) = \phi_k(z) = 0$ \ whenever \ $-z \in \RR_+$ \ and
    \ $k \in \NN$;
  \item[(iv)]
   $\phi_k''(x - y) (\sqrt{x} - \sqrt{y})^2 \leq 2/k$ \ for all
    \ $x, y \in \RR_+$ \ and \ $k \in \NN$.
 \end{enumerate}
For a construction of such functions, see, e.g., the proof of Theorem 3.1 of
 Ma \cite{Ma}.
Let \ $\bY_t := \bX_t - \bX_t'$ \ for all \ $t \in \RR_+$.
\ By \eqref{SDE_X_YW_spec}, and using that
 \begin{align*}
  \int_0^t\int_{\cR_0} g(\bX_{s-},\br)\, N(\dd s,\dd \br)
     = \int_0^t\int_{\cR_0} \br \, N(\dd s,\dd \br)
     = \int_0^t\int_{\cR_0} g(\bX_{s-}',\br)\, N(\dd s,\dd \br),
 \end{align*}
 we have
 \begin{align*}
  Y_{t,i}
  &= Y_{0,i}
     + \int_0^t
        \Bigl( \beta_i - \beta_i' + \be_i^\top \bD \bY_s \Bigr)
        \, \dd s
     + \int_0^t
        \sqrt{2c_i} \Bigl(\sqrt{X_{s,i}} - \sqrt{X_{s,i}'}\Bigr) \, \dd W_{s,i} \\
   &\quad
     + \sum_{j=1}^d
        \int_0^t \int_{\cR_{j,0}}
         \bigl(\bbone_{\{u\leq X_{s-,j}\}} - \bbone_{\{u\leq X_{s-,j}'\}}\bigr)
         z_i \bbone_{\{\|\bz\|<1\}}
         \, \tN(\dd s, \dd\br) \\
   &\quad
     + \sum_{j=1}^d
        \int_0^t \int_{\cR_{j,1}}
         \bigl(\bbone_{\{u\leq X_{s-,j}\}} - \bbone_{\{u\leq X_{s-,j}'\}}\bigr)
         z_i \bbone_{\{\|\bz\|\geq1\}}
         \, N(\dd s, \dd\br)
 \end{align*}
 for all \ $t \in \RR_+$ \ and \ $i \in \{1, \ldots, d\}$.
\ For each \ $m \in \NN$, \ put
 \[
   \tau_m := \inf\Bigl\{ t \in \RR_+
                         : \max_{i \in \{1, \ldots, d\}} \max\{X_{t,i}, X_{t,i}'\}
                           \geq m \Bigr\} .
 \]
By It\^o's formula (which can be used since \ $\bX$ \ and \ $\bX'$ \ are
 adapted to the same filtration \ $(\cF_t)_{t\in\RR_+}$), \ we obtain
 \[
   \phi_k(Y_{t\land\tau_m,i}) = \phi_k(Y_{0,i}) + \sum_{\ell=1}^6 I_{i,m,k,\ell}(t)
 \]
 for all \ $t \in \RR_+$, \ $i \in \{1, \ldots, d\}$ \ and \ $k, m \in \NN$,
 \ where
 \begin{align*}
  &I_{i,m,k,1}(t)
   := \int_0^{t\land\tau_m}
       \phi_k'(Y_{s,i})
       \Bigl( \beta_i - \beta_i' + \be_i^\top \bD \bY_s \Bigr)
       \, \dd s , \\
  &I_{i,m,k,2}(t)
   := \int_0^{t\land\tau_m}
       \phi_k'(Y_{s,i}) \sqrt{2c_i} \Bigl(\sqrt{X_{s,i}} - \sqrt{X_{s,i}'}\Bigr)
       \, \dd W_{s,i} , \\
  &I_{i,m,k,3}(t)
   := \frac{1}{2}
      \int_0^{t\land\tau_m}
       \phi_k''(Y_{s,i}) 2 c_i \Bigl(\sqrt{X_{s,i}} - \sqrt{X_{s,i}'}\Bigr)^2
       \, \dd s ,\\
  &I_{i,m,k,4}(t)
   := \sum_{j=1}^d
       \int_0^{t\land\tau_m} \int_{\cR_{j,0}}
        \Bigl[ \phi_k\bigl( Y_{s-,i}
                            + (\bbone_{\{u\leq X_{s-,j}\}}
                               - \bbone_{\{u\leq X_{s-,j}'\}}) z_i \bigr)
               - \phi_k(Y_{s-,i}) \Bigr]
        \bbone_{\{\|\bz\|<1\}}
        \tN(\dd s, \dd\br) \\
  &= \sum_{j=1}^d
       \int_0^{t\land\tau_m} \int_{\cR_{j,0}}
        \Bigl[ \phi_k(Y_{s-,i} + z_i) - \phi_k(Y_{s-,i}) \Bigr]
        \bbone_{\{\|\bz\|<1\}}
        \bbone_{\{X_{s-,j}'<u\leq X_{s-,j}\}} \bbone_{\{Y_{s-,j}>0\}}
        \, \tN(\dd s, \dd\br) \\
  &\quad
      + \sum_{j=1}^d
         \int_0^{t\land\tau_m} \int_{\cR_{j,0}}
          \Bigl[ \phi_k(Y_{s-,i} - z_i) - \phi_k(Y_{s-,i}) \Bigr]
          \bbone_{\{\|\bz\|<1\}}
          \bbone_{\{X_{s-,j}<u\leq X_{s-,j}'\}} \bbone_{\{Y_{s-,j}<0\}}
          \, \tN(\dd s, \dd\br) ,
 \end{align*}
 \begin{align*}
  &I_{i,m,k,5}(t)
   := \sum_{j=1}^d
       \int_0^{t\land\tau_m} \int_{U_d} \int_{U_1}
        \Bigl[ \phi_k\bigl(Y_{s-,i}
                           + (\bbone_{\{u\leq X_{s-,j}\}}
                              - \bbone_{\{u\leq X_{s-,j}'\}}) z_i \bigr)
               - \phi_k(Y_{s-,i}) \\
  &\phantom{I_{i,m,k,5}(t)
            := \sum_{j=1}^d \int_0^{t\land\tau_m} \int_{U_d} \int_{U_1} \Bigl[}
               - \phi_k'(Y_{s-,i})
                 (\bbone_{\{u\leq X_{s-,j}\}} - \bbone_{\{u\leq X_{s-,j}'\}}) z_i \Bigr]
        \bbone_{\{\|\bz\|<1\}}
        \, \dd s \, \mu_j(\dd\bz) \, \dd u \\
  &= \sum_{j=1}^d
      \int_0^{t\land\tau_m} \int_{U_d} \int_{U_1}
       \Bigl[ \phi_k(Y_{s-,i} + z_i) - \phi_k(Y_{s-,i})
              - \phi_k'(Y_{s-,i}) z_i \Bigr]\\
  &\phantom{= \sum_{j=1}^d \int_0^{t\land\tau_m} \int_{U_d} \int_{U_1}\; }
        \;\times\bbone_{\{\|\bz\|<1\}}
        \bbone_{\{X_{s-,j}'<u\leq X_{s-,j}\}} \bbone_{\{Y_{s-,j}>0\}}
        \, \dd s \, \mu_j(\dd\bz) \, \dd u \\
  &\quad
      + \sum_{j=1}^d
         \int_0^{t\land\tau_m} \int_{U_d} \int_{U_1}
          \Bigl[ \phi_k(Y_{s-,i} - z_i) - \phi_k(Y_{s-,i})
                 + \phi_k'(Y_{s-,i}) z_i \Bigr]\\
  &\phantom{= +\sum_{j=1}^d \int_0^{t\land\tau_m} \int_{U_d} \int_{U_1}\; }
          \;\times\bbone_{\{\|\bz\|<1\}}
          \bbone_{\{X_{s-,j}<u\leq X_{s-,j}'\}} \bbone_{\{Y_{s-,j}<0\}}
          \, \dd s \, \mu_j(\dd\bz) \, \dd u , \\
  &I_{i,m,k,6}(t)
   := \sum_{j=1}^d
       \int_0^{t\land\tau_m} \int_{\cR_{j,1}}
        \Bigl[ \phi_k\bigl( Y_{s-,i}
                            + (\bbone_{\{u\leq X_{s-,j}\}}
                               - \bbone_{\{u\leq X_{s-,j}'\}}) z_i \bigr)
               - \phi_k(Y_{s-,i}) \Bigr]
        \bbone_{\{\|\bz\|\geq1\}}
        N(\dd s, \dd\br) ,
 \end{align*}
 where we used that
 \begin{align}\label{seged}
  \bbone_{\{u\leq X_{s-,j}\}} - \bbone_{\{u\leq X_{s-,j}'\}}
  = \begin{cases}
     1  & \text{if \ $Y_{s-,j} > 0$ \ and \ $X_{s-,j}' < u \leq X_{s-,j}$,} \\
     -1 & \text{if \ $Y_{s-,j} < 0$ \ and \ $X_{s-,j} < u \leq X_{s-,j}'$,} \\
     0  & \text{otherwise.}
    \end{cases}
 \end{align}
Using formula (3.8) in Chapter II in Ikeda and Watanabe \cite{IkeWat},
 the last integral can be written as
 \ $I_{i,m,k,6}(t) = I_{i,m,k,6,1}(t) + I_{i,m,k,6,2}(t)$, \ where
 \[
   I_{i,m,k,6,1}(t)
   := \sum_{j=1}^d
       \int_0^{t\land\tau_m} \int_{\cR_{j,1}}
        \Bigl[ \phi_k\bigl( Y_{s-,i}
                            + (\bbone_{\{u\leq X_{s-,j}\}}
                               - \bbone_{\{u\leq X_{s-,j}'\}}) z_i \bigr)
               - \phi_k(Y_{s-,i}) \Bigr]
        \bbone_{\{\|\bz\|\geq1\}}
        \tN(\dd s, \dd\br)
 \]
 \[
   I_{i,m,k,6,2}(t)
   := \sum_{j=1}^d
       \int_0^{t\land\tau_m} \int_{U_d} \int_{U_1}
        \Bigl[ \phi_k\bigl( Y_{s-,i}
                            + (\bbone_{\{u\leq X_{s-,j}\}}
                               - \bbone_{\{u\leq X_{s-,j}'\}}) z_i \bigr)
               - \phi_k(Y_{s-,i}) \Bigr]
        \bbone_{\{\|\bz\|\geq1\}}
        \, \dd s \, \mu_j(\dd\bz) \, \dd u ,
 \]
 since, for each \ $j \in \{1, \ldots, d\}$,
 \begin{align*}
  &\EE\left( \int_0^{t\land\tau_m} \int_{U_d} \int_{U_1}
              \Bigl| \phi_k\bigl( Y_{s-,i}
                                  + (\bbone_{\{u\leq X_{s-,j}\}}
                                  - \bbone_{\{u\leq X_{s-,j}'\}}) z_i \bigr)
                     - \phi_k(Y_{s-,i}) \Bigr|
              \bbone_{\{\|\bz\|\geq1\}}
              \, \dd s \, \mu_j(\dd\bz) \, \dd u \right) \\
  &=\EE\left( \int_0^{t\land\tau_m} \int_{U_d} \int_{U_1}
               \Bigl| \phi_k(Y_{s-,i} + z_i) - \phi_k(Y_{s-,i}) \Bigr|
               \bbone_{\{\|\bz\|\geq1\}} \bbone_{\{X_{s-,j}'<u\leq X_{s-,j}\}}
               \bbone_{\{Y_{s-,j}>0\}}
               \, \dd s \, \mu_j(\dd\bz) \, \dd u \right) \\
  &\quad
     + \EE\left( \int_0^{t\land\tau_m} \int_{U_d} \int_{U_1}
                  \Bigl| \phi_k(Y_{s-,i} - z_i) - \phi_k(Y_{s-,i}) \Bigr|
                  \bbone_{\{\|\bz\|\geq1\}} \bbone_{\{X_{s-,j}<u\leq X_{s-,j}'\}}
                  \bbone_{\{Y_{s-,j}<0\}}
                  \, \dd s \, \mu_j(\dd\bz) \, \dd u \right)
  \end{align*}
  \begin{align*}
  &\leq\EE\left( \int_0^{t\land\tau_m} \int_{U_d}
                  z_i \bbone_{\{\|\bz\|\geq1\}} |Y_{s-,j}|
                  \, \dd s \, \mu_j(\dd\bz) \right)
   \leq 2 m t \int_{U_d} z_i \bbone_{\{\|\bz\|\geq1\}} \, \mu_j(\dd\bz)
   < \infty ,
 \end{align*}
 where we used that, by properties (ii) and (iii) of the function \ $\phi_k$,
 \ we have \ $\phi_k'(u) \in [0, 1]$ \ for all \ $u \in \RR$, \ and hence, by
 mean value theorem,
 \begin{equation}\label{(ii)}
  - z \leq \phi_k(y - z) - \phi_k(y) \leq 0 \leq \phi_k(y + z) - \phi_k(y)
      \leq z ,
  \qquad y \in \RR , \quad z \in \RR_+ , \quad k \in \NN .
 \end{equation}

One can check that the process
 \ $\left(I_{i,m,k,2}(t) + I_{i,m,k,4}(t)
          + I_{i,m,k,6,1}(t)\right)_{t\in\RR_+}$
 \ is a martingale.
Indeed, by properties (ii) and (iii) of the function \ $\phi_k$ \ and the
 definition of \ $\tau_m$,
 \begin{align*}
  \EE\left( \int_0^{t\land\tau_m}
             \left(\phi_k'(Y_{s,i}) \sqrt{2c_i}
                   \Bigl(\sqrt{X_{s,i}} - \sqrt{X_{s,i}'}\Bigr)\right)^2
             \, \dd s  \right)
   &\leq 2c_i \EE\left( \int_0^{t\land\tau_m}
                         (X_{s,i} + X_{s,i}') \, \dd s \right) \\
   &\leq 4 c_i m t < \infty ,
 \end{align*}
 hence, by Ikeda and Watanabe \cite[page 55]{IkeWat},
  \ $\left(I_{i,m,k,2}(t)\right)_{t\in\RR_+}$ \ is a martingale.
Next we show
 \begin{align*}
  \EE\left( \int_0^{t\land\tau_m} \int_{U_d} \int_{U_1}
             \left| \phi_k(Y_{s-,i} + z_i) - \phi_k(Y_{s-,i}) \right|^2
             \bbone_{\{\|\bz\|<1\}}
             \bbone_{\{X_{s-,j}'<u\leq X_{s-,j}\}} \bbone_{\{Y_{s-,j}>0\}}
             \, \dd s \, \mu_j(\dd \bz) \, \dd u \right)
  < \infty ,
 \end{align*}
 and
 \begin{align*}
  \EE\left( \int_0^{t\land\tau_m} \int_{U_d} \int_{U_1}
             \left| \phi_k(Y_{s-,i} - z_i) - \phi_k(Y_{s-,i}) \right|^2
             \bbone_{\{\|\bz\|<1\}}
             \bbone_{\{X_{s-,j}<u\leq X_{s-,j}'\}} \bbone_{\{Y_{s-,j}<0\}}
             \, \dd s \, \mu_j(\dd \bz) \, \dd u \right)
  < \infty
 \end{align*}
 for all \ $j \in \{1, \ldots, d\}$, \ which yield that the functions
 \begin{align*}
  \RR_+\times U_d\times U_1\times \Omega\ni (s,\bz,u,\omega)
      \mapsto & (\phi_k(Y_{s-,i}(\omega) + z_i) - \phi_k(Y_{s-,i}(\omega)))
              \bbone_{\{\|\bz\|<1\}} \\
             & \times\bbone_{\{X_{s-,j}'(\omega)<u\leq X_{s-,j}(\omega)\}} \bbone_{\{Y_{s-,j}(\omega)>0\}}
               \bbone_{\{s\leq \tau_m(\omega)\}}
 \end{align*}
 and
 \begin{align*}
   \RR_+\times U_d\times U_1\times \Omega\ni (s,\bz,u,\omega)
      \mapsto & (\phi_k(Y_{s-,i}(\omega) - z_i) - \phi_k(Y_{s-,i}(\omega)))
               \bbone_{\{\|\bz\|<1\}} \\
              & \times\bbone_{\{X_{s-,j}(\omega)<u\leq X_{s-,j}'(\omega)\}} \bbone_{\{Y_{s-,j}(\omega)<0\}}
                \bbone_{\{s\leq \tau_m(\omega)\}}
 \end{align*}
 belong to the class \ $\bF_p^2$, \ and then
 \ $\left(I_{i,m,k,4}(t)\right)_{t\in\RR_+}$ \ is a martingale, again by page 62 in Ikeda and Watanabe \cite{IkeWat}.
By \eqref{z2} and \eqref{(ii)},
 \begin{multline*}
  \EE\left( \int_0^{t\land\tau_m} \int_{U_d} \int_{U_1}
             \left| \phi_k(Y_{s-,i} + z_i) - \phi_k(Y_{s-,i}) \right|^2
             \bbone_{\{\|\bz\|<1\}}
             \bbone_{\{X_{s-,j}'<u\leq X_{s-,j}\}} \bbone_{\{Y_{s-,j}>0\}}
             \, \dd s \, \mu_j(\dd \bz) \, \dd u \right) \\
  \begin{aligned}
   &\leq \EE\left( \int_0^{t\land\tau_m} \int_{U_d} \int_{U_1}
                    z_i^2 \bbone_{\{\|\bz\|<1\}} \bbone_{\{X_{s-,j}'<u\leq X_{s-,j}\}}
                    \bbone_{\{Y_{s-,j}>0\}}
                    \, \dd s \, \mu_j(\dd \bz) \, \dd u \right) \\
   &= \EE\left( \int_0^{t\land\tau_m} \int_{U_d}
                 z_i^2 \bbone_{\{\|\bz\|<1\}} Y_{s-,j} \bbone_{\{Y_{s-,j}>0\}}
                 \, \dd s \, \mu_j(\dd \bz) \right)
    \leq 2mt\int_{U_d} z_i^2 \bbone_{\{\|\bz\|<1\}} \, \mu_j(\dd\bz)
    < \infty .
  \end{aligned}
 \end{multline*}
In the same way one can get the finiteness of the other expectation.
Finally, we show
 \begin{align*}
  \EE\left( \int_0^{t\land\tau_m} \int_{U_d} \int_{U_1}
             \left| \phi_k(Y_{s-,i} + z_i) - \phi_k(Y_{s-,i}) \right|
             \bbone_{\{\|\bz\|\geq1\}}
             \bbone_{\{X_{s-,j}'<u\leq X_{s-,j}\}} \bbone_{\{Y_{s-,j}>0\}}
             \, \dd s \, \mu_j(\dd \bz) \, \dd u \right)
  < \infty ,
 \end{align*}
 and
 \begin{align*}
  \EE\left( \int_0^{t\land\tau_m} \int_{U_d} \int_{U_1}
             \left| \phi_k(Y_{s-,i} - z_i) - \phi_k(Y_{s-,i}) \right|
             \bbone_{\{\|\bz\|\geq1\}}
             \bbone_{\{X_{s-,j}<u\leq X_{s-,j}'\}} \bbone_{\{Y_{s-,j}<0\}}
             \, \dd s \, \mu_j(\dd \bz) \, \dd u \right)
  < \infty
 \end{align*}
 for all \ $j \in \{1, \ldots, d\}$, \ which yield that the functions
 \begin{align*}
    \RR_+\times U_d\times U_1\times \Omega\ni (s,\bz,u,\omega)
      \mapsto &  ( \phi_k(Y_{s-,i}(\omega) + z_i) - \phi_k(Y_{s-,i}(\omega)) )
                 \bbone_{\{\|\bz\|\geq1\}}\\
             &\times\bbone_{\{X_{s-,j}'(\omega)<u\leq X_{s-,j}(\omega)\}} \bbone_{\{Y_{s-,j}(\omega)>0\}}
              \bbone_{\{s\leq \tau_m(\omega)\}}
 \end{align*}
 and
 \begin{align*}
 \RR_+\times U_d\times U_1\times \Omega\ni (s,\bz,u,\omega)
      \mapsto & \left( \phi_k(Y_{s-,i}(\omega) - z_i) - \phi_k(Y_{s-,i}(\omega)) \right)
                \bbone_{\{\|\bz\|\geq1\}}\\
             &\times \bbone_{\{X_{s-,j}(\omega)<u\leq X_{s-,j}'(\omega)\}} \bbone_{\{Y_{s-,j}(\omega)<0\}}
              \bbone_{\{s\leq \tau_m(\omega)\}}
 \end{align*}
 belong to the class \ $\bF_p^1$, \ and then \ $\left(I_{i,m,k,6,1}(t)\right)_{t\in\RR_+}$ \ is a martingale,
 again by Ikeda and Watanabe \cite[page 62]{IkeWat}.
By \eqref{help_Li_page45_equiv} and \eqref{(ii)},
 \begin{multline*}
  \EE\left( \int_0^{t\land\tau_m} \int_{U_d} \int_{U_1}
             \left| \phi_k(Y_{s-,i} + z_i) - \phi_k(Y_{s-,i}) \right|
             \bbone_{\{\|\bz\|\geq1\}}
             \bbone_{\{X_{s-,j}'<u\leq X_{s-,j}\}} \bbone_{\{Y_{s-,j}>0\}}
             \, \dd s \, \mu_j(\dd \bz) \, \dd u \right) \\
  \begin{aligned}
   &\leq \EE\left( \int_0^{t\land\tau_m} \int_{U_d} \int_{U_1}
                    z_i \bbone_{\{\|\bz\|\geq1\}} \bbone_{\{X_{s-,j}'<u\leq X_{s-,j}\}}
                    \bbone_{\{Y_{s-,j}>0\}}
                    \, \dd s \, \mu_j(\dd \bz) \, \dd u \right) \\
   &= \EE\left( \int_0^{t\land\tau_m} \int_{U_d}
                 z_i \bbone_{\{\|\bz\|\geq1\}} Y_{s-,j} \bbone_{\{Y_{s-,j}>0\}}
                 \, \dd s \, \mu_j(\dd \bz) \right)
    \leq 2mt\int_{U_d} z_i \bbone_{\{\|\bz\|\geq1\}} \, \mu_j(\dd\bz)
    < \infty ,
  \end{aligned}
 \end{multline*}
 and the finiteness of the other expectation can be shown in the same way.

Using the assumption \ $\Bbeta \leq \Bbeta'$, \ the property that the matrix
 \ $\bD$ \ has non-negative off-diagonal entries and the properties (ii)
 and (iii), we obtain
 \begin{align*}
  I_{i,m,k,1}(t)
  &= \int_0^{t\land\tau_m}
      \phi_k'(Y_{s,i})
      \biggl( \beta_i - \beta_i' + \sum_{j=1}^d d_{i,j} Y_{s,j} \biggr)
      \, \dd s \\
  &\leq \int_0^{t\land\tau_m}
         \phi_k'(Y_{s,i})
         \biggl( d_{i,i} Y_{s,i}
                 + \sum_{j\in\{1,\ldots,d\}\setminus\{i\}} d_{i,j} Y_{s,j}^+ \biggr)
         \bbone_{\RR_+}(Y_{s,i}) \, \dd s \\
  &\leq \int_0^{t\land\tau_m}
         \biggl( |d_{i,i}| Y_{s,i}^+
                 + \sum_{j\in\{1,\ldots,d\}\setminus\{i\}} d_{i,j} Y_{s,j}^+ \biggr)
         \, \dd s
   = \sum_{j=1}^d |d_{i,j}| \int_0^{t\land\tau_m} Y_{s,j}^+ \, \dd s .
 \end{align*}
By (iv),
 \[
   I_{i,m,k,3}(t) \leq (t\land\tau_m) c_i \frac{2}{k} \leq \frac{2 c_i t}{k} .
 \]
Now we estimate
 \begin{align*}
  &I_{i,m,k,5}(t)
   = \sum_{j=1}^d
      \int_0^{t\land\tau_m} \int_{U_d}
        \Bigl[ \phi_k(Y_{s-,i} + z_i) - \phi_k(Y_{s-,i})
               - \phi_k'(Y_{s-,i}) z_i \Bigr]
        \bbone_{\{\|\bz\|<1\}} Y_{s-,j} \bbone_{\{Y_{s-,j}>0\}}
        \, \dd s \, \mu_j(\dd\bz) \\
  &\qquad
      + \sum_{j=1}^d
         \int_0^{t\land\tau_m} \int_{U_d}
          \Bigl[ \phi_k(Y_{s-,i} - z_i) - \phi_k(Y_{s-,i})
                 + \phi_k'(Y_{s-,i}) z_i \Bigr]
          \bbone_{\{\|\bz\|<1\}} (-Y_{s-,j}) \bbone_{\{Y_{s-,j}<0\}}
          \, \dd s \, \mu_j(\dd\bz) .
 \end{align*}
By \eqref{(ii)} and (iii), we obtain
 \[
   \int_0^{t\land\tau_m} \int_{U_d}
     \Bigl[ \phi_k(Y_{s-,i} - z_i) - \phi_k(Y_{s-,i}) + \phi_k'(Y_{s-,i}) z_i \Bigr]
     \bbone_{\{\|\bz\|<1\}} (-Y_{s-,j}) \bbone_{\{Y_{s-,j}<0\}}
     \, \dd s \, \mu_j(\dd\bz)
   \leq 0
 \]
 for all \ $i, j \in \{1, \ldots, d\}$.
\ By \eqref{help3},
 \ $\int_{U_d} z_i \bbone_{\{\|\bz\|<1\}} \, \mu_j(\dd z) < \infty$
 \ for all \ $i, j \in \{1, \ldots, d\}$ \ with \ $i \ne j$, \ hence using
 (iii), we obtain
 \begin{align*}
  I_{i,m,k,5}(t)
  &\leq \int_0^{t\land\tau_m} \int_{U_d}
         \Bigl[ \phi_k(Y_{s-,i} + z_i) - \phi_k(Y_{s-,i})
                - \phi_k'(Y_{s-,i}) z_i \Bigr]
         \bbone_{\{\|\bz\|<1\}} Y_{s-,i}^+ \, \dd s \, \mu_i(\dd\bz) \\
  &\quad
        + \sum_{j\in\{1,\ldots,d\}\setminus\{i\}}
           \int_0^{t\land\tau_m} \int_{U_d}
            \Bigl[ \phi_k(Y_{s-,i} + z_i) - \phi_k(Y_{s-,i}) \Bigr]
            \bbone_{\{\|\bz\|<1\}} Y_{s-,j}^+ \, \dd s \, \mu_j(\dd\bz) .
 \end{align*}
By \eqref{(ii)}, for \ $i \ne j$,
 \[
   \int_0^{t\land\tau_m} \int_{U_d}
    \Bigl[ \phi_k(Y_{s-,i} + z_i) - \phi_k(Y_{s-,i}) \Bigr]
    \bbone_{\{\|\bz\|<1\}} Y_{s-,j}^+ \, \dd s \, \mu_j(\dd\bz) \\
   \leq \int_0^{t\land\tau_m} Y_{s,j}^+ \, \dd s
        \int_{U_d} z_i \bbone_{\{\|\bz\|<1\}} \, \mu_j(\dd\bz) .
 \]
Applying (iv) with \ $y = 0$, \ we have \ $x \phi_k''(x) \leq 2/k$ \ for all
 \ $x \in \RR_+$ \ and \ $k \in \NN$.
 \ By Taylor's theorem, for all \ $y \in \RR_{++}$, \ $z \in \RR_+$ \ and
 \ $k \in \NN$, \ there exists some
 \ $\vartheta = \vartheta(y, z) \in [0,1]$ \ such that
 \[
   \phi_k(y + z) - \phi_k(y) - \phi_k'(y) z
   = \phi_k''(y + \vartheta z) \frac{z^2}{2}
   \leq \frac{2 z^2}{2k(y + \vartheta z)}
   \leq \frac{z^2}{k y} .
 \]
Hence, using \eqref{z2}, we obtain
 \begin{align*}
  &\int_0^{t\land\tau_m} \int_{U_d}
    \Bigl[ \phi_k(Y_{s-,i} + z_i) - \phi_k(Y_{s-,i}) - \phi_k'(Y_{s-,i}) z_i \Bigr]
    \bbone_{\{\|\bz\|<1\}} Y_{s-,i}^+ \, \dd s \, \mu_i(\dd\bz) \\
  &\qquad\qquad
   \leq \int_0^{t\land\tau_m} \int_{U_d}
         \frac{z_i^2}{k Y_{s-,i}} \bbone_{\{Y_{s-,i}>0\}}
         \bbone_{\{\|\bz\|<1\}} Y_{s-,i}^+ \, \dd s \, \mu_i(\dd\bz)
   \leq \frac{t}{k} \int_{U_d} z_i^2 \bbone_{\{\|\bz\|<1\}} \, \mu_i(\dd\bz) .
 \end{align*}
Using \eqref{seged}, one can easily check that
 \begin{align*}
  I_{i,m,k,6,2}(t)
  &= \sum_{j=1}^d
      \int_0^{t\land\tau_m} \int_{U_d}
        \Bigl[ \phi_k(Y_{s-,i} + z_i) - \phi_k(Y_{s-,i}) \Bigr]
        \bbone_{\{\|\bz\|\geq1\}} Y_{s-,j} \bbone_{\{Y_{s-,j}>0\}}
        \, \dd s \, \mu_j(\dd\bz) \\
  &\quad
      + \sum_{j=1}^d
         \int_0^{t\land\tau_m} \int_{U_d}
          \Bigl[ \phi_k(Y_{s-,i} - z_i) - \phi_k(Y_{s-,i}) \Bigr]
          \bbone_{\{\|\bz\|\geq1\}} (-Y_{s-,j}) \bbone_{\{Y_{s-,j}<0\}}
          \, \dd s \, \mu_j(\dd\bz) .
 \end{align*}
By \eqref{(ii)}, we obtain
 \[
   \int_0^{t\land\tau_m} \int_{U_d}
     \Bigl[ \phi_k(Y_{s-,i} - z_i) - \phi_k(Y_{s-,i}) \Bigr]
     \bbone_{\{\|\bz\|\geq1\}} (-Y_{s-,j}) \bbone_{\{Y_{s-,j}<0\}}
     \, \dd s \, \mu_j(\dd\bz)
   \leq 0
 \]
 for all \ $i, j \in \{1, \ldots, d\}$.
\ By \eqref{help_Li_page45_equiv},
 \ $\int_{U_d} z_i \bbone_{\{\|\bz\|\geq1\}} \, \mu_j(\dd z) < \infty$ \ for all
 \ $i, j \in \{1, \ldots, d\}$, \ thus applying \eqref{(ii)}, we obtain
 \begin{align*}
  I_{i,m,k,6,2}(t)
  &\leq \sum_{j=1}^d
         \int_0^{t\land\tau_m} \int_{U_d}
         \Bigl[ \phi_k(Y_{s-,i} + z_i) - \phi_k(Y_{s-,i}) \Bigr]
         \bbone_{\{\|\bz\|\geq1\}} Y_{s-,j}^+ \, \dd s \, \mu_j(\dd\bz) \\
  &\leq \sum_{j=1}^d
         \int_0^{t\land\tau_m} Y_{s,j}^+ \, \dd s
         \int_{U_d} z_i \bbone_{\{\|\bz\|\geq1\}} \, \mu_j(\dd\bz) .
 \end{align*}
Summarizing, we have
 \begin{align}\label{help_comparison}
  \begin{split}
  \phi_k(Y_{t\land\tau_m,i})
  &\leq \phi_k(Y_{0,i})
        + C_i \sum_{j=1}^d \int_0^{t\land\tau_m} Y_{s,j}^+ \, \dd s
        + \frac{2 c_i t}{k}
        + \frac{t}{k}
          \int_{U_d} z_i^2 \bbone_{\{\|\bz\|<1\}} \, \mu_i(\dd\bz) \\
  &\quad
        + I_{i,m,k,2}(t) + I_{i,m,k,4}(t) + I_{i,m,k,6,1}(t) , \qquad t \in \RR_+ ,
  \end{split}
 \end{align}
 where
 \begin{align*}
  C_i := \max_{j\in\{1,\ldots,d\}} |d_{i,j}|
         + \max_{j\in\{1,\ldots,d\}\setminus\{i\}} \int_{U_d} z_i  \, \mu_j(\dd\bz)
         + \int_{U_d} z_i \bbone_{\{\|\bz\|\geq1\}} \, \mu_i(\dd\bz) .
 \end{align*}
By (iii), we obtain \ $\PP(\phi_k(Y_{0,i}) \leq 0) = 1$,
 \ $i \in \{1, \ldots, d\}$.
\ By (i), the non-negativeness of \ $\phi_k$ \ and monotone convergence
 theorem yield \ $\EE(\phi_k(Y_{t\land\tau_m,i})) \to \EE(Y_{t\land\tau_m,i}^+)$ \ as
 \ $k \to \infty$ \ for all \ $t \in \RR_+$, \ $m\in\NN$, \ and
 \ $i \in \{1, \ldots, d\}$.
\ We have
 \ $\int_0^{t\land\tau_m} Y_{s,j}^+ \, \dd s
    \leq \int_0^t Y_{s\land\tau_m,j}^+ \, \dd s$,
 \ hence taking the expectation of \eqref{help_comparison} and letting
 \ $k \to \infty$, \ we obtain
 \[
   \EE\biggl(\sum_{i=1}^d Y_{t\land\tau_m,i}^+\biggr)
   \leq C \int_0^t \EE\biggl(\sum_{i=1}^d Y_{s\land\tau_m,i}^+\biggr) \dd s ,
 \]
 with \ $C := \sum_{i=1}^d C_i$.
\ By Gronwall's inequality, we conclude
 \[
   \EE\biggl(\sum_{i=1}^d Y_{t\land\tau_m,i}^+\biggr) = 0
 \]
 for all \ $t \in \RR_+$ \ and \ $m \in \NN$.
\ Hence \ $\PP(\bX_{t\land\tau_m,i} \leq \bX'_{t\land\tau_m,i} )=1$ \ for all
 \ $t \in \RR_+$, \ $m \in \NN$ \ and \ $i \in \{1, \ldots, d\}$, \ and then
 \ $\PP(\text{$\bX_{t\land\tau_m,i} \leq \bX'_{t\land\tau_m,i}$ for all
              $m \in \NN$})
    = 1$
 \ for all \ $t \in \RR_+$ \  and \ $i \in \{1, \ldots, d\}$.
\ Since \ $\bX$ \ and \ $\bX'$ \ have c\`adl\`ag trajectories, these
 trajectories are bounded almost surely on \ $[0, T]$ \ for all
 \ $T \in \RR_+$, \ hence \ $\tau_m \as \infty$ \ as \ $m \to \infty$.
\ This yields \ $\PP(\bX_t \leq \bX'_t )=1$ \ for all \ $t \in \RR_+$.
\ Since the set of non-negative rational numbers \ $\QQ_+$ \ is countable, we
 obtain \ $\PP(\text{$\bX_t \leq \bX'_t$ for all $t \in \QQ_+$}) = 1$.
\ Using again that \ $\bX$ \ and \ $\bX'$ \ have c\`adl\`ag trajectories
 almost surely, we get
 \ $\PP(\text{$\bX_t \leq \bX_t'$ for all $t \in \RR_+$}) = 1$.
\proofend

\begin{Rem}
We note that Dawson and Li \cite[Theorem 2.3]{DawLi2} provided a comparison
 theorem for SDEs with jumps in a much more general setting, but only for
 1-dimensional processes.
\proofend
\end{Rem}

Consider the following objects:
 \begin{enumerate}
  \item[(E1)]
   a probability space \ $(\Omega, \cF, \PP)$;
  \item[(E2)]
   a \ $d$-dimensional standard Brownian motion \ $(\bW_t)_{t\in\RR_+}$;
  \item[(E3)]
   a stationary Poisson point process \ $p$ \ on \ $V$ \ with characteristic
    measure \ $m$ \ given in \eqref{m};
  \item[(E4)]
   a random vector \ $\bxi$ \ with values in \ $\RR_+^d$, \ independent of
    \ $\bW$ \ and \ $p$.
 \end{enumerate}

\begin{Rem}\label{dRM_strong}
Note that if conditions (E1)--(E4) are satisfied, then \ $\bxi$, \ $\bW$ \ and
 \ $p$ \ are automatically mutually independent according to Remark
 \ref{dRM}.
\proofend
\end{Rem}

Provided that the objects (E1)--(E4) are given, let
 \ $(\cF^{\bxi,\bW\!,\,p}_t)_{t\in\RR_+}$ \ be the augmented filtration generated by
 \ $\bxi$, \ $\bW$ \ and \ $p$, \ i.e., for each \ $t \in \RR_+$,
 \ $\cF^{\bxi,\bW\!,\,p}_t$ \ is the $\sigma$-field generated by
 \ $\sigma(\bxi; \, \bW_s, s\in[0,t]; \, p(s), s\in(0,t]\cap D(p))$ \ and by
 the \ $\PP$-null sets from
 \ $\sigma(\bxi; \, \bW_s, s\in\RR_+; \, p(s), s\in\RR_{++}\cap D(p))$ \ (which
 is similar to the definition in Karatzas and Shreve \cite[page 285]{KarShr}).
One can check that \ $(\cF^{\bxi,\bW\!,\,p}_t)_{t\in\RR_+}$ \ satisfies the usual
 hypotheses, \ $(\bW_t)_{t\in\RR_+}$ \ is a standard
 $(\cF^{\bxi,\bW\!,\,p}_t)_{t\in\RR_+}$-Brownian motion, and \ $p$ \ is a stationary
 \ $(\cF^{\bxi,\bW\!,\,p}_t)_{t\in\RR_+}$-Poisson point process on \ $V$
 \ with characteristic measure \ $m$, \ see, e.g., Barczy et al.\
 \cite{BarLiPap}.

\begin{Def}\label{Def_strong_solution2}
Suppose that the objects \textup{(E1)--(E4)} are given.
An \ $\RR_+^d$-valued strong solution of the SDE \eqref{SDE_X_YW_spec} on
 \ $(\Omega, \cF, \PP)$ \ and with respect to the standard Brownian motion
 \ $\bW$, \ the stationary Poisson point process \ $p$ \ and initial value
 \ $\bxi$, \ is an \ $\RR_+^d$-valued
 \ $(\cF^{\bxi,\bW\!,\,p}_t)_{t\in\RR_+}$-adapted c\`{a}dl\`{a}g process
 \ $(\bX_t)_{t \in \RR_+}$ \ with \ $\PP(\bX_0 = \bxi) = 1$ \ satisfying
 \textup{(D4)(b)--(e)}.
\end{Def}

Clearly, if \ $(\bX_t)_{t \in \RR_+}$ \ is an \ $\RR_+^d$-valued strong solution,
 then
 \ $\bigl( \Omega, \cF, (\cF^{\bxi,\bW\!,\,p}_t)_{t\in\RR_+}, \PP, \bW, p,
           \bX \bigr)$
 \ is an \ $\RR_+^d$-valued weak solution.

\begin{Thm}\label{strong_solution_old}
Let \ $(d, \bc, \Bbeta, \bB, \nu, \bmu)$ \ be a set of admissible parameters
 in the sense of Definition \ref{Def_admissible} such that the moment
 condition \eqref{moment_condition_1} holds.
Suppose that objects \textup{(E1)--(E4)} are given.
If \ $\EE(\|\bxi\|) < \infty$, \ then there is a pathwise unique
 \ $\RR_+^d$-valued strong solution to the SDE \eqref{SDE_X_YW_spec} with
 initial value \ $\bxi$, \ and the solution is a CBI process with
 parameters \ $(d, \bc, \Bbeta, \bB, \nu, \bmu)$.
\end{Thm}

\noindent
\textbf{Proof.}
The pathwise uniqueness among \ $\RR_+^d$-valued weak solutions follows from
 Lemma \ref{comparison_beta}.
Then, by Theorem 5.5 in Barczy et al. \cite{BarLiPap} (Yamada-Watanabe type
 result for SDEs with jumps) and Theorem \ref{CBI_SDE}, we conclude that the
 SDE \eqref{SDE_X_YW_spec} has a pathwise
 unique \ $\RR_+^d$-valued strong solution.
\proofend

\section{Special cases}
\label{section_special_cases}

In this section we specialize our results to dimension 1 and 2.
Moreover, we consider a special case of the SDE \eqref{SDE_X_YW_spec} with
 \ $\nu = 0$, \ $\mu_i = 0$, \ $i \in \{1, \ldots, d\}$, \ i.e., without
 integrals with respect to (compensated) Poisson random measures, and another
 special case with \ $\bc = \bzero$, \ i.e., without integral with respect to
 a Wiener process.

First we rewrite the SDE \eqref{SDE_X_YW_spec} in a form which is more
 comparable with the results of Li \cite[Theorem 9.31]{Li} (one-dimensional
 case) and Ma \cite[Theorem 3.2]{Ma} (two-dimensional case).

For each \ $j \in \{0, 1, \ldots, d\}$, \ the thinning \ $p_j$ \ of \ $p$
 \ onto \ $\cR_j$ \ is again a stationary \ $(\cF_t)_{t\in\RR_+}$-Poisson point
 process on \ $\cR_j$, \ and its characteristic measure is the restriction
 \ $m|_{\cR_j}$ \ of \ $m$ \ onto \ $\cR_j$ \ (this can be checked calculating
 its conditional Laplace transform, see Ikeda and Watanabe
 \cite[page 44]{IkeWat}).
Using these Poisson point processes, we obtain the useful decomposition
 \begin{equation}\label{deco}
  \begin{aligned}
   &\int_0^t \int_{V_0} f(\bX_{s-}, \br) \, \tN(\dd s, \dd\br)
    + \int_0^t \int_{V_1} g(\bX_{s-}, \br) \, N(\dd s, \dd\br) \\
   &= \sum_{j=1}^d
       \int_0^t \int_{\cR_{j,0}}
        \bz \bbone_{\{u\leq X_{s-,j}\}}
        \, \tN_j(\dd s, \dd\br) \\
   &\quad
      + \sum_{j=1}^d
         \int_0^t \int_{\cR_{j,1}}
          \bz \bbone_{\{u\leq X_{s-,j}\}}
          \, N_j(\dd s, \dd\br)
      + \int_0^t \int_{\cR_0}
         \br \, M(\dd s, \dd\br) ,
  \end{aligned}
 \end{equation}
 where, for each \ $j \in \{1, \ldots, d\}$, \ $N_j(\dd s, \dd\br)$ \ is the
 counting measure of \ $p_j$ \ on \ $\RR_{++} \times \cR_j$,
 \ $\tN_j(\dd s, \dd\br)
    := N_j(\dd s, \dd\br) - \dd s \, (\mu_j(\dd \bz) \, \dd u)$,
 \ and \ $M(\dd s, \dd\br)$ \ is the counting measure of \ $p_0$ \ on
 \ $\RR_{++} \times \cR_0$.
\ Indeed,
 \begin{equation*}
  \begin{aligned}
  \int_0^t \int_{\cR'} F(s, \br) \, \tN(\dd s, \dd\br)
  &= \int_0^t \int_{\cR'} F(s, \br) \, \tN'(\dd s, \dd\br) ,
  \qquad F \in \bF_p^{2,loc} , \\
  \int_0^t \int_{\cR'} G(s, \br) \, N(\dd s, \dd\br)
  &= \int_0^t \int_{\cR'} G(s, \br) \, N'(\dd s, \dd\br) , \qquad G \in \bF_p ,
  \end{aligned}
 \end{equation*}
 are valid for the thinning \ $p'$ \ of \ $p$ \ onto any measurable subset
 \ $\cR' \subset \cR$, \ where \ $N'(\dd s, \dd\br)$ \ denotes the counting
 measure of the stationary \ $(\cF_t)_{t\in\RR_+}$-Poisson point process
 \ $p'$, \ and
 \ $\tN'(\dd s, \dd\br)
    := N'(\dd s, \dd\br) - \bbone_{\{\br\in\cR'\}} \dd s \, m(\dd\br)$.

Remark that for any \ $\RR_+^d$-valued weak solution of the SDE
 \eqref{SDE_X_YW_spec}, the Brownian motion \ $\bW$ \ and the stationary Poisson
 point processes \ $p_j$, \ $j \in \{0, 1, \ldots, d\}$ \ are mutually
 independent according again to Theorem 6.3 in Chapter II of Ikeda and
 Watanabe \cite{IkeWat}.
Indeed, the intensity measures of \ $p_j$, \ $j \in \{0, 1, \ldots, d\}$, \ are
 deterministic, and condition (6.11) of this theorem is satisfied, because
 \ $p_j$, \ $j \in \{0, 1, \ldots, d\}$, \ live on disjoint subsets of \ $\cR$.

For \ $d = 1$, \ applying \eqref{deco}, the SDE \eqref{SDE_X_YW_spec} takes
 the form
 \begin{align*}
  X_t
  &= X_0 + \int_0^t (\beta + d X_s) \, \dd s
     + \int_0^t \sqrt{2 c X_s^+} \, \dd W_s \\
   &\quad
      + \int_0^t \int_{\cR_{1,0}} z \bbone_{\{u\leq X_{s-}\}} \, \tN_1(\dd s, \dd r)
      + \int_0^t \int_{\cR_{1,1}} z \bbone_{\{u\leq X_{s-}\}} \, N_1(\dd s, \dd r)
      + \int_0^t \int_{\cR_0} r \, M(\dd s, \dd r)
 \end{align*}
 for \ $t\in\RR_+$, \ where \ $\beta \in \RR_+$,
 \ $d = \tb - \int_0^\infty z \bbone_{\{z\geq1\}} \, \mu_1(\dd z)$,
 \ $\tb = b + \int_0^\infty (z - 1)^+ \, \mu_1(\dd z)$,
 \ $b \in \RR$, \ $c \in \RR_+$,
 \ $\cR_{1,0} = \{0\} \times \{ z \in \RR_{++} : z < 1\} \times \RR_{++}$,
 \ $\cR_{1,1} = \{0\} \times \{ z \in \RR_{++} : z \geq 1\} \times \RR_{++}$,
 \ $\cR_0 = \RR_{++} \times \{(0, 0)\}$.
\ We have
 \begin{align*}
  I_0
  &:=\int_0^t \int_{\cR_{1,0}} z \bbone_{\{u\leq X_{s-}\}} \, \tN_1(\dd s, \dd r)
    = \int_0^t \int_0^\infty \int_0^\infty
       z \bbone_{\{z<1\}} \bbone_{\{u\leq X_{s-}\}}
       \, \toN_1(\dd s, \dd z, \dd u) , \\
  I_1
  &:= \int_0^t \int_{\cR_{1,1}} z \bbone_{\{u\leq X_{s-}\}} \, N_1(\dd s, \dd r)
    = \int_0^t \int_0^\infty \int_0^\infty
       z \bbone_{\{z\geq1\}} \bbone_{\{u\leq X_{s-}\}}
       \, \oN_1(\dd s, \dd z, \dd u) , \\
  I_2
  &:= \int_0^t \int_{\cR_0} r \, M(\dd s, \dd r)
    = \int_0^t \int_0^\infty z \, \oM(\dd s, \dd z) ,
 \end{align*}
 where \ $\oN_1$ \ and \ $\oM$ \ are Poisson random measures on
 \ $\RR_{++} \times \RR_{++}^2$ \ and on \ $\RR_{++} \times \RR_{++}$ \ with
 intensity measures \ $\dd s \, \mu_1(\dd z) \, \dd u$ \ and
 \ $\dd s \, \nu(\dd z)$, \ respectively, and
 \ $\toN_1(\dd s, \dd z, \dd u)
    := \oN_1(\dd s, \dd z, \dd u) - \dd s \, \mu_1(\dd z) \, \dd u$.
\ Under the moment conditions \eqref{help_Li_page45_equiv},
 \[
   I_0 + I_1
   = \int_0^t \int_0^\infty \int_0^\infty
      z \bbone_{\{u\leq X_{s-}\}} \, \toN_1(\dd s, \dd z, \dd u)
     + \int_0^t X_s \, \dd s
       \int_0^\infty z \bbone_{\{z\geq1\}} \, \mu_1(\dd z) .
 \]
Consequently, the SDE \eqref{SDE_X_YW_spec} can be rewritten in the form
 \begin{align*}
  X_t
  &= X_0 + \int_0^t (\beta + \tb X_s) \, \dd s
     + \int_0^t \sqrt{2 c X_s^+} \, \dd W_s \\
   &\quad
      + \int_0^t \int_0^\infty \int_0^\infty
         z \bbone_{\{u\leq X_{s-}\}} \, \toN_1(\dd s, \dd z, \dd u)
      + \int_0^t \int_0^\infty z \, \oM(\dd s, \dd z) ,
      \qquad t\in\RR_+,
 \end{align*}
 hence, taking into account the form \eqref{CBI_inf_gen_1} of the
 infinitesimal generator of the process \ $(X_t)_{t\in\RR_+}$, \ we obtain
 equation (9.46) of Li \cite{Li}.

In a similar way, for \ $d = 2$, \ applying \eqref{deco}, the SDE
 \eqref{SDE_X_YW_spec} takes the form
 \begin{align*}
  \bX_t
  &= \bX_0 + \int_0^t (\Bbeta + \bD \bX_s) \, \dd s
     + \sum_{i=1}^2
        \int_0^t \sqrt{2 c_i X_{s,i}^+} \be_i \be_i^\top \, \dd \bW_s
     + \int_0^t \int_{\cR_0} \br \, M(\dd s, \dd \br) \\
   &\quad
     + \sum_{j=1}^2
        \int_0^t \int_{\cR_{j,0}}
         \bz \bbone_{\{u\leq X_{s-,j}\}} \, \tN_j(\dd s, \dd \br)
     + \sum_{j=1}^2
        \int_0^t \int_{\cR_{j,1}}
         \bz \bbone_{\{u\leq X_{s-,j}\}} \, N_j(\dd s, \dd \br)
 \end{align*}
 for \ $t\in\RR_+$, \ where \ $\Bbeta \in \RR_+^2$, \ $\bD$ \ is given in
 \eqref{tbeta_tB}, \ $(c_1, c_2)^\top \in \RR_+^2$,
 \begin{align*}
  \cR_0 &= U_2 \times \{(0, 0, 0)\} \times \{(0, 0, 0)\} , \\
  \cR_{1,0} &= \{(0, 0)\} \times \{ \bz \in U_2 : \|\bz\| < 1 \}
              \times \RR_{++} \times \{(0, 0, 0)\} , \\
  \cR_{2,0} &= \{(0, 0)\} \times \{(0, 0, 0)\}
              \times \{ \bz \in U_2 : \|\bz\| < 1 \} \times \RR_{++} , \\
  \cR_{1,1} &= \{(0, 0)\} \times \{ \bz \in U_2 : \|\bz\| \geq 1 \}
              \times \RR_{++} \times \{(0, 0, 0)\} , \\
  \cR_{2,1} &= \{(0, 0)\} \times \{(0, 0, 0)\}
              \times \{ \bz \in U_2 : \|\bz\| \geq 1 \} \times \RR_{++} .
 \end{align*}
For each \ $j \in \{1, 2\}$, \ we have
 \begin{align*}
  I_{j,0}
  &:=\int_0^t \int_{\cR_{j,0}}
      \bz \bbone_{\{u\leq X_{s-,j}\}} \, \tN_j(\dd s, \dd \br)
    = \int_0^t \int_{U_2} \int_0^\infty
       \bz \bbone_{\{\|\bz\|<1\}} \bbone_{\{u\leq X_{s-,j}\}}
       \, \toN_j(\dd s, \dd \bz, \dd u) , \\
  I_{j,1}
  &:= \int_0^t \int_{\cR_{j,1}}
       \bz \bbone_{\{u\leq X_{s-,j}\}} \, N_j(\dd s, \dd \br)
    = \int_0^t \int_{U_2} \int_0^\infty
       \bz \bbone_{\{\|\bz\|\geq1\}} \bbone_{\{u\leq X_{s-,j}\}}
       \, \oN_j(\dd s, \dd \bz, \dd u) , \\
  I_2
  &:= \int_0^t \int_{\cR_0} \br \, M(\dd s, \dd \br)
    = \int_0^t \int_{U_2} \bz \, \oM(\dd s, \dd \bz) ,
 \end{align*}
 where \ $\oN_j$ \ and \ $\oM$ \ are Poisson random measures on
 \ $\RR_{++} \times U_2 \times \RR_{++}$ \ and on \ $\RR_{++} \times U_2$
 \ with intensity measures \ $\dd s \, \mu_j(\dd\bz) \, \dd u$ \ and
 \ $\dd s \, \nu(\dd\bz)$, \ respectively, and
 \ $\toN_j(\dd s, \dd \bz, \dd u)
    := \oN_j(\dd s, \dd \bz, \dd u) - \dd s \, \mu_j(\dd \bz) \, \dd u$.
\ Under the moment conditions \eqref{help_Li_page45_equiv},
 \[
   I_{j,0} + I_{j,1}
   = \int_0^t \int_{U_2} \int_0^\infty
      \bz \bbone_{\{u\leq X_{s-,j}\}} \, \toN_j(\dd s, \dd \bz, \dd u)
     + \int_0^t X_{s,j} \, \dd s
       \int_{U_2} \bz \bbone_{\{\|\bz\|\geq1\}} \, \mu_j(\dd \bz) .
 \]
Consequently, the SDE \eqref{SDE_X_YW_spec} can be rewritten in the form
 \begin{align*}
  \bX_t
  &= \bX_0 + \int_0^t (\Bbeta + \tbB \bX_s) \, \dd s
     + \sum_{i=1}^2
        \int_0^t \sqrt{2 c_i X_{s,i}^+} \, \dd W_{s,i} \, \be_i \\
   &\quad
      + \sum_{j=1}^2
         \int_0^t \int_{U_2} \int_0^\infty
          \bz \bbone_{\{u\leq X_{s-,j}\}} \, \toN_j(\dd s, \dd \bz, \dd u)
      + \int_0^t \int_{U_2} \bz \, \oM(\dd s, \dd \bz) , \qquad t\in\RR_+.
 \end{align*}
Due to \eqref{help3}, we have
 \begin{align*}
  X_{t,1} & = X_{0,1}
             + \int_0^t \left(\beta_1 + \tb_{1,1} X_{s,1}
                            + \left(\tb_{1,2} - \int_{U_2} z_1\, \mu_2(\dd \bz) \right) X_{s,2}\right)\,\dd s
             + \int_0^t\sqrt{2c_1 X_{s,1}^+}\,\dd W_{s,1}\\
          &\quad  + \int_0^t \int_{U_2} \int_0^\infty
              z_1 \bbone_{\{u\leq X_{s-,1}\}} \, \toN_1(\dd s, \dd \bz, \dd u)
            + \int_0^t \int_{U_2} \int_0^\infty
              z_1 \bbone_{\{u\leq X_{s-,2}\}} \, \oN_2(\dd s, \dd \bz, \dd u) \\
           &\quad + \int_0^t \int_{U_2} z_1 \, \oM(\dd s, \dd \bz) , \qquad t\in\RR_+,
 \end{align*}
 and
 \begin{align*}
  X_{t,2} &= X_{0,2}
            + \int_0^t \left(\beta_2 + \left(\tb_{2,1} - \int_{U_2} z_2 \,\mu_1(\dd \bz) \right) X_{s,1}
                             +  \tb_{2,2} X_{s,2}\right)\,\dd s
            + \int_0^t\sqrt{2c_2 X_{s,2}^+}\,\dd W_{s,2}\\
        &\quad + \int_0^t \int_{U_2} \int_0^\infty
              z_2 \bbone_{\{u\leq X_{s-,2}\}} \, \toN_2(\dd s, \dd \bz, \dd u)
            + \int_0^t \int_{U_2} \int_0^\infty
              z_2 \bbone_{\{u\leq X_{s-,1}\}} \, \oN_1(\dd s, \dd \bz, \dd u)\\
        &\quad   + \int_0^t \int_{U_2} z_2 \, \oM(\dd s, \dd \bz) , \qquad t\in\RR_+.
 \end{align*}
In the special case \ $\nu = 0$, \ we obtain equations (2.1) and (2.2) of Ma \cite{Ma}.
Indeed, due to \eqref{help3}, one can rewrite the infinitesimal generator \eqref{CBI_inf_gen_1}
 of the process \ $(\bX_t)_{t\in\RR_+}$ \ in the following form
 \begin{align*}
   (\cA_{\bX} f)(\bx)
   &= \sum_{i=1}^2 c_i x_i f_{i,i}''(x)
      + \sum_{i=1}^2
         x_i
         \int_{U_2}
          \bigl(f(\bx + \bz) - f(\bx) - z_if_i'(x)\bigr)
          \, \mu_i(\dd \bz) \\
   &\quad
      + \langle \Bbeta + \tbB \bx, \Bf'(\bx) \rangle
      + \int_{U_2} \bigl(f(\bx + \bz) - f(\bx)\bigr) \, \nu(\dd \bz)\\
   &\quad - x_1f_2'(x)\int_{U_2}z_2\,\mu_1(\dd z)
          - x_2f_1'(x)\int_{U_2}z_1\,\mu_2(\dd z) \\
  &= \sum_{i=1}^2 c_i x_i f_{i,i}''(x)
      + \sum_{i=1}^2
         x_i
         \int_{U_2}
          \bigl(f(\bx + \bz) - f(\bx) - z_if_i'(x)\bigr)
          \, \mu_i(\dd \bz) \\
   &\quad
      + \langle \Bbeta + \widetilde\tbB \bx, \Bf'(\bx) \rangle
      + \int_{U_2} \bigl(f(\bx + \bz) - f(\bx)\bigr) \, \nu(\dd \bz)
 \end{align*}
 for \ $f \in C^2_\cc(\RR_+^d,\RR)$ \ and \ $\bx \in \RR_+^d$, \ where
 \[
   \widetilde\tbB
   := \begin{bmatrix}
       \tb_{1,1} & \tb_{1,2} - \int_{U_2}z_1\,\mu_2(\dd \bz) \\[1mm]
       \tb_{2,1} - \int_{U_2}z_2\,\mu_1(\dd \bz) & \tb_{2,2}
      \end{bmatrix}.
 \]
This form of the infinitesimal generator \ $\cA_{\bX}$ \ is readily comparable with
 the corresponding one in Ma \cite[equation (1.5)]{Ma}.

In what follows, we consider a special form of the SDE \eqref{SDE_X_YW_spec}
 without integrals with respect to (compensated) Poisson random measures.
Namely, if \ $\nu = 0$, \ $\mu_i = 0$, \ $i \in \{1, \ldots, d\}$, \ then
 the SDE \eqref{SDE_X_YW_spec} takes the form
 \begin{align*}
   \bX_t
   &= \bX_0 + \int_0^t b(\bX_s) \, \dd s
      + \int_0^t \sigma(\bX_s) \, \dd \bW_s \\
   & = \bX_0 + \int_0^t (\Bbeta + \bB \bX_s)\,\dd s
       + \sum_{i=1}^d \int_0^t \sqrt{2 c_i X_{s,i}} \be_i \be_i^\top
       \,\dd \bW_s,  \qquad
   t \in \RR_+ ,
 \end{align*}
 and consequently,
 \begin{align*}
   X_{t,i} = \int_0^t \Bigg(\beta_i + \sum_{j=1}^d b_{i,j} X_{s,j}\Bigg)\dd t
             + \int_0^t \sqrt{2c_i X_{s,i}}\,\dd W_{s,i},
   \qquad t\in\RR_+, \quad i\in\{1,\ldots,d\}.
 \end{align*}
If \ $\bB$ \ is diagonal, then the process \ $(\bX_t)_{t\in\RR_+}$ \ is known to
 be a multi-factor Cox-Ingersoll-Ross process, see, e.g.,
 Jagannathan et al. \cite{JagKapSun}.

Finally, Theorem \ref{strong_solution_old} is valid also if the SDE
 \eqref{SDE_X_YW_spec} does not contain integral with respect to a Wiener
 process, i.e., if \ $\bc = \bzero$.
\ We note that in the proof of Theorem \ref{CBI_SDE} we applied Theorem 7.1'
 in Chapter II of Ikeda and Watanabe \cite{IkeWat}, which is valid in case
 \ $\bc = \bzero$ \ as well.

\vspace*{5mm}

\appendix

\vspace*{5mm}

\noindent{\bf\Large Appendix}

\section{Extension of a probability space}
\label{section_ext}

We recall the definition of extensions of probability spaces, see, e.g.,
 Ikeda and Watanabe \cite[Chapter II, Definition 7.1]{IkeWat}.

\begin{Def}\label{Def_extension}
We say that a filtered probability space
 \ $(\tOmega, \tcF, (\tcF_t)_{t\in\RR_+}, \tPP)$ \ is an extension of a filtered
 probability space \ $(\Omega, \cF, (\cF_t)_{t\in\RR_+}, \PP)$, \ if there exists
 an \ $\tcF/\cF$-measurable mapping \ $\pi :\tOmega \to \Omega$ \ such that
 \ $\pi^{-1}(\cF_t) \subset \tcF_t$ \ for all \ $t \in \RR_+$,
 \ $\PP(A) = \tPP(\pi^{-1}(A))$ \ for all \ $A \in \cF$, \ and
 \ $\tEE(\tX \mid \tcF_t)(\tomega) = \EE(X \mid \cF_t)(\pi(\tomega))$
 \ $\tPP$-almost surely for each essentially bounded
 \ ($\cF/\cB(\RR^d)$-measurable) random variable \ $X : \Omega \to \RR^d$,
 \ where we set \ $\tX(\tomega) := X(\pi(\tomega))$, \ $\tomega \in \tOmega$.
\end{Def}

\begin{Rem}\label{Rem_extension}
With the notations of Definition \ref{Def_extension},
 if \ $(\bX_t)_{t\in\RR_+}$ \ is an \ $\RR^d$-valued \ $(\cF_t)_{t\in\RR_+}$-adapted stochastic process,
 then \ $(\tbX_t)_{t\in\RR_+}$ \ is \ $(\tcF_t)_{t\in\RR_+}$-adapted.
Indeed, for each \ $t \in \RR_+$ \ and \ $B \in \cB(\RR^d)$, \ we have
 \begin{align*}
  \tbX_t^{-1}(B)
  = \{ \tomega \in \tOmega : \tbX_t(\tomega) \in B\}
  = \{ \tomega \in \tOmega : \bX_t(\pi(\tomega)) \in B\}
  = \pi^{-1}(\bX_t^{-1}(B))
  \in \tcF_t ,
 \end{align*}
 since \ $\bX_t^{-1}(B) \in \cF_t$.
\proofend
 \end{Rem}

\begin{Lem}\label{ext_BM}
Let \ $(\Omega, \cF, (\cF_t)_{t\in\RR_+}, \PP)$ \ be a filtered probability
 space, and let \ $(\bW_t)_{t\in\RR_+}$ \ be a $d$-dimensional
 \ $(\cF_t)_{t\in\RR_+}$-Brownian motion.
Let \ $(\tOmega, \tcF, (\tcF_t)_{t\in\RR_+}, \tPP)$ \ be an extension of
 \ $(\Omega, \cF, (\cF_t)_{t\in\RR_+}, \PP)$ \ with the mapping
 \ $\pi :\tOmega \to \Omega$.
\ Let \ $\tbW_t(\tomega) := \bW_t(\pi(\tomega))$ \ for all
 \ $\tomega \in \tOmega$ \ and \ $t \in \RR_+$.
\ Then \ $(\tbW_t)_{t\in\RR_+}$ \ is a $d$-dimensional
 \ $(\tcF_t)_{t\in\RR_+}$-Brownian motion.
\end{Lem}

\noindent
\textbf{Proof.}
According to Ikeda and Watanabe \cite[Chapter I, Definition 7.2]{IkeWat}, we
 have to check that the process \ $(\tbW_t)_{t\in\RR_+}$ \ has continuous
 trajectories, it is \ $(\tcF_t)_{t\in\RR_+}$-adapted, and satisfies
 \[
   \tEE(\exp\{\ii \langle \bu, \tbW_t - \tbW_s\rangle\} \mid \tcF_s)
   = \ee^{-(t-s)\|\bu\|^2/2} \qquad \text{$\tPP$-almost surely}
 \]
 for every \ $\bu \in \RR^d$ \ and \ $s, t \in \RR_+$ \ with \ $s < t$.
\ Clearly, \ $\RR_+ \ni t \mapsto \tbW_t(\tomega) = \bW_t(\pi(\tomega))$ \ is
 continuous for all \ $\tomega \in \tOmega$.
\ By Remark \ref{Rem_extension}, \ $(\tbW_t)_{t\in\RR_+}$ \ is
 \ $(\tcF_t)_{t\in\RR_+}$-adapted.
Finally, for every \ $\bu \in \RR^d$ \ and \ $s, t \in \RR_+$ \ with \ $s < t$,
 \[
   \tEE(\exp\{\ii \langle \bu, \tbW_t - \tbW_s\rangle\} \mid \tcF_s)(\tomega)
   = \EE(\exp\{\ii \langle \bu, \bW_t - \bW_s\rangle\}
         \mid \cF_s)(\pi(\tomega))
   = \ee^{-(t-s)\|\bu\|^2/2}
 \]
$\tPP$-almost surely, since we have \ $\xi(\omega) = c$ \ $\PP$-almost surely
 with \ $\xi := \EE(\exp\{\ii \langle \bu, \bW_t - \bW_s\rangle\} \mid \cF_s)$
 \ and $c := \ee^{-(t-s)\|\bu\|^2/2}$, \ which implies
 \ $\xi(\pi(\tomega)) = c$ \ $\tPP$-almost surely, because
 \ $\tPP(\{ \tomega \in \tOmega : \xi(\pi(\tomega)) = c \})
    = \tPP(\pi^{-1}(\xi^{-1}(\{c\}))) = \PP(\xi^{-1}(\{c\})) = 1$.
\proofend

\begin{Lem}\label{ext_BM_PRN}
Let \ $(\Omega, \cF, (\cF_t)_{t\in\RR_+}, \PP)$ \ be a filtered probability
 space, let \ $(\bW_t)_{t\in\RR_+}$ \ be a $d$-dimensional
 \ $(\cF_t)_{t\in\RR_+}$-Brownian motion, and let \ $p$ \ be a stationary
 \ $(\cF_t)_{t\in\RR_+}$-Poisson point process on
 \ $V = \RR_+^d \times (\RR_+^d \times \RR_+)^d$ \ with characteristic measure
 \ $m$, \ where \ $m$ \ is given in \eqref{m}.
Let
 \[
   \cG_t := \bigcap_{\vare>0} \sigma\left(\cF_{t+\vare} \cup \cN\right) ,
   \qquad t \in \RR_+ ,
 \]
 where \ $\cN$ \ denotes the collection of null sets under the probability
 measure \ $\PP$.
\ Then \ $(\bW_t)_{t\in\RR_+}$ \ is a $d$-dimensional
 \ $(\cG_t)_{t\in\RR_+}$-Brownian motion, and \ $p$ \ is a stationary
 \ $(\cG_t)_{t\in\RR_+}$-Poisson point process \ on \ $V$ \ with
 characteristic measure \ $m$.
\end{Lem}

\noindent
\textbf{Proof.}
The proof is essentially the same as the proof of Lemma A.5 in
 Barczy et al.~\cite{BarLiPap}.
\proofend

\section*{Acknowledgements}
We would like to thank the referee for his/her comments that helped us to improve the presentation
 of the Introduction.

\end{document}